\documentclass[11pt]{article}
\usepackage{graphicx}

\usepackage{subfiles}

\usepackage{algorithm}
\usepackage{algorithmicx}
\usepackage{algpseudocode}
\usepackage{latexsym,amsmath,cases,amsfonts,amscd, amsthm, dsfont}
\usepackage{extarrows}
\usepackage{bm,color}
\usepackage{epsfig,verbatim,graphics}
\usepackage{caption}
\usepackage{subfigure}
\usepackage[percent]{overpic}
\usepackage{changebar}
\usepackage{multirow}
\usepackage{setspace}
\usepackage{hyperref}
\usepackage{cite} 
\usepackage{indentfirst}

\usepackage{yhmath}%
 \usepackage{booktabs} %
 
 \usepackage{tikz}
\usepackage{verbatim}
\usetikzlibrary{arrows,backgrounds,decorations,shapes}
\usetikzlibrary{decorations.pathreplacing}
 \numberwithin{equation}{section}

\hypersetup{
    colorlinks=true,
    linkcolor=blue,
    citecolor=blue,
    urlcolor=blue
}
\graphicspath{{./}{./figures/}}
\allowdisplaybreaks

\newtheoremstyle{plainNoItalics}{}{}{\normalfont}{}{\bfseries}{.}{ }{}

\theoremstyle{plain}
\newtheorem{theorem}{Theorem}[section]

\counterwithin{figure}{section}
\theoremstyle{plainNoItalics}
\newtheorem{remark}[theorem]{Remark}

\newtheorem{propo}{Proposition}[section]

\newcommand{\be}{\begin{eqnarray}}
\newcommand{\ee}{\end{eqnarray}}
\newcommand{\beno}{\begin{eqnarray*}}
\newcommand{\eeno}{\end{eqnarray*}}

\usetikzlibrary{calc,trees,positioning,arrows,chains,shapes.geometric,%
	decorations.pathreplacing,decorations.pathmorphing,shapes,%
	matrix,shapes.symbols}
\tikzset{
	>=stealth',
	punktchain/.style={
		rectangle, 
		rounded corners, 
		draw=black, very thick,
		text width = 8em, 
		minimum height=2em, 
		text centered
	},
	yuan/.style={
		ellipse, 
		rounded corners, 
		draw=black, very thick,
		text width = 8em, 
		minimum height=2em, 
		text centered
	},
	line/.style={draw, thick, <-},
	element/.style={
		tape,
		top color=white,
		bottom color=blue!50!black!60!,
		minimum width=8em,
		draw=blue!40!black!90, very thick,
		text width=10em, 
		minimum height=3.5em, 
		text centered, 
		on chain},
	every join/.style={->, thick,shorten >=1pt},
	decoration={brace},
	tuborg/.style={decorate},
	tubnode/.style={midway, right=2pt},
	punkt/.style={
		rectangle, 
		rounded corners, 
		draw=black, very thick,
		text width=12em, 
		minimum height=3em, 
		text centered},
}

\makeatletter

\newcommand{\Rmnum}[1]{\expandafter\@slowromancap\romannumeral #1@}
\makeatother

\begin{document}
\baselineskip=1.5pc

\begin{center}
{\bf
  A high-order energy-conserving semi-Lagrangian discontinuous Galerkin method
for the Vlasov-Amp$\mathbf{\grave{e}}$re system
}
\end{center}

\vspace{.2in}
\centerline{
Xiaofeng Cai
  \footnote{Research Center of Mathematics, Advanced Institute of Natural Sciences, Beijing Normal University at Zhuhai, Zhuhai
519087, China} \footnote{Guangdong Provincial/Zhuhai Key Laboratory of Interdisciplinary Research and Application for Data Science, Beijing Normal-Hong Kong Baptist University, Zhuhai 519087, China. E-mail: xfcai@bnu.edu.cn. }
  ,
  Qingtao Li
  \footnote{School of Mathematical Sciences, East China Normal University, Shanghai, 200241, China. \\
  E-mail: 52265500022@stu.ecnu.edu.cn.}
  ,
  Hongtao Liu%
  \footnote{Center for Mathematical Plasma Astrophysics, Department of Mathematics, KU Leuven, Leuven, 3001, Belgium. E-mail: hongtao.liu@kuleuven.be.  }
  ,
  Haibiao Zheng%
  \footnote{School of Mathematical Sciences, Key Laboratory of MEA (Ministry of Education) \& Shanghai Key Laboratory of PMMP, East China Normal University, Shanghai 200241, China. \\
  E-mail: hbzheng@math.ecnu.edu.cn.}
}

\bigskip
\noindent
{\bf Abstract.}
In this paper, we propose a high-order energy-conserving  semi-Lagrangian discontinuous Galerkin
(ECSLDG) method for the Vlasov–Ampère system. The method employs a semi-Lagrangian discontinuous Galerkin scheme for spatial discretization of the Vlasov equation, achieving high-order accuracy while removing the Courant–Friedrichs–Lewy (CFL) constraint. To ensure energy conservation and eliminate the need to resolve the plasma period, we adopt an energy-conserving time discretization introduced by Liu et al. [J. Comput. Phys., 492 (2023), 112412]. Temporal accuracy is further enhanced through a high-order operator splitting strategy, yielding a method that is high-order accurate in both space and time. The resulting ECSLDG scheme is unconditionally stable and conserves both mass and energy at the fully discrete level, regardless of spatial or temporal resolution. Numerical experiments demonstrate the accuracy, stability, and conservation properties of the proposed method. In particular, the method achieves more accurate enforcement of Gauss’s law and improved numerical fidelity over low-order schemes, especially when using a large CFL number.

\vfill
\newpage

\section{Introduction}
The Vlasov equation serves as the theoretical foundation of collisionless plasma modeling, describing the evolution of the particle distribution function in a six‑dimensional phase space (three physics and three velocity dimensions). By self-consistent coupling with electromagnetic fields, it captures essential collective plasma phenomena, such as kinetic instabilities and wave–particle interactions \cite{inan2010principles}. Under the assumption of a non-magnetized or weakly magnetized plasma, coupling the Vlasov equation with Amp\`{e}re’s law yields the Vlasov–Amp\`{e}re (VA) model, one of fundamental kinetic models for collisionless electrostatic plasma.  However, its intrinsic high dimensionality and inherent nonlinearity generally make analytical solutions intractable, necessitating the development of robust and accurate numerical methods.

Numerical approaches for solving the Vlasov equation broadly fall into two categories: particle-based and grid-based approaches. 
The particle-in-cell (PIC) method,  a representative particle-based technique, models the plasma using finite macro-particles whose trajectories evolve under self-consistent electromagnetic fields \cite{GL2012}. Owing to its simplicity and efficiency, PIC has been widely applied and provides reasonable results at moderate computational cost\cite{SMGL2010,LG2018,quan2023influence,fathurrohim2025}. Nevertheless, as a stochastic method, it suffers from inherent numerical noise.
Alternatively, the grid based method directly solve the Vlasov equation in the full phase grids, which are deterministic and free from statistical noise, allowing for capturing fine-scale structures with high-order accuracy \cite{LQCS2020}.
Various grid-based solvers have been developed \cite{JQAC2010,LE2019,CYFF2014,NCLEE2015,NCLE2020,BJCF2020,xiong2014,liu2017,xiao2021,pu2024,Liu2019}, with comprehensive reviews available in \cite{FFES2003,dimarco2014}. 

Among the grid-based methods,  one promising method is the semi-Lagrangian (SL) method \cite{liu2023combined,YANG2021110632,GRANDGIRARD2006395,li2025semi,ZHENG2022114973}, which tracks particle trajectories in the spirit of Lagrangian methods while discretizing the distribution function on a fixed Eulerian grid. This strategy avoids strict Courant-Friedrichs-Lewy (CFL) constraints, enabling computationally efficient simulations of long-term plasma evolution. Among SL methods, the semi-Lagrangian discontinuous Galerkin (SLDG) approach has emerged as a powerful strategy \cite{cai2024ap,CGQ2019}. It inherits the local conservation, compact stencil, and high-order accuracy of discontinuous Galerkin (DG) methods, while also benefiting from the large time step capability of semi-Lagrangian techniques. These properties make SLDG method especially attractive for multidimensional Vlasov simulations, where it can accurately capture filamentation and steep gradients through the flexible use of high-order basis functions \cite{CBSQ2021,QS2011,RJSD2011,EL2019}.  


Most grid-based solvers employ explicit time integration schemes, which decouple the evolution of particles and fields. In each time step, particles are advanced using the previous fields, followed by field updates based on the new particle distribution. While this approach is straightforward to implement, it is constrained by a strict stability condition requiring the time step to resolve the plasma period. The necessity of resolving the smallest scales induced by plasma period, even if they are not of primary interest, makes explict methods computationally expensive for large-scale simulations. Moreover, even with small time steps, explicit schemes may suffer from numerical heating or cooling, which can compromise long-term accuracy and stability \cite{SMGL2010,liu2025ap}. To overcome these limitations, fully implicit kinetic methods have been developed \cite{cheng2014energy}. 
However, their efficiency and robustness are often limited by the complexity and cost of solving nonlinear systems.

A major direction in the development of kinetic solvers is to achieve energy conservation while efficiently relaxing the stability constraints of explicit schemes \cite{GL2017,chizhonkov2024,LIU2023112412}. Notable advances include semi-implicit energy-conserving particle-based approaches such as ECSim \cite{GL2017} and GEMPIC \cite{kormann2021energy}. 
However, extending these benefits to grid-based solvers remains a challenge, as grid-based methods and PIC methods fundamentally differ in their approach, even though both aim to achieve similar goals.
Recent progress includes the work of Yin et al. \cite{YZW2023}, who proposed an energy-conserving method for the Vlasov–Maxwell system based on the regularized moment method. Liu et al. \cite{LIU2023112412} proposed a energy-conserving semi-Lagrangian (ECSL)  method for VA system, which retains the efficiency of explicit schemes while inheriting the stability and conservation properties of implicit methods.
This framework has subsequently been extended to the full Vlasov–Maxwell system \cite{liu2025ecsl}.
Extending the ECSL framework to the SLDG method is particularly attractive, as it enhances spatial accuracy through local reconstruction and is well-suited for efficient parallelization in high-dimensional kinetic simulations.

The main contribution of this paper is the development of a high-order energy-conserving SLDG
(ECSLDG) method for the VA system.  The method integrates the energy-conserving time discretization of the ECSL framework \cite{LIU2023112412} with the high-order spatial accuracy of the DG method, achieving high-order temporal accuracy by adopting a suitable fourth-order operator splitting scheme.
Specifically, the VA system is decomposed into two energy-conserving subsystems. In each subsystem, the multidimensional Vlasov equation is reduced to a set of one-dimensional (1D) advection equations, which are efficiently solved using the SLDG scheme. This strategy ensures mass conservation, eliminates CFL restrictions, and maintains high-order spatial accuracy. Moreover, inspired by \cite{LIU2023112412}, we semi-implicitly couple the Amp\`{e}re equation with the moments of the Vlasov equation to solve the electric field, effectively removing time step constraints imposed by the plasma period. A suitable fourth-order composition method is employed for time integration, effectively balancing accuracy and stability among several established schemes for the entire VA system. 
The ECSLDG method thus simultaneously achieves unconditional stability, exact conservation of total energy and mass, and high-order accuracy in both space and time. Finally, several numerical experiments are presented to demonstrate the accuracy, efficiency, and conservation properties of the proposed method. 

 The rest of the paper is organized as follows. Section \ref{sec:level2} introduces the VA system. In Section \ref{sec:level3}, we present the high-order ECSLDG method. In Section \ref{section:numerical}, we present the numerical results. Finally, a summary is given in Section \ref{section:conclusion}.

\section{The Vlasov-Amp\`{e}re system \label{sec:level2}}
In this section, we review the VA system, which is given by
\begin{align}
&\partial_t f_s+\bm{v}_s\cdot \nabla_{\bm{x}}f_s+\frac{q_s}{m_s}\bm{E}\cdot \nabla_{\bm{v}_s}f_s=0, \label{Vlasov}\\
&\varepsilon \partial_t \bm{E}= -\bm{J}, \label{Ampere}
\end{align}
where $f_s(\bm{x},\bm{v},t)$ describes the velocity distribution function of species $s$ (choosing $e$ for electrons and $i$ for ions) as they move in a $d$-dimensional($d=1,2,3$) velocity space with velocity $\bm{v}_s$ at position $\bm{x}$ and time $t$. The parameters $q_s$  and $m_s$ denote the charge and mass of spieces $s$, $\bm{E}$ is the electric field, and $\epsilon$ is the permittivity. The total current density is given by $\bm{J}=\sum_s q_sp_s$, where  $p_s= \int_{\Omega_{\bm{v}}} f_s\bm{v}_s d\bm{v}$. The domain is defined as $\Omega=\Omega_{\bm{x}} \times \Omega_{\bm{v}}$, where $\Omega_{\bm{x}}=\mathbb{R}^d$  represents the physics space and $\Omega_{\bm{v}}=\mathbb{R}^d$ represents the velocity domain. For simplicity, we consider periodic boundary conditions and assume $f(\bm{x},\bm{v},t)$ has compact support in the velocity space. 

Integrating Eq. \eqref{Vlasov} over $\bm v_s$ and summing over species yields the continuity equation
\begin{align}\label{con-eq}
\partial_t\rho+ \nabla \cdot \bm{J}=0,
\end{align}
where $\rho= \sum_s q_s n_s$ is the total charge density, and $n_s=  \int_{\Omega_{\bm{v}}} f_s d\bm{v}$ is the particle number density. Combining Eq. (\ref{Ampere}) with Eq. (\ref{con-eq}) leads to the Poisson equation,
\begin{align}\label{poisson}
-\varepsilon \Delta \phi=\rho,
\end{align}
where the electric field $\bm{E}$ and the electric potential $\phi$ satisfy $\nabla \phi=-\bm{E}$. It is well known that the Vlasov-Poisson (VP) system is constituted by Eq. (\ref{Vlasov}) and Eq. (\ref{poisson}). When the charge continuity equation (\ref{con-eq}) is satisfied, the VA system is equivalent to VP system.

In this paper, we focus on electron dynamics with $q_e=-e$, assuming ions with charge $q_i=e$ form a uniform neutral background. The Debye length $\lambda$ and the electron plasma frequency $\omega_p$ are defined by
\begin{align*}
\lambda=\biggl(\frac{\varepsilon K_B T_e}{e^2n_e}\biggr)^{1/2}, && \omega_p=\biggl( \frac{n_ee^2}{\varepsilon m_e}\biggr)^{1/2},
\end{align*}
where $K_B$ is the Boltzmann constant and $T_e$ is the electron temperature.

In order to normalize the VA system, the following dimensionless variables are defined
\begin{align*}
\bar{\bm{x}}=\frac{\bm{x}}{{x}_0}, && \bar{T}=\frac{T}{T_0}, && \bar{m}=\frac{m_e}{m_0}, && \bar{n}=\frac{n}{n_0},\\
\bar{\bm{v}}=\frac{\bm{v}}{{v}_0}, && \bar{t}=\frac{t}{t_0}, && \bar{f}=\frac{f}{f_0}, && \bar{\bm{E}}=\frac{\bm{E}}{\mathit{E}_0},
\end{align*}
where ${x}_0,~T_0,~m_0,~n_0$ are  reference length, temperature, mass, number density. Besides, we choose ${v}_0=\sqrt{K_BT_0/m_0},~t_0=x_0/v_0,~{E}_0=K_BT_0/e{x}_0,~f_0=n_0/{v}_0^d$. 

In this paper, unless otherwise stated, we choose the characteristic parameters as $m_0=m_e,~n_0=n_e,~T_0=T_e$. Additionally, the Debye length and plasma frequency can be further normalized as $\bar{\lambda}=\lambda/x_0,\bar{\omega}_p=1/\bar{\lambda}$. As a result, the dimensionless form of the VA system becomes
\begin{equation}\label{VA-non}
    \begin{aligned}
        &\partial_t f+\bm{v}\cdot \nabla_{\bm{x}}f-\bm{E}\cdot \nabla_{\bm{v}}f=0,\\
        &\lambda^2 \partial_t \bm{E}= -\bm{J}.
    \end{aligned}
\end{equation}
It should be noted that, for convenience, we have omitted the overbars on all dimensionless variables. Additionally, the subscript
$s$ on the distribution function $f$ is dropped, since only electron dynamics are considered in this study. Accordingly, the dimensionless form of Eq. (\ref{poisson}) is given by
\begin{align}\label{poisson-non}
-\lambda^2 \Delta \phi=1-n.
\end{align}
Note that Eq.~(\ref{poisson-non}) is used to provide the initial electric field $\bm{E}$ for the VA system in the numerical experiments presented later.

It is well known that the VA system conserves the total energy $E_{total}$
\begin{align*}
E_{total}=\frac{1}{2} \int_{\Omega_{\bm{x}}}\int_{\Omega_{\bm{v}}} f\bm{v}^2d\bm{v} d\bm{x} +\frac{\lambda^2}{2} \int_{\Omega_{\bm{x}}} \bm{E}^2d\bm{x},
\end{align*}
which is composed of the kinetic energy $E_K= \frac{1}{2} \int_{\Omega_{\bm{x}}}\int_{\Omega_{\bm{v}}} f\bm{v}^2d\bm{v} d\bm{x}$  and the electric energy $\textit{E}_{\bm{E}}= \frac{\lambda^2}{2} \int_{\Omega_{\bm{x}}} \bm{E}^2d\bm{x}$. Furthermore, it is notable that any function taking the form of $ \int_{\Omega_{\bm{x}}}\int_{\Omega_{\bm{v}}}\bm{G}(f)d\bm{v} d\bm{x}$ is conserved. In particular, $L_1= \int_{\Omega_{\bm{x}}}\int_{\Omega_{\bm{v}}}fd\bm{v} d\bm{x}$ represents the total mass, while $P= \int_{\Omega_{\bm{x}}}\int_{\Omega_{\bm{v}}} f\bm{v}d\bm{v} d\bm{x}$ represents the total momentum. 
However, developing numerical methods that can accurately and simultaneously preserve mass, momentum, and total energy remains a significant challenge. Furthermore, traditional explicit kinetic solvers impose a severe stability constraint on the time step size \( \Delta t \), which must resolve the normalized Debye length \( \lambda \) to maintain accuracy and stability.

In the following section, we aim to construct a kinetic scheme for the VA system (\ref{VA-non}) that conserves both mass and total energy, while also relaxing  strict time step restrictions associated with explicit methods.
\section{The ECSLDG method for Vlasov-Amp\texorpdfstring{$\mathbf{\grave{e}}$}{e}re system \label{sec:level3}}
In this section, we propose the ECSLDG method for the VA system. First, we review the two-dimensional (2D) SLDG method based on the operator splitting method for solving linear transport problems. Then, we introduce the ECSL method. Subsequently, we present the proposed ECSLDG method for the VA system. Finally, we extend the ECSLDG method to achieve high-order accuracy in time.

\subsection{2D SLDG with operator splitting for linear transport problems}
\subsubsection{1D SLDG method}
To introduce the 2D SLDG algorithm with operator splitting, the one-dimensional linear transport problem is expanded upon as follows
\begin{align}\label{1d-transport}
u_t + (a(x,t)u)_x = 0,~~~~(x,t) \in [x_l,x_r]\times(0,T].
\end{align}

Let $x_l=x_{\frac{1}{2}}< \cdots < x_{N+\frac{1}{2}}=x_r$ be a partition of the interval $[x_l,x_r]$, denoted by $\mathcal{E}_h$. 
Define $I_j=[x_{j-\frac{1}{2}},x_{j+\frac{1}{2}}]$ ($j=1,\cdots,N$) and $\Delta x_j=x_{j+\frac{1}{2}}-x_{j-\frac{1}{2}}$. The space $\mathcal{D}_k(\mathcal{E}_h)$ is defined as the space of piecewise discontinuous polynomials of degree $k$
\begin{align*}
\mathcal{D}_k(\mathcal{E}_h) = \{ v:v|_{I_j} \in P^k(I_j),~~~~\forall ~j=1,\cdots,N\},
\end{align*}
where $P^k(I_j)$ is the space of polynomials of degree $k$ on $I_j$. In addition, let $0=t_0<t_1<\cdots<t_M=T$, and,
for notational simplicity in this exposition, set a uniform time step $\Delta t=T/M$  (in numerical experiments, $\Delta t$ may vary).

To find suitable test functions $\psi(x,t),~t\in[t_m,t_{m+1}]$ ($m=0,\cdots,M-1$), the adjoint problem of Eq.~(\ref{1d-transport}) is considered
\begin{equation}\label{adjoint}
\begin{aligned}
&\partial_t \psi+a(x,t)\psi_x=0,\\
&\psi(t=t_{m+1})=\Psi(x)\in \mathcal{D}_k(\mathcal{E}_h).
\end{aligned}
\end{equation}

\begin{figure}[htbp]
    \centering
    \begin{tikzpicture}
	tikzstyle {every node }=[font=very small];
	\draw[line width=1pt] (0,-1) -- (8,-1);
	\draw[line width=1pt] (0,-3) -- (8,-3);
	\draw[line width=0.5pt] (3,-2.9) -- (3,-3.1);
	\draw[line width=0.5pt] (6,-2.9) -- (6,-3.1);
	\draw[line width=0.5pt] (3,-0.9) -- (3,-1.1);
	\draw[line width=0.5pt] (6,-0.9) -- (6,-1.1);
	\coordinate [label=$x_{j-{\frac{1}{2}}}$] (O) at (3,-1);
	\coordinate [label=$x_{j+\frac{1}{2}}$] (O) at (6,-1);
	   \node (A) at(1,-3.3) {$x^\star_{j-\frac{1}{2}}$};
        \node (B) at(4,-3.3) {$x^\star_{j+\frac{1}{2}}$};
         \draw[decorate,decoration={brace,mirror},,line width=1pt,blue] (1,-3) -- (4,-3);
         \node at(2.5,-3.5) {$I_j^{\star}$};
        \node at (9,-1) {$t_{m+1}$};
        \node at(9,-3) {$t_m$};
        \draw[loosely dashed,  line width=0.5pt,->] (3,-1) -- (1,-3);
        \draw[loosely dashed,  line width=0.5pt,  <-] (4,-3) -- (6,-1);
        \filldraw[fill=blue!10] (3,-1) --(1,-3) -- (4,-3) -- (6,-1)--cycle;
        \node at(3.5,-2) {$\tilde{I}_{j}(t)$};	
        \node at(4.5,-0.75) {$I_{j}$};	
    \end{tikzpicture}
    \caption{The integral region $\tilde{I}_j(t)\times [t_m,t_{m+1}] $ related to space-time variables.}
    \label{inte-re}
\end{figure}
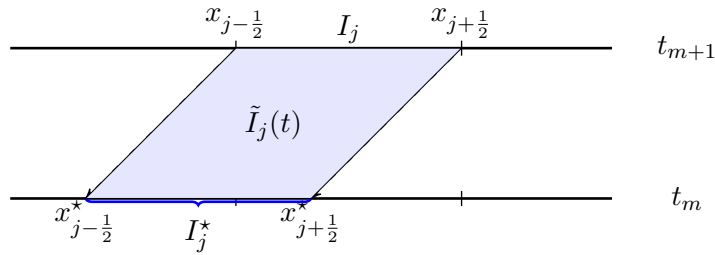

Combining Eq.~(\ref{1d-transport}) with Eq.~(\ref{adjoint}) and employing the Reynolds transport theorem, it can be obtained in \cite{Guo2014ACS} that
\begin{align}\label{int-inv}
\displaystyle \frac{d}{dt}\int_{\tilde{I}_j(t)}u(x,t)\psi(x,t)dx=0,
\end{align}
where the construction of the region $\tilde{I}_j(t)$ is as follows: the two endpoints of the cell $I_j$ at $t_{m+1}$ are backtracked along the characteristic curves $\frac{dx}{dt}=a(x,t)$ to time $t_m$ to obtain the corresponding endpoints and determine the cell $I_j^{\star}=[x_{j-\frac{1}{2}}^{\star},x_{j+\frac{1}{2}}^{\star}]$. The region $\tilde{I}_j(t)$ is enclosed by the cell $I_j$, the cell $I_j^{\star}$, and the two characteristic curves (see Fig.~\ref{inte-re}).

By Eq.~(\ref{int-inv}), the fully discrete SLDG scheme
seeks  $u_h^{m+1} \in \mathcal{D}_k(\mathcal{E}_h)$ such that, for any $\Psi(x) \in \mathcal{D}_k(\mathcal{E}_h)$, 
\begin{align}\label{SL-int}
\displaystyle \int_{I_j}u_h^{m+1}\Psi(x)dx = \displaystyle \int_{I_j^{\star}}u_h^m\psi(x,t_m)dx,
\end{align}
where $u_h^m$ is the DG solution at $t_{m}$. The key to updating $u_h^{m+1}$ lies in computing the right-hand side of Eq.~(\ref{SL-int}), 
which proceeds as follows:

\begin{enumerate}

    \item  Interpolating the test function $\psi(x,t_m)$:

    \begin{itemize}
        \item [\textbf{a.}] Assuming $k=2$, three points $x_{j,q}$ ($q=1,2,3$), such as Gauss-Lobatto points, are selected in cell $I_j$, which are then made to follow the characteristic curves to locate the feet $x_{j,q}^{\star}$ at time $t_m$ as shown in Fig. \ref{step-sldg}(a). The feet are obtained by calculating the subsequent equation with a suitable numerical integrator such as a fourth-order Runge-Kutta method in \cite{EGC2006}
\begin{equation}\label{ODE}
\begin{aligned}
&\frac{dx}{dt}=a(x,t),\\
&x(t_{m+1})=x_{j,q}.
\end{aligned}
\end{equation}
        \item[\textbf{b.}] Since the test function $\psi(x,t_m)$ satisfies $\frac{d\psi}{dt} = 0$ along characteristic curves $\frac{dx}{dt}=a(x,t)$ (i.e., remains constant along these curves), it follows that $\psi(x_{j,q}^{\star}, t_m) = \Psi(x_{j,q})$. Using these relationship, we can obtain a unique third-order test function $\psi^{\star}(x)$ that interpolates $\psi(x,t_m)$ on the upstream cell $I_j^{\star}$.
        \begin{figure}[htbp]
	\begin{center}
		\begin{minipage}{0.4\linewidth}
			\begin{tikzpicture}
	           tikzstyle {every node }=[font=very small];
	           \draw[line width=1pt] (0,-1) -- (6.5,-1);
	           \draw[line width=1pt] (0,-3) -- (6.5,-3);
	           \draw[line width=0.5pt] (3,-2.9) -- (3,-3.1);
	           \draw[line width=0.5pt] (3,-0.9) -- (3,-1.1);
	           \draw[line width=0.5pt] (6,-0.9) -- (6,-1.1);
	           \coordinate [label=$x_{j-\frac{1}{2}}$($\textcolor{red}{x_{j,1}}$)] (O) at (3,-1);
                \fill[black] (3,-1) circle (2pt);
	           \coordinate [label=$x_{j+\frac{1}{2}}$($\textcolor{red}{x_{j,3}}$)] (O) at (6,-1);
                \coordinate [label=$\textcolor{red}{x_{j,2}}$] (O) at (4.5,-0.9);
                \fill[black] (6,-1) circle (2pt);
                \fill[black] (4.5,-1) circle (2pt);
                \fill[blue] (1,-3) circle (2pt);
                \fill[blue] (4,-3) circle (2pt);
                \fill[blue] (2.5,-3) circle (2pt);
	        \node (A) at(1,-3.3) {$x^\star_{j-\frac{1}{2}}$($\textcolor{red}{x^\star_{j,1}}$)};
                \node (B) at(4,-3.3) {$x^\star_{j+\frac{1}{2}}$($\textcolor{red}{x^\star_{j,3}}$)};
                \node (C) at(2.5,-3.3) {$\textcolor{red}{x^\star_{j,2}}$};
                \node at (6.9,-1) {$t_{m+1}$};
                \node at(6.75,-3) {$t_m$};
                \draw[loosely dashed, line width=1pt, -] (3,-1) ..controls (1.8,-1.5) .. (1,-3);
                \draw[loosely dashed, line width=1pt, -] (4.5,-1) ..controls (3.3,-1.5) .. (2.5,-3);
                \draw[loosely dashed, line width=1pt,  -] (6,-1) ..controls (4.8,-1.5) .. (4,-3);
            \end{tikzpicture}
            \centerline{(a)}
		\end{minipage}
        \hspace{1cm}
		\begin{minipage}{0.4\linewidth}
			\begin{tikzpicture}
	           tikzstyle {every node }=[font=very small];
	           \draw[line width=1pt] (0,-1) -- (6.5,-1);
	           \draw[line width=1pt] (0,-3) -- (6.5,-3);
	           \draw[line width=0.5pt] (3,-2.9) -- (3,-3.1);
	           \draw[line width=0.5pt] (3,-0.9) -- (3,-1.1);
	           \draw[line width=0.5pt] (6,-0.9) -- (6,-1.1);
	           \coordinate [label=$x_{j-\frac{1}{2}}$] (O) at (3,-1);
	           \coordinate [label=$x_{j+\frac{1}{2}}$] (O) at (6,-1);
	           \node (A) at(0.5,-2.7) {$x^\star_{j-\frac{1}{2}}$};
                \node (B) at(3.5,-2.7) {$x^\star_{j+\frac{1}{2}}$};
                \draw[decorate,decoration={brace,mirror},line width=1pt, red] (1,-3) -- (3,-3);
                \draw[decorate,decoration={brace,mirror},line width=1pt, red] (3,-3) -- (4,-3);
                \node (C) at(2,-3.4) {$I^\star_{j,1}$};
                \node (D) at(3.5,-3.4) {$I^\star_{j,2}$};
                \node at (6.9,-1) {$t_{m+1}$};
                \node at(6.75,-3) {$t_m$};
                \draw[loosely dashed, line width=0.5pt,->] (3,-1) ..controls (1.8,-1.5) .. (1,-3);
                \draw[loosely dashed, line width=0.5pt,  ->] (6,-1) ..controls (4.8,-1.5) .. (4,-3);
            \end{tikzpicture}
            \centerline{(b)}
		\end{minipage}
	\end{center}
	\caption{Schematic illustration for one-dimensional SLDG schemes.}
	\label{step-sldg}
\end{figure}
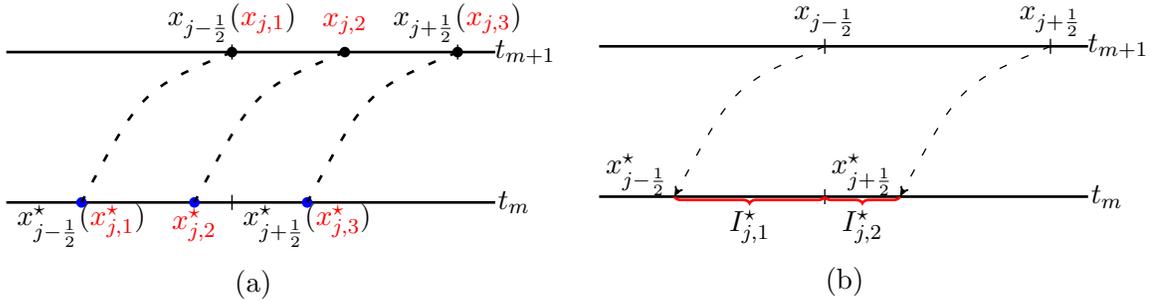

    \end{itemize}
   
    \item  Calculating the right-hand side of Eq.~(\ref{inte-re}):

    Searching the information of the overlap to the upstream cell $I_j^{\star}$ and the grid at time $t_{m+1}$. For instance, the number of cells at $t_{m+1}$ involved with $I_j^{\star}$, denoted as $l$, is evident in Fig.~\ref{step-sldg}(b), where $l=2$, that is, $I_j^{\star} = I_{j,1}^\star \cup I_{j,2}^\star$. Consequently, there is
\begin{align*}
\displaystyle \int_{I_j^{\star}} u_h^m\psi(x,t_m)dx \approx \sum_{i=1}^l \int_{I_{j,i}^{\star}} u_h^m\psi^\star(x)dx.
\end{align*}
\end{enumerate}

In addition to the 1D SLDG framework described above, other variants are detailed in \cite{QS2011,cai2022comparison}.

\subsubsection{2D SLDG with operator splitting }
We now turn to the 2D SLDG method with Strang splitting (SS).
 Consider the   linear problem
\begin{align}\label{2d-linear}
\frac{\partial u}{\partial t} +(a(x,y,t)u)_x +(b(x,y,t)u)_y = 0,~~~~(x,y,t)\in [x_l,x_r]\times[y_b,y_t]\times[0,T],
\end{align}
with periodic boundary conditions and initial conditions. Here, $(a(x,y,t),b(x,y,t))$  denotes the given velocity field.

Let $\mathcal{E}_{x,h}:=\{x_{i-\frac{1}{2}}\}_1^{N_x+1}$ be a partition of the interval $[x_l,x_r]$, and $\mathcal{E}_{y,h}:=\{y_{j-\frac{1}{2}}\}_1^{N_y+1}$ be a partition of the interval $[y_b,y_t]$. Then the   regular rectangular partition $\mathcal{E}_h$ of the computational domain $[x_l,x_r]\times[y_b,y_t]$ is expressed as $\mathcal{E}_h = \mathcal{E}_{x,h} \otimes \mathcal{E}_{y,h}$. Let $\mathcal{D}_k(\mathcal{E}_{x,h})$ and $\mathcal{D}_k(\mathcal{E}_{y,h})$ be denoted by the space of piecewise discontinuous polynomials of degree $k$ on $\mathcal{E}_{x,h}$ and $\mathcal{E}_{y,h}$, respectively. In the dimensional splitting setting, $Q^k(\mathcal{E}_h)= \mathcal{D}_k(\mathcal{E}_{x,h}) \otimes \mathcal{D}_k(\mathcal{E}_{y,h})$ is denoted by the space of piecewise discontinuous polynomials on $\mathcal{E}_h$. Therefore, there are $(k+1)^2$ degrees of freedom per computational cell.

Eq.~(\ref{2d-linear}) can be split into two 1D problems in a conservative form in \cite{QS2011}:
\begin{align}
    &\frac{\partial u}{\partial t} +(a(x,y,t)u)_x = 0,\label{2d-1d-1}\\
    &\frac{\partial u}{\partial t} +(b(x,y,t)u)_y = 0\label{2d-1d-2}.
\end{align}

The 2D SLDG framework proceeds as follows:

\begin{enumerate}
    \item {\bf Quadrature nodes.}

 On each rectangular cell
$[x_{i-\frac{1}{2}},x_{i+\frac{1}{2}}]\times [y_{j-\frac{1}{2}},y_{j+\frac{1}{2}}]$
select the $(k+1)^2$ tensor-product Gaussian nodes
$ \{   (x_{i,p},y_{j,q})  | p,q=1,\cdots,k+1  \} $.

 \item {\bf Strang splitting.}

 Eq.~(\ref{2d-1d-1})-Eq.~(\ref{2d-1d-2}) is computed using Strang Splitting over a time step $\Delta t$:

\begin{description}
    \item[a.] {\bf $x-$direction half step.}

For each fixed $y_{j,q}$ ($j=1,\cdots,N_y$, $q=1,\cdots,k+1$), 
evolve
\begin{align*}
    \frac{\partial u}{\partial t} +(a(x,y_{j,q},t)u)_x = 0
\end{align*}
by the 1D SLDG method for a half time step of size $\Delta t/2$.

   \item[b.] {\bf $y-$direction full step.}

For each fixed $x_{i,p}$ ($i=1,\cdots,N_x$, $p=1,\cdots,k+1$), 
evolve
\begin{align*}
    \frac{\partial u}{\partial t} +(b(x_{i,p},y,t)u)_y = 0
\end{align*}
 by the 1D SLDG method for a full step of size $\Delta t$.

\item[c.] {\bf $x-$direction half step.}

 Repeat step a. for another half step $\Delta t/2$.
    
\end{description}

\end{enumerate}

For  detailed discussions on the operator-splitting SLDG approach, see \cite{Guo2014ACS,cai2022comparison}.

\subsection{ECSLDG for VA system} 

Here, we proposed the ECSLDG method for the VA system under 1D1V setting,  which is composed of the ECSL method \cite{LIU2023112412} and the SLDG method. The computational domain is $\Omega = \Omega_x \times \Omega_v=[0,L] \times [-v_m,v_m]$ ($v_m$ should be sufficiently large to ensure that $f(x,v,0)$ is compactly supported on $[-v_m,v_m]$).

Let $0=x_{\frac{1}{2}}<x_{\frac{3}{2}}<\cdots<x_{N_x+\frac{1}{2}}=L$ be a partition in the $x$-direction, denoted by $\mathcal{E}_{x}$, and $-v_m=v_{\frac{1}{2}}<v_{\frac{3}{2}}<\cdots<v_{N_v+\frac{1}{2}}=v_m$ be a partition in the $v$-direction, denoted by $\mathcal{E}_v$. Let $\mathcal{D}_k(\mathcal{E}_{x})$ and $\mathcal{D}_k(\mathcal{E}_{v})$ be denoted by the space of piecewise discontinuous polynomials of degree $k$ on $\mathcal{E}_{x}$ and $\mathcal{E}_{v}$, respectively. The grid is defined as 
\begin{align*}
	&C_{i,j} := [x_{i-\frac{1}{2}},x_{i+\frac{1}{2}}]\times [v_{j-\frac{1}{2}},v_{j+\frac{1}{2}}],&~\\
	&X_i:=[x_{i-\frac{1}{2}},x_{i+\frac{1}{2}}],&V_j:=[v_{j-\frac{1}{2}},v_{j+\frac{1}{2}}],~~~~i=1,\cdots,N_x,~j=1,\cdots,N_v,
\end{align*}
and define $\Delta x_i: =x_{i+\frac{1}{2}}-x_{i-\frac{1}{2}}$ and $\Delta v_j:=v_{j+\frac{1}{2}}-v_{j-\frac{1}{2}}$. There are $(k+1)^2$ Gaussian quadrature points on $C_{i,j} $.

To ensure the total energy conservation of the VA system, an operator splitting technique is ingeniously employed in \cite{cheng2014energy} to decompose the VA system into two subsystems\textrm{:}
the Hamiltonian $\mathcal{H}_f$ system
\begin{equation}\label{Hf}
\begin{aligned}
&\partial_tf+v\cdot \nabla_xf=0,\\
&\lambda^2\partial_t E=0,
\end{aligned}
\end{equation}
and the Hamiltonian $\mathcal{H}_E$ system
\begin{equation}\label{HE}
\begin{aligned}
&\partial_tf-E\cdot \nabla_vf=0,\\
&\lambda^2\partial_tE=-J.
\end{aligned}
\end{equation}


Note that system (\ref{Hf}) represents the particle transport equation, which can be solved using the conservative semi-Lagrangian method. The key idea of the ECSL method \cite{LIU2023112412} is to extract the first moment of the Vlasov equation from system (\ref{HE}) and couple it with the Amp\`{e}re equation. By applying the implicit midpoint method, we derive the following linear system
\begin{equation}
\begin{array}{l}
J^{k+1}= J^k +\Delta t  E^{k+1/2}n^{k+1/2},\\
{\lambda ^2}  E^{k+1}= {\lambda ^2}E^{k} -\Delta t  J^{k+1/2},
\end{array}
\label{eq:12}
\end{equation}
where ${{{J}}^{m + 1/2}} = ({{{J}}^{m + 1}} + {{{J}}^m})/2$ and ${{{E}}^{m + 1/2}} = ({{{E}}^{m + 1}} + {{{E}}^m})/2$. Besides, $n^{k+1/2}= n^{k}$ since $n^{k+1}= n^{k}$ holds true for system (\ref{HE}). By solving Eq.~(\ref{eq:12}), we can easily obtain the $E^{m+1}$ as follow,
\begin{align*}
	E^{m+1}=\frac{\lambda^2-\theta}{\lambda^2+\theta}E^m-\frac{\Delta t}{\lambda^2+\theta}J^m,
\end{align*}
where $\theta = \frac{1}{4}n^m\Delta t^2$. Then, the distribution function can be updated by using
the semi-Lagrangian method.


We now extend the core idea of ECSL method to a discontinuous Galerkin formulation and introduce the ECSLDG method for solving the Vlasov–Ampère system as follows:
\begin{enumerate}
    \item {\bf Quadrature nodes.}
    
     On each rectangular cell
$C_{i,j}$ select the $(k+1)^2$ tensor-product Gaussian nodes $ \{   (x_{i,p},v_{j,q})  | p,q=1,\cdots,k+1  \} $.
    \item {\bf Strang splitting.}

    Combining 1D SLDG, the ECSLDG method for solving the VA system from $t_m$ to $t_{m+1}$ is as follows:
    \begin{itemize}
        \item [\bf{Step 1:}] {\bf $x-$direction half step.}

            For each fixed $v_{j,q}$ ($j=1,\cdots,N_y$, $q=1,\cdots,k+1$),
            \begin{itemize}
                \item [a.] 
                For $i=1,\cdots,N_x$, seek $f_h^*(x,v_{j,q}) \in \mathcal{D}_k(\mathcal{E}_x)$ such that
                \begin{align}\label{Hf-SLDG}
                    \displaystyle \int_{X_i}f_h^{*}(x,v_{j,q})\Psi(x)dx=\int_{X_i^\star}f_h^m(x,v_{j,q})\Psi(x+v_{j,q}\Delta t/2)dx,~~~~\forall\Psi(x)\in \mathcal{D}_k(\mathcal{E}_x).
                \end{align}
                \item [b.] Compute Eq.~(\ref{Hf-SLDG}) to obtain the distribution function $f_h^{*}(x,v_{j,q})$ for a half time step of size $\Delta t/2$.
            \end{itemize}
        \item [\bf{Step 2:}] {\bf $v-$direction full step.}

        For each fixed $x_{i,p}$ ($i=1,\cdots,N_x$, $p=1,\cdots,k+1$), 
        \begin{itemize}
            \item [a.] Update the electric field $E_h^{m+1}(x_{i,p})$ and obtain the electric field $E_h^{m+1/2}(x_{i,p})$ for solving the subsystem Eq.~(\ref{HE}):
                \begin{align*}
	               &E_h^{m+1}(x_{i,p})= \frac{\lambda^2-\theta(x_{i,p})}{\lambda^2+\theta(x_{i,p})}E_h^{m}(x_{i,p})-\frac{\Delta t}{\lambda^2+\theta(x_{i,p})}J_h^*(x_{i,p}),\\
	               &E_h^{m+1/2}(x_{i,p}) =(E_h^m(x_{i,p})+E_h^{m+1}(x_{i,p}))/2,
                \end{align*}
                where $\displaystyle \theta(x_{i,p})=\frac{\Delta t^2}{4}n_h^*(x_{i,p})=\frac{\Delta t^2}{4}\sum\limits_{j=1}^{N_v}\sum\limits_{q=1}^{k+1}w_{j,q}f_h^*(x_{i,p},v_{j,q})$ ($w_{j,q}$ is the corresponding Gaussian weight of Gaussian point $v_{j,q}$) and the current density $J_h^*(x_{i,p})=\displaystyle-\sum\limits_{j=1}^{N_v}\sum\limits_{q=1}^{k+1}w_{j,q}f_h^*(x_{i,p},v_{j,q})v_{j,q}$ are obtained from the distribution function $f_h^{*}(x,v_{j,q})$ of \textbf{Step 1}.
                \item[b.]  
               For $j=1,\cdots,N_v$, seek $f_h^{**}(x_{i,p},v) \in \mathcal{D}_k(\mathcal{E}_v)$ such that
               \begin{align}\label{HE-SLDG}
                   \displaystyle \int_{V_j}f_h^{**}(x_{i,p},v)\Psi(v)dv=\int_{V_j^\star}f_h^*(x_{i,p},v)\Psi(v - E_h^{m+1/2}(x_{i,p})\Delta t)dv,~~~~\forall\Psi(v)\in \mathcal{D}_k(\mathcal{E}_v).
               \end{align}
               \item[c.] Compute Eq. (\ref{HE-SLDG}) to obtain the distribution function function $f_h^{**}(x_{i,p},v)$ for a time step $\Delta t$.
        \end{itemize}
        \item[\bf{Step 3:}] {\bf $x-$direction half step.}
        
        Repeat \textbf{Step 1} for another half step $\Delta t/2$.
    \end{itemize}
\end{enumerate}





Here, we use $\mathcal{H}_f^{SLDG}(\Delta t/2)$ to denote \textbf{Step 1}, and $\mathcal{H}_E^{SLDG}(\Delta t)$ to denote \textbf{Step 2}. Then, we can obtain a second-order ECSLDG methods for the VA system by using the Strang splitting method as follows
\begin{align}\label{2or-SLDG}
    \mathcal{H}^{SLDG}_{SS}(\Delta t)=\mathcal{H}_f^{SLDG}(\Delta t/2)\mathcal{H}_E^{SLDG}(\Delta t)\mathcal{H}_f^{SLDG}(\Delta t/2).
\end{align}

\begin{remark}
    To implement the algorithm, we need to obtain the initial information of the electric field. Solve Eq.~(\ref{poisson-non}) using the LDG method\cite{ADBFC2002,CBS1998,CPCBP2000} to obtain the initial values of the electric field $E_h^0(x_{i,p})$ at $k+1$ Gaussian nodes $x_{i,p}$ ($i=1,\cdots,N_x$, $p=1,\cdots,k+1$) in the $x$-direction for each cell. It is noted that all functions with the subscript $h$ in this paper denote the numerical solution.
\end{remark}
\begin{remark}
     Gaussian numerical integration with $k+1$ Gaussian points is exact for polynomials of degree less than $2k+2$.
\end{remark}

The introduction to ECSL method mentioned above only covers the parts necessary for our method. For a more detailed explanation of this technology, please refer to \cite{LIU2023112412}.

\subsection{Properties of the ECSLDG}
In this subsection, we present some properties of the ECSLDG method. Without loss of generality, we focus here on the 1D1V VA system with periodic boundary condition. 


\begin{propo}\label{mass-prop}
	(Total particle number conservation). The ECSLDG method preserves the total particle number of the VA system, i.e.,
	\begin{align*}
			L_1^{m+1} = \displaystyle \int_{0}^L\int_{-v_m}^{v_m} f_h^{m+1}(x,v) dvdx=\int_{0}^L\int_{-v_m}^{v_m} f_h^{m}(x,v) dvdx=L_1^m,~~~m=0,1,\cdots,M-1.
	\end{align*} 
	\begin{proof}
		For \textbf{Step 1}, 
        taking $\Psi(x)=1$ in Eq.~(\ref{Hf-SLDG}), we get 
		\begin{align*}
            \displaystyle \int_{X_i}f_h^{*}(x,v_{j,q})dx=  \int_{X_i^{\star}}f_h^m(x,v_{j,q})dx,
		\end{align*}
        where $X_i^{\star}=[x_{i-\frac{1}{2}}-v_{j,q}\Delta t/2,x_{i+\frac{1}{2}}-v_{j,q}\Delta t/2]$. 
        \newline
        Under the previous framework of the SLDG method, and without loss of generality, we assume $v_{j,q}>0$ and $v_{j,q}\Delta t/2<\Delta x$. Summing above equality from $i=1$ to $i=N_x$, we get 
         \begin{align*}
             &\displaystyle \int_{0}^Lf^*_h(x,v_{j,q})dx=\int_{-v_{j,q}\Delta t/2}^{L-v_{j,q}\Delta t/2}f_h^{m}(x,v_{j,q})dx
             =\sum_{i=1}^{N_x}\int_{X_i^{\star}}f_h^m(x,v_{j,q})dx
             = \sum_{i=1}^{N_x}\sum_{l=1}^2 \int_{X_{i,l}^{\star}} f_h^m(x,v_{j,q})dx\\             
             &=\int_{0}^Lf_h^m(x,v_{j,q})dx= \int_{0}^Lf_h^*(x,v_{j,q})dx.
         \end{align*}         
        Integrating in the $v$-direction, we obtain
		\begin{align*}
			&\displaystyle \int_{-v_m}^{v_m}\int_{0}^Lf^*_h(x,v)dxdv\\
            &=\sum_{j=1}^{N_v}\sum_{q=1}^{k+1}w_{j,q}\int_{0}^Lf_h^*(x,v_{j,q})dx=\sum_{j=1}^{N_v}\sum_{q=1}^{k+1}w_{j,q}\int_{0}^Lf_h^m(x,v_{j,q})dx=\int_{-v_m}^{v_m}\int_{0}^L f_h^m(x,v)dxdv.
		\end{align*}
		For \textbf{Step 2}, taking $\Psi(v)= 1$ in Eq.~(\ref{HE-SLDG}), we have
            \begin{align}\label{sl-he}
			\displaystyle \int_{V_j}f_h^{**}(x_{i,p},v)dv=\int_{V_j^{\star}}f_h^{*}(x_{i,p},v)dv,
		\end{align}
		where $V_j^{\star}=[v_{j-\frac{1}{2}}+E_h^{m+1/2}(x_{i,p})\Delta t,v_{j+\frac{1}{2}}+E_h^{m+1/2}(x_{i,p})\Delta t]$.
        \newline
            Summing Eq.~(\ref{sl-he}) from $j=1$ to $j=N_v$ and then integrating in the $x$-direction, we have
		\begin{align*}
			\displaystyle \int_{-v_m}^{v_m}\int_{0}^Lf^{**}_h(x,v)dxdv=\int_{-v_m}^{v_m}\int_{0}^L f_h^{*}(x,v)dxdv.
		\end{align*}
        It is obvious that
		\begin{align*}
			\displaystyle \int_{-v_m}^{v_m}\int_{0}^Lf^{m+1}_h(x,v)dxdv=\int_{-v_m}^{v_m}\int_{0}^L f_h^{**}(x,v)dxdv.
		\end{align*}
		In summary, we have
		\begin{align*}
			\int_{0}^L\int_{-v_m}^{v_m} f_h^{m+1}(x,v) dvdx=\int_{0}^L\int_{-v_m}^{v_m} f_h^{m}(x,v) dvdx.
		\end{align*}
	\end{proof}	
\end{propo}
    
\begin{propo}\label{momentum-prop}
	(Momentum conservation). The ECSLDG method preserves the momentum of the VA system, i.e.,
	\begin{align*}
		P^{m+1}=\displaystyle \int_{0}^L\int_{-v_m}^{v_m} f_h^{m+1}(x,v)vdvdx=\int_{0}^L\int_{-v_m}^{v_m} f_h^m(x,v)vdvdx=P^m,~~~~m=0,\cdots,M-1
	\end{align*}
	 if the following condition holds
	\begin{align*}
			\displaystyle \int_{0}^L n_h^{*}(x)E_h^{m+1/2}(x)dx=0,
	\end{align*}
	where $n_h^*(x)=\displaystyle \int_{-v_m}^{v_m}f_h^*(x,v)dv$.
		\begin{proof}
			For $\textbf{Step 1}$, it is obvious that
			\begin{align*}
				\int_{0}^L\int_{-v_m}^{v_m} f_h^*(x,v)vdvdx=\int_{0}^L\int_{-v_m}^{v_m}f_h^m(x,v)vdvdx.
			\end{align*}
			For $\textbf{Step 2}$, taking $\Psi(v)=v$ in Eq.~(\ref{HE-SLDG}), we have
            \begin{align*}
                 &\displaystyle \int_{V_j}f_h^{**}(x_{i,p},v)vdv
                 =\int_{V_{j}^{\star}} f_h^*(x_{i,p},v)(v-E_h^{m+1/2}(x_{i,p})\Delta t)dv\\
                 &=\int_{V_j^{\star}}f_h^*(x_{i,p},v)vdv-E_h^{m+1/2}(x_{i,p})\Delta t\int_{V_j^{\star}}f_h^*(x_{i,p},v)dv.
            \end{align*}
			Without loss of generality, we assume $E_h^{m+1/2}(x_{i,p})<0$ and $|E_h^{m+1/2}(x_{i,p})\Delta t|<\Delta x$. And then summing the above equality from $j=1$ to $j=N_v$,  we get
            \begin{align*}
                &\displaystyle \int_{-v_m+E_h^{m+1/2}(x_{i,p})\Delta t}^{v_m+E_h^{m+1/2}(x_{i,p})\Delta t} f_h^{*}(x_{i,p},v)vdv = \sum_{j=1}^{N_v}\int_{V_j^{\star}}f_h^{*}(x_{i,p},v)vdv
                 =\sum_{j=1}^{N_v}\sum_{l=1}^2\int_{V_{j,l}^{\star}}f_h^*(x_{i,p},v)vdv\\
                 &=\int_{-v_m}^{v_m}f_h^*(x_{i,p},v)vdv + \int_{V_{1,1}^{\star}}f_h^*(x_{i,p},v)v dv-\int_{V_{N_v}^2}f_h^*(x_{i,p},v)v dv=\int_{-v_m}^{v_m}f_h^*(x_{i,p},v)vdv,
            \end{align*}
            where $V_{N_v}^2$ represents the right part of the two segments into which the cell $V_{N_v}$ is divided by the upstream grid node $v_{N_v+\frac{1}{2}}^{\star}$.
            \newline
            It is worth that to ensure the validity of the method, we set the initial values $f_h^0(x_{i,p},v)$ of cell $V_1$ and $V_{N_v}$ to 0. Then, it is obvious that $f_h^{m+1}(x_{i,p},v)=0$ on cell $V_1$ and $V_{N_v}$.\\
			Combining the above equalities and using the derivation of Proposition \ref{mass-prop}, we have
            \begin{align*}
                \int_{-v_m}^{v_m} f_h^{**}(x_{i,p},v)vdv = \int_{-v_m}^{v_m}f_h^*(x_{i,p},v)vdv - E_h^{m+1/2}(x_{i,p})\Delta t\int_{-v_m}^{v_m}f_h^*(x_{i,p},v)dv.
            \end{align*}
            Integrating in the $x$-direction, we get 
            \begin{align*}
                &\int_{0}^L\int_{-v_m}^{v_m}f_h^{**}(x,v)vdvdx = \sum_{i=1}^{N_x}\sum_{p=1}^{k+1}w_{i,p}\int_{-v_m}^{v_m}f_h^{**}(x_{i,p},v)vdv\\
                &= \sum_{i=1}^{N_x}\sum_{p=1}^{k+1}w_{i,p}\biggl(\int_{-v_m}^{v_m}f_h^*(x_{i,p},v)vdv - E_h^{m+1/2}(x_{i,p})\Delta t\int_{-v_m}^{v_m}f_h^*(x_{i,p},v)dv\biggr)\\
                &=\int_{0}^L\int_{-v_m}^{v_m}f_h^*(x,v)vdvdx-\Delta t\int_{0}^Ln_h^*(x)E_h^{m+1/2}(x)dx,
            \end{align*}
            where $w_{i,p}$ is the the corresponding Gaussian weight of Gaussian point $x_{i,p}$ and $n_h^*(x)=\displaystyle\int_{-v_m}^{v_m}f_h^*(x,v)dv$.\\
			Similarly, for $\textbf{Step 3}$, we have
			\begin{align*}
				\int_{0}^L\int_{-v_m}^{v_m}f_h^{m+1}(x,v)vdvdx=\int_{0}^L\int_{-v_m}^{v_m}f_h^{**}(x,v)vdvdx
			\end{align*}
			In conclusion, we can obtain 
			\begin{align*}
				\int_{0}^L\int_{-v_m}^{v_m}f_h^{m+1}(x,v)vdvdx = \int_{0}^L\int_{-v_m}^{v_m}f_h^{m}(x,v)vdvdx-\Delta t\int_{0}^Ln_h^{*}(x)E_h^{m+1/2}(x)dx.
			\end{align*}
		\end{proof}
\end{propo}

\begin{propo}\label{te-prop}
	(Total energy conservation). The ECSLDG method preserves the total energy of the VA system, i.e.,
	\begin{align*}
			TE^{m+1} &=\displaystyle \frac{1}{2}\int_{0}^L\int_{-v_m}^{v_m} f_h^{m+1}(x,v)v^2dvdx+\frac{\lambda^2}{2}\int_{0}^L(E_h^{m+1}(x))^2dx\\
			&=\displaystyle \frac{1}{2}\int_{0}^L\int_{-v_m}^{v_m} f_h^{m}(x,v)v^2dvdx+\frac{\lambda^2}{2}\int_{0}^L(E_h^{m}(x))^2dx\\
			&=TE^m,~~~~m=0,\cdots,M-1.
		\end{align*}
		\begin{proof}
			For $\textbf{Step 1}$, it is obvious that
			\begin{align*}
				\frac{1}{2}\int_{0}^L\int_{-v_m}^{v_m} f_h^*(x,v)v^2dvdx=\frac{1}{2}\int_{0}^L\int_{-v_m}^{v_m}f_h^m(x,v)v^2dvdx.
			\end{align*}
			For $\textbf{Step 2}$, taking $\Psi(v)=v^2$ in Eq.~(\ref{HE-SLDG}) and using the derivation of Proposition \ref{momentum-prop}, we have
			\begin{align*}
                &\frac{1}{2}\int_{V_j}f_h^{**}(x_{i,p},v)v^2dv=\frac{1}{2}\int_{V_j^{\star}}f_h^*(x_{i,p},v)(v-E_h^{m+1/2}(x_{i,p})\Delta t)^2dv\\
                &=\frac{1}{2}\int_{V_j^{\star}}f_h^*(x_{i,p},v)v^2dv\\
                &~~~~~-\frac{\Delta t}{2}E_h^{m+1/2}(x_{i,p})\biggl(\int_{V_j^\star}f_h^{\star}(x_{i,p},v)vdv+\int_{V_j^{\star}}f_h^*(x_{i,p},v)(v-E_h^{m+1/2}(x_{i,p})\Delta t)dv\biggr)\\
				&=\frac{1}{2}\int_{V_j^{\star}}f_h^{*}(x_{i,p},v)v^2dv -\frac{\Delta t}{2}E_h^{m+1/2}(x_{i,p})\biggl(\int_{V_j^{\star}}f_h^*(x_{i,p},v)vdv+\int_{V_j}f_h^{**}(x_{i,p},v)vdv\biggr).
			\end{align*}
			 Summing the above equality from $j=1$ to $j=N_v$ and using the assumption in the derivation of Proposition \ref{momentum-prop}, we have
             \begin{align*}
                 &\frac{1}{2}\int_{-v_m+E_h^{m+1/2}(x_{i,p})\Delta t}^{v_m+E_h^{m+1/2}(x_{i,p})\Delta t} f_h^{*}(x_{i,p},v)v^2dv = \frac{1}{2}\sum_{j=1}^{N_v}\int_{V_j^{\star}}f_h^{*}(x_{i,p},v)v^2dv=\sum_{j=1}^{N_v}\sum_{l=1}^2\int_{V_{j,l}^{\star}}f_h^*(x_{i,p},v)v^2dv\\
                 &=\frac{1}{2}\int_{-v_m}^{v_m}f_h^{*}(x_{i,p},v)v^2dv.
             \end{align*}
             Combining the above equalities, we have
             \begin{align*}
                 \frac{1}{2}\int_{-v_m}^{v_m}f_h^{**}(x_{i,p},v)v^2dv = \frac{1}{2}\int_{-v_m}^{v_m} f_h^*(x_{i,p},v)v^2dv +\Delta tE_h^{m+1/2}(x_{i,p})J_h^{m+1/2}(x_{i,p}),
             \end{align*}
			 where $J_h^{m+1/2}(x_{i,p})=(J_h^*(x_{i,p})+J_h^{**}(x_{i,p}))/2$.
			 
			 Integrating the above equation in the $x$-direction, we obtain
			\begin{align*}
				&\frac{1}{2}\int_{0}^L\int_{-v_m}^{v_m}f_h^{**}(x,v)v^2dvdx=\frac{1}{2}\sum_{i=1}^{N_x}\sum_{p=1}^{k+1}w_{i,p}\int_{-v_m}^{v_m}f_h^{**}(x_{i,p},v)v^2dv\\
                &=\sum_{i=1}^{N_x}\sum_{p=1}^{k+1}w_{i,p}\biggl(\frac{1}{2}\int_{-v_m}^{v_m}f_h^*(x_{i,p},v)v^2dv+\Delta tE_h^{m+1/2}(x_{i,p})J_h^{m+1/2}(x_{i,p})\biggr)\\
                &=\frac{1}{2}\int_{0}^L\int_{-v_m}^{v_m}f_h^*(x,v)v^2dvdx+\Delta t\int_{0}^LE_h^{m+1/2}(x)J_h^{m+1/2}(x)dv.
			\end{align*}			
			 Multiplying both sides of Eq.~(\ref{eq:12}) by $E_h^{m+1/2}(x_{i,p})$ and 
                    then integrating in the $x$-direction, we have
			\begin{align*}
				&\displaystyle \frac{\lambda^2}{2}\int_{0}^L (E_h^{m+1}(x))^2dx =\frac{\lambda^2}{2}\sum_{i=1}^{N_x}w_{i,p}(E_h^{m+1}(x_{i,p}))^2\\ &=\sum_{i=1}^{N_x}w_{i,p}\biggl(\frac{\lambda^2}{2} (E_h^m(x_{i,p}))^2 - \Delta tE_h^{m+1/2}(x_{i,p})J_h^{m+1/2}(x_{i,p}) \biggr)\\
                &=\frac{\lambda^2}{2}\int_{0}^L (E_h^{m}(x))^2dx -\Delta t
				\int_{0}^L E_h^{m+1/2}(x)J_h^{m+1/2}(x)dx.
			\end{align*}
			For $\textbf{Step 3}$, it is obvious that
			\begin{align*}
				\frac{1}{2}\int_{0}^L\int_{-v_m}^{v_m} f_h^{m+1}(x,v)v^2dvdx=\frac{1}{2}\int_{0}^L\int_{-v_m}^{v_m}f_h^{**}(x,v)v^2dvdx.
			\end{align*}
			To sum up, we have
			\begin{align*}
				&\displaystyle \frac{1}{2}\int_{0}^L\int_{-v_m}^{v_m} f_h^{m+1}(x,v)v^2dvdx+\frac{\lambda^2}{2}\int_{0}^L(E_h^{m+1}(x))^2dx\\
				&=\displaystyle \frac{1}{2}\int_{0}^L\int_{-v_m}^{v_m} f_h^{m}(x,v)v^2dvdx+\frac{\lambda^2}{2}\int_{0}^L(E_h^{m}(x))^2dx.
			\end{align*}
		\end{proof}
\end{propo}

\begin{propo}\label{L2-prop}
	($L^2$ stability). The ECSLDG method satisfy
	\begin{align*}
		\displaystyle \int_{0}^L\int_{-v_m}^{v_m} |f_h^{m+1}(x,v)|^2dvdx \le \int_{0}^L\int_{-v_m}^{v_m} |f_h^{m}(x,v)|^2dvdx.
	\end{align*}
	\begin{proof}
		For $\textbf{Step 1}$, taking $\Psi(x)=f_h^*(x,v_{j,q})$ in Eq.~(\ref{Hf-SLDG}), we have 
        \begin{align*}
            &\displaystyle \int_{X_i}|f_h^*(x,v_{j,q})|^2dx=\int_{X_i^{\star}} f_h^m(x,v_{j,q})f_h^*(x+v_{j,q}\Delta t/2,v_{j,q})dx\\
            &\le \frac{1}{2}\biggl(\int_{X_i}|f_h^*(x,v_{j,q})|^2dx+\int_{X_i^{\star}}|f_h^m(x,v_{j,q})|^2dx\biggr).
        \end{align*}
		Thus, summing the above inequality from $i=1$ to $i=N_x$ and integrating over the velocity space, we have
		\begin{align*}
			\int_{-v_m}^{v_m}\int_{0}^L|f_h^*(x,v_{j,q})|^2dxdv \le \int_{-v_m}^{v_m}\int_{0}^L|f_h^m(x,v_{j,q})|^2dxdv.
		\end{align*}
		Similarly, we can obtain the following conclusion in order
		\begin{align*}
			&\int_{-v_m}^{v_m}\int_{0}^L|f_h^{**}(x,v_{j,q})|^2dxdv \le \int_{-v_m}^{v_m}\int_{0}^L|f_h^*(x,v_{j,q})|^2dxdv,\\
			&\int_{-v_m}^{v_m}\int_{0}^L|f_h^{m+1}(x,v_{j,q})|^2dxdv \le \int_{-v_m}^{v_m}\int_{0}^L|f_h^{**}(x,v_{j,q})|^2dxdv.
		\end{align*}
		Therefore, we have
		\begin{align*}
			\displaystyle \int_{0}^L\int_{-v_m}^{v_m} |f_h^{m+1}(x,v)|^2dvdx \le \int_{0}^L\int_{-v_m}^{v_m} |f_h^{m}(x,v)|^2dvdx.
		\end{align*}
	\end{proof}
\end{propo}

\subsection{High order ECSLDG for VA \label{hgspva}}

The development of high-order time integration schemes is critical for enhancing the accuracy and robustness of kinetic plasma simulations. In this work, we systematically evaluate and integrate advanced splitting methods into the ECSLDG framework to achieve temporal accuracy beyond second order. A key challenge lies in identifying schemes that balance high-order precision with computational efficiency and physical fidelity, as not all high-order methods can perform as well as expected in practical applications\cite{EGC2006}. 
To begin with, a review of some work on splitting methods is provided.

Consider the following differential equation
\begin{align}\label{des}
	\dot{x} = X(x),
\end{align}
where $x\in \mathbb{R}^I$ and $X$ is a vector field on $\mathbb{R}^I$. Applying the constant formula to Eq.~(\ref{des}) yields the result that $x(t)=\exp(tX)(x(0))$. The splitting methods apply when $X= \displaystyle\sum \limits^I_{i=1} X_i$ and each of $X_i$ can be exactly integrated or approximated.

According to \cite{Trotter1959}, a first-order integrator can be obtained
\begin{align}\label{1or-splitting}
    \varphi(\Delta t)=\exp(\Delta tX_1)\cdots \exp(\Delta t X_I).
\end{align}
Applied to the basic composition Eq.~(\ref{1or-splitting}), a second-order symmetric integrator
\begin{align}\label{2or-splitting}
    SS(\Delta t)=\exp(\frac{1}{2}\Delta tX_1)\cdots \exp(\frac{1}{2}\Delta tX_{I-1})\exp(\Delta tX_I)\exp(\frac{1}{2}\Delta tX_{I-1})\cdots \exp(\frac{1}{2}\Delta tX_1)
\end{align}
 which is called Type S with symmetric stages, can be obtained by using the Strang compisition\cite{SG1968}
\begin{align*}
    SS(\Delta t)=\varphi(\Delta t/2) \varphi^{-1}(-\Delta t/2),
\end{align*}
where $\varphi^{-1}(-\Delta t)$ denotes the adjoint of $\varphi(\Delta t)$ and 
\begin{align*}
    \varphi^{-1}(-\Delta t)=\exp(\Delta tX_I)\cdots \exp(\Delta t X_1).
\end{align*}

To better maintain physical properties, a common practice is to improve the accuracy of the splitting scheme. Refer to \cite{McLachlan_Quispel_2002}, a high order splitting method can be derived
\begin{align}\label{eqss}
    SS_{2m+1}(\Delta t)= (SS(\alpha \Delta t))^m(SS(\beta \Delta t))(SS(\alpha \Delta t))^m,
\end{align}
where $\alpha = 1/(2m-(2m)^{1/3})$, $\beta = 1-2m\alpha$. 

Therefore, we can derive a fourth-order ECSLDG method for the VA system by using splitting method $SS_{2m+1}$ as follows
\begin{align}\label{4or-SS-SLDG}
    \mathcal{H}^{SLDG}_{SS_{2m+1}}(\Delta t) = (\mathcal{H}^{SLDG}_{SS}(\alpha\Delta t))^m(\mathcal{H}^{SLDG}_{SS}(\beta \Delta t))(\mathcal{H}^{SLDG}_{SS}(\alpha \Delta t))^m.
\end{align}
\begin{remark}
    One time step $\Delta t$ of executing $\mathcal{H}^{SLDG}_{SS_{2m+1}}$ follows the specific procedure: first, $\mathcal{H}^{SLDG}_{SS}$ method is executed $m$ times with a time step of $\alpha \Delta t$; then, $\mathcal{H}^{SLDG}_{SS}$ method is executed $1$ time with a time step of $\beta\Delta t$; finally, $\mathcal{H}^{SLDG}_{SS}$ method is executed $m$ times again with a time step of $\alpha \Delta t$.
\end{remark}
An alternative fourth-order scheme $10Lie$(4th,10Lie) is as below\cite{EGC2006} 
\begin{align}\label{4or-10lie}
    10Lie(\Delta t) = \varphi(a_5\Delta t)\varphi^{-1}(-b_5\Delta t)\cdots\varphi(a_1\Delta t)\varphi^{-1}(-b_1\Delta t),
\end{align}
 where 
 \begin{align*}
     a_1=b_5=\displaystyle \frac{146+5\sqrt{19}}{540}, && a_2=b_4=\displaystyle \frac{-2+10\sqrt{19}}{135},\\
     a_3=b_3=\frac{1}{5}, && a_4=b_2=\displaystyle \frac{-23-20\sqrt{19}}{270},\\
     a_5=b_1= \displaystyle \frac{14-\sqrt{19}}{108}.&&~
 \end{align*}


Due to similarity with Eq.~(\ref{2or-SLDG}), we can obtain a first-order ECSLDG method for the VA system by using splitting method Eq.~(\ref{1or-splitting}) as follows
\begin{align}
    \mathcal{H}^{SLDG}_{\varphi}(\Delta t)=\mathcal{H}_f^{SLDG}(\Delta t)\mathcal{H}_E^{SLDG}(\Delta t).
\end{align}
Then, we can derive another fourth-order ECSLDG method for the VA system by using splitting method $10Lie$ as follows
\begin{align}\label{10lie}
    \mathcal{H}^{SLDG}_{10Lie}(\Delta t) = \mathcal{H}^{SLDG}_{\varphi}(a_5\Delta t)\mathcal{H}^{SLDG}_{\varphi^{-1}}(-b_5\Delta t)\cdots \mathcal{H}^{SLDG}_{\varphi}(a_1\Delta t)\mathcal{H}^{SLDG}_{\varphi^{-1}}(-b_1\Delta t). 
\end{align}
\begin{remark}
    Since each subsystem is solved either exactly or based on an energy conserving discretization, the high-order ECSLDG method still satisfies Proposition \ref{mass-prop}-Proposition \ref{L2-prop}.
\end{remark}
\begin{remark}
Among all the forth-order splitting schemes mentioned above, $SS_3$ has lowest computational cost but large truncation error affecting stability\cite{McLachlan_Quispel_2002}. Increasing $m$ improves stability until $m\approx19$ but raises computational cost. Based on our simulation tests, we found that the fourth-order $10Lie$ scheme provides the best balance between stability and efficiency.
\end{remark}
\begin{remark}
The proposed method cannot guarantee Gauss's law. However, numerical tests show that the ECSLDG method, using the fourth-order $10Lie$ scheme significantly reduces Gauss's law residuals compared to the second-order scheme, especially when a large time step is used.
\end{remark}

\section{Numerical results}
\label{section:numerical}

We evaluate the performance of the proposed ECSLDG method for the VA system using four 1D1V numerical experiments. The assessment is carried out from the following perspectives:

\textbf{1. Spatial and temporal accuracy.} In the simulation of weak Landau damping in section \ref{test1}, we evaluate the temporal and spatial accuracy of the ECSLDG method combined with different high-order splitting schemes discussed in section \ref{hgspva}. The spatial accuracy is verified using the time-reversibility property of the VA system, while the temporal error is assessed following the approach in \cite{LCLC2021}. For simplicity, we report results using the fourth-order splitting scheme $SS_3$, which offers the lowest computational cost among the forth-order methods.

\textbf{2. Conservation properties.}
We compare conservation properties of the ECSLDG method with those of a non-energy-conserving, conservative semi-Lagrangian discontinuous Galerkin method (CSLDG-AE). The CSLDG-AE scheme shares the same computational cycle as ECSLDG, differing only in the electric field solver, where the electric field is updated by
\begin{align*}
    \mathit{E}^{m+1} = \mathit{E}^m - \frac{\Delta t}{\lambda^2}\mathit{J}^m,
\end{align*}
In the weak Landau damping case, this comparison illustrates the ECSLDG method conserves total energy. In the two-stream instability test in section \ref{test3} or \ref{test4}, we  examine how the conservation properties depend on temporal or spatial resolution. The following diagnostics are used:

\begin{itemize}
    \item Electric energy $E_{\bm{E}}$:
\begin{align*}
    E_{\bm{E}}(t_m) = \frac{\lambda^2}{2} \sum_{i=1}^{N_x}\sum_{p=1}^{k+1}w_{i,p} (E_h^m(x_{i,p}))^2
\end{align*}

\item Mass $L_1$:
\begin{align*}
    L_1(t_m) = \displaystyle \sum_{i=1}^{N_x}\sum_{p=1}^{k+1}\sum_{j=1}^{N_v}\sum_{q=1}^{k+1} w_{i,p}w_{j,q}f_h^{m}(x_{i,p},v_{j,q})
\end{align*}

\item Momentum $P$:
\begin{align*}
    P(t_m) = \displaystyle \sum_{i=1}^{N_x}\sum_{p=1}^{k+1}\sum_{j=1}^{N_v}\sum_{q=1}^{k+1} w_{i,p}w_{j,q}f_h^{m}(x_{i,p},v_{j,q})v_{j,q}.
\end{align*}

\item Total energy $E_{total}$:
\begin{align*}
    E_{total}(t_m) = \frac{1}{2}\displaystyle \sum_{i=1}^{N_x}\sum_{p=1}^{k+1}\sum_{j=1}^{N_v}\sum_{q=1}^{k+1} w_{i,p}w_{j,q}f_h^{m}(x_{i,p},v_{j,q})v_{j,q}^2 + \frac{\lambda^2}{2} \sum_{i=1}^{N_x}\sum_{p=1}^{k+1}w_{i,p} (E_h^m(x_{i,p}))^2.
\end{align*}
\end{itemize}


\textbf{3. Temporal accuracy effects.}  In the two-stream instability simulation in section \ref{test3}, we investigate how temporal accuracy influences physical quantities by examining macroscopic variables such as the electric field $E$ and number density $n$, as well as the microscopic distribution function $f$, using both second-order and fourth-order time integration schemes.

\textbf{4. Gauss’s law preservation.} In the simulations of strong Landau damping in section \ref{test2},  we evaluate enforcement of Gauss’s law  using the following formula:
\begin{align*}
    ||\nabla \cdot E-\rho||_2 = \displaystyle \biggl(\int_0^L (\nabla \cdot E-\rho)^2dx \biggr)^{1/2}.
\end{align*}

\textbf{5. Unconditional stability.} In the simulation of two-stream instability in section \ref{test4}, we test the stability of the ECSLDG method under small Debye length $\lambda$ to verify that its performance is not limited by this scale, highlighting the unconditional stability of the scheme.


    

Unless otherwise stated, we adopt the following time step
\[
\Delta t= \frac{CFL}{\frac{v_m}{\Delta x}+\frac{max|E|}{\Delta v}},
\]
and use quadratic piecewise discontinuous polynomials for spatial discretization in all simulations.
\begin{table}[http]
\centering
\caption{Weak Landau damping: $T=0.5$ and $\rm{CFL}=0.1$. The spatial errors and convergence orders of the ECSLDG method.}
\label{1d1v-w-space}
\begin{tabular}{llllllllll}
\hline
Mesh            &  & $||f||_{L^1}$ & Order &  & $||f||_{L^2}$ & Order &  & $||E||_{L^2}$ & Order \\ \hline
                &  & \multicolumn{8}{l}{$P^1$}                                                 \\ \cline{3-10}
$40\times 40$   &  & 1.30E-02      &       &  & 1.90E-03      &       &  & 4.82E-06      &       \\
$60\times 60$   &  & 5.32E-03      & 2.20  &  & 7.83E-04      & 2.19  &  & 2.15E-06      & 1.99  \\
$80\times 80$   &  & 2.93E-03      & 2.07  &  & 4.37E-04      & 2.03  &  & 1.24E-06      & 1.91  \\
$100\times 100$ &  & 1.96E-03      & 1.80  &  & 2.94E-04      & 1.78  &  & 7.88E-07      & 2.04  \\ \hline
                &  & \multicolumn{8}{l}{$P^2$}                                                 \\ \cline{3-10}
$40\times 40$   &  & 3.30E-03      &       &  & 5.55E-04      &       &  & 5.64E-07      &       \\
$60\times 60$   &  & 1.23E-03      & 2.44  &  & 2.09E-04      & 2.40  &  & 2.08E-07      & 2.46  \\
$80\times 80$   &  & 5.88E-04      & 2.55  &  & 1.01E-04      & 2.53  &  & 9.72E-08      & 2.65  \\
$100\times 100$ &  & 3.28E-04      & 2.62  &  & 5.62E-05      & 2.63  &  & 5.35E-08      & 2.68  \\ \hline
                &  & \multicolumn{8}{l}{$P^3$}                                                 \\ \cline{3-10}
$40\times 40$   &  & 7.79E-05      &       &  & 1.26E-05      &       &  & 7.85E-10      &       \\
$60\times 60$   &  & 1.26E-05      & 4.50  &  & 2.03E-06      & 4.50  &  & 1.48E-10      & 4.11  \\
$80\times 80$   &  & 3.63E-06      & 4.32  &  & 5.82E-07      & 4.34  &  & 5.58E-11      & 3.39  \\
$100\times 100$ &  & 1.42E-06      & 4.20  &  & 2.28E-07      & 4.19  &  & 2.23E-11      & 4.12  \\ \hline
\end{tabular}
\end{table}
\subsection{Weak Landau damping \label{test1}}
In this section, weak Landau damping is presented to investigate the accuracy and conservative properties of the proposed ECSLDG method. The initial distribution function is given by
\begin{align}\label{wld}
f_0(x,v)=\displaystyle \frac{1}{\sqrt{2\pi}}(1+\alpha \cos(\kappa x))\exp(-\frac{v^2}{2}),
\end{align}
where $\alpha=0.01$, $\kappa=0.5$. The computational domain is set to 
$[0,{2\pi}/{\kappa}]\times [-v_m,v_m]$.
Unless otherwise stated, all simulations are performed using the ECSLDG method in combination with the fourth-order time-splitting scheme $SS_3$.


In Table \ref{1d1v-w-space}, we present the $L^1$ and $L^2$ errors of the velocity distribution function $f$, along with the $L^2$ error of the electric field $E$ for the ECSLDG method at varying spatial resolutions.
The spatial errors are computed using the time-reversibility technique described in \cite{CAI2018529}.
Here, we set $T = 0.5$, $\rm{CFL}=0.1$, $v_m=10$ and use piecewise discontinuous polynomial spaces of degree $k = 1, 2, 3$ for spatial discretization. 
It can be clearly observed that the spatial convergence order of the ECSLDG method is consistent with the theoretical expectation. 
\begin{figure}[H]
	\begin{center}
	\begin{minipage}{0.5\linewidth}			             
        \centerline{\includegraphics[width=1\linewidth]{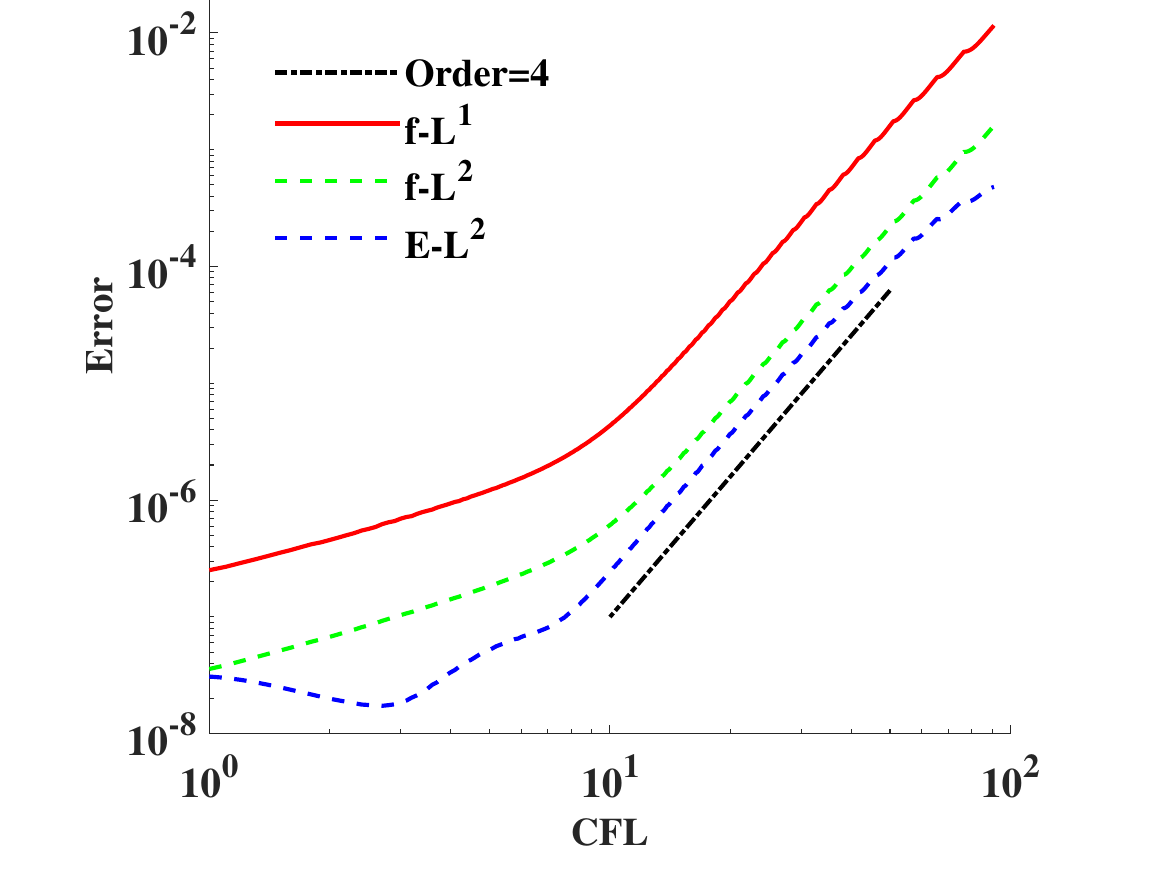}}
	\end{minipage}
	\end{center}
	\caption{Weak Landau damping: $N_x=N_v=128$ and $T=5$. $L^1$ and $L^2$ errors of the velocity distribution function $f$, and $L^2$ errors of the electric field $E$ for the ECSLDG method under different CFL numbers.}
	\label{temporal_order_w}
\end{figure}


To verify temporal convergence, we generate a reference solution using a small CFL number of 0.01 \cite{LCLC2021} on a fixed $128\times128$ grid. As shown in Fig.~\ref{temporal_order_w}, the ECSLDG method achieves fourth‐order accuracy in time, in agreement with theoretical predictions. All fourth‐order splitting methods introduced in Section \ref{hgspva} display similar convergence behavior; here, for brevity, we report only the results obtained with the $SS_3$ scheme.

\begin{figure}[H]
	\begin{center}
		\begin{minipage}{0.49\linewidth}
			\centerline{\includegraphics[width=1\linewidth]{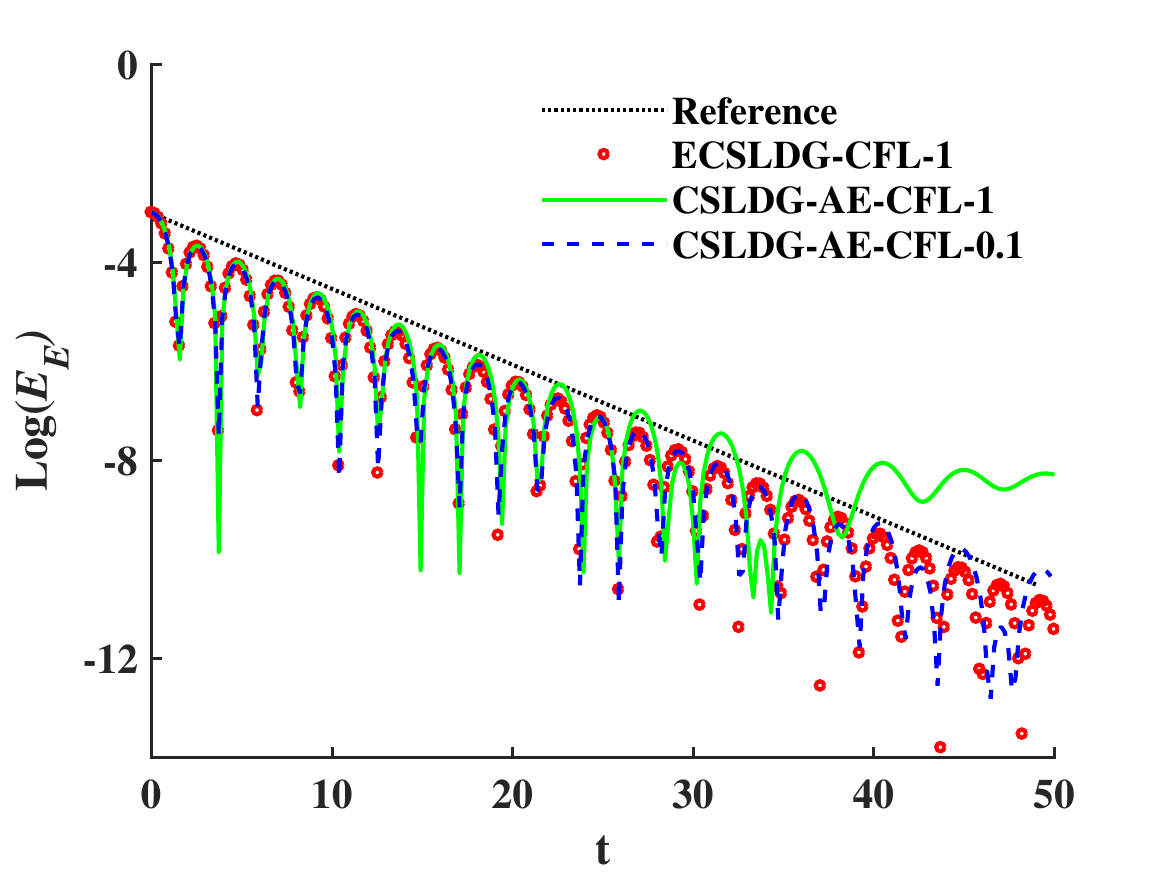}}
			\centerline{\textbf{(a)}}
		\end{minipage}
		\hfill
		\begin{minipage}{0.49\linewidth}
			\centerline{\includegraphics[width=1\linewidth]{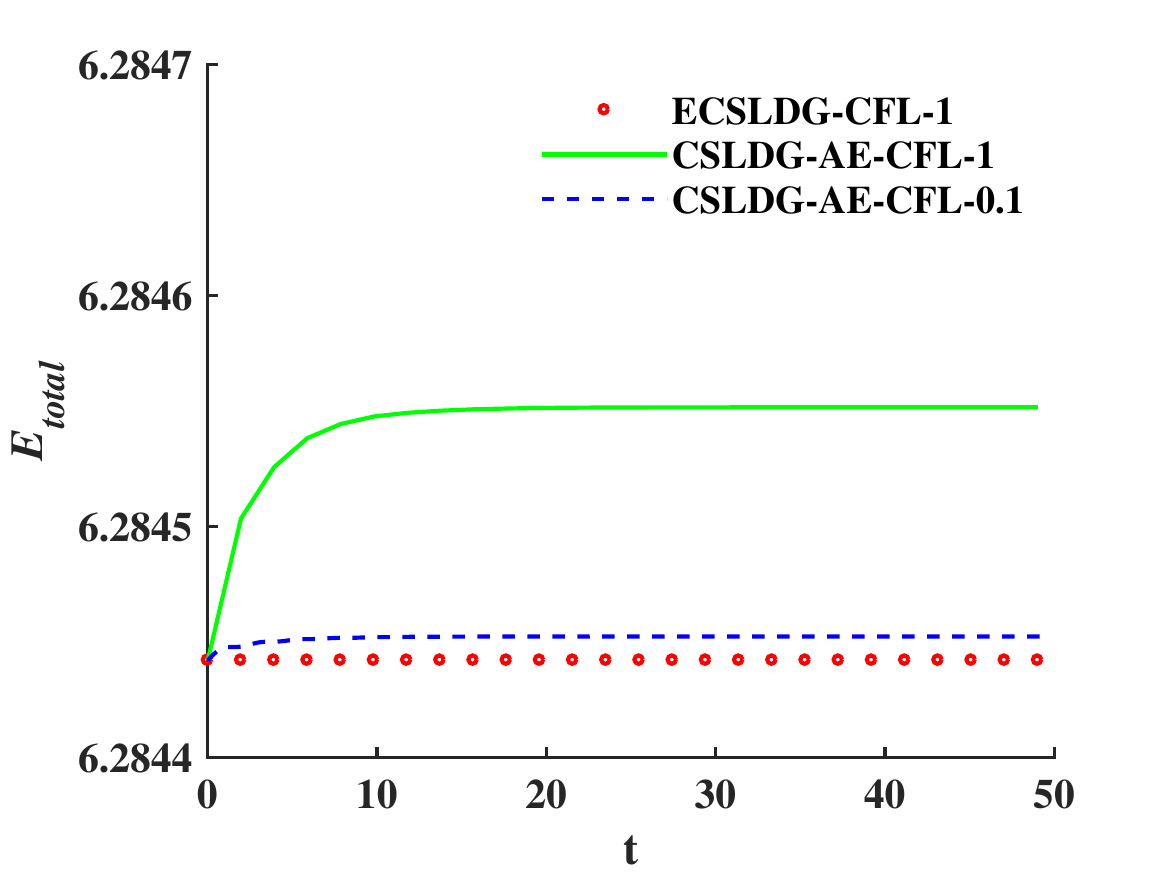}}
			\centerline{\textbf{(b)}}
		\end{minipage}
	\end{center}
	\caption{Weak Landau damping: $N_x=N_v=128$ and $T=50$. Time evolution of the electric energy (a) and the total energy (b).}
	\label{electric-w}
\end{figure}

To further evaluate energy conservation, we compare the ECSLDG method with a non-energy-conserving variant, CSLDG-AE. 
In Fig.~\ref{electric-w}, we present the time evolution of the electric energy $E_{\bm{E}}$ and the total energy $E_{total}$ using  $128\times 128$ elements and quadratic polynomial approximation. 
As shown in the Fig.~\ref{electric-w}, the decay rate simulated by the CSLDG-AE method deviates  from the theoretical value at $T=50$ when using a CFL number of 1, resulting in a 0.01\% error in the total energy. Even with a smaller CFL number of 0.1, the CSLDG-AE method still exhibits the energy drift (approximately $ 10^{-5}$). Encouragingly, the ECSLDG method with $\rm{CFL} = 1$ accurately captures the theoretical decay rate $\gamma=-0.1533$\cite{Filbet2001} and conserves the total energy.

The above arguments demonstrate that the ECSLDG method, combined with a high-order splitting scheme, achieves high accuracy in both space and time, while conserving the total energy.

\subsection{Strong Landau damping \label{test2}}
In this section, we employ the ECSLDG method to simulate strong Landau damping, aiming to demonstrate the necessity of using high-order temporal schemes. Specifically, we assess the accuracy and conservation properties of the ECSLDG method when combined with different high-order splitting schemes. Based on the parameter settings for the initial velocity distribution function in Section~\ref{test1}, the perturbation parameter is set to $\alpha=0.5$. Unless otherwise specified, simulations are performed with $v_m=10$ and a grid of $128 \times 128$.

\begin{figure}[htbp]
	\begin{center}
		\begin{minipage}{0.49\linewidth}
			\centerline{\includegraphics[width=1\linewidth]{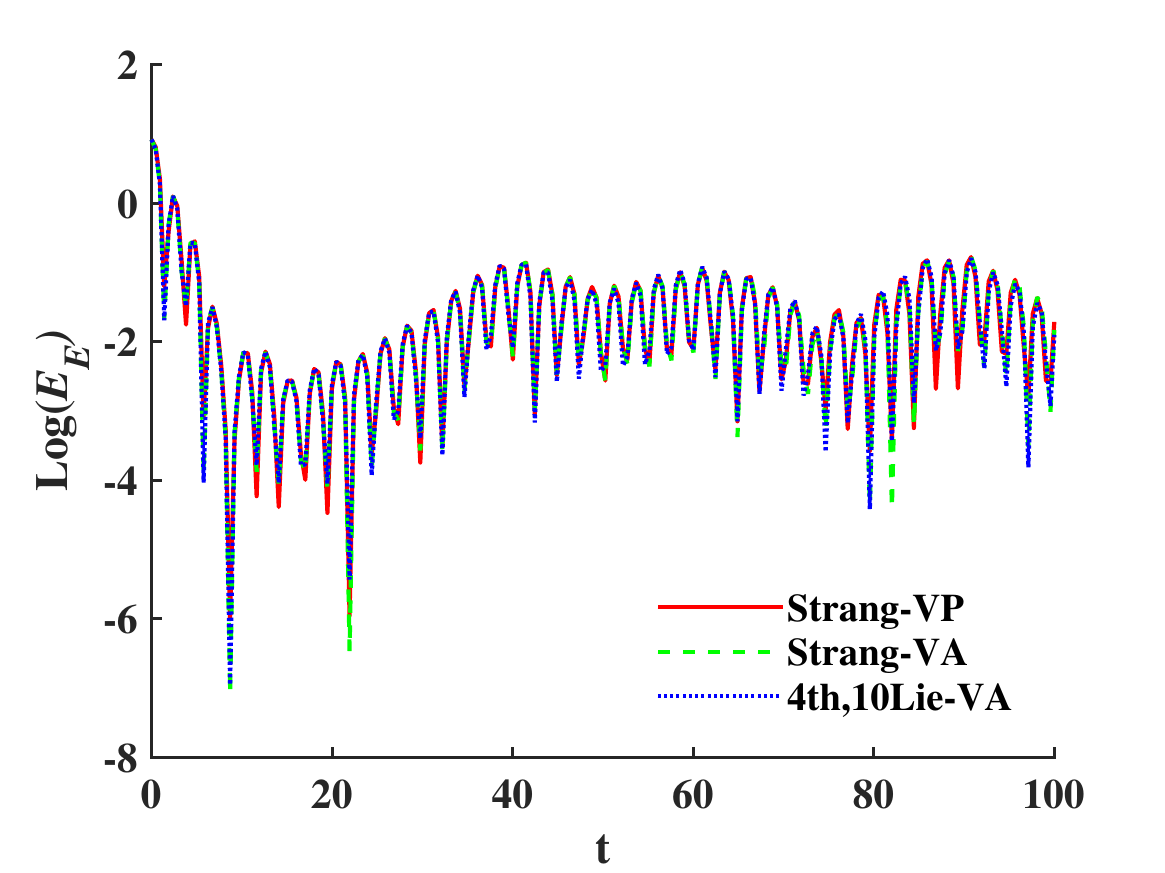}}
			\centerline{\textbf{(a)}}
		\end{minipage}
		\hfill
		\begin{minipage}{0.49\linewidth}
			\centerline{\includegraphics[width=1\linewidth]{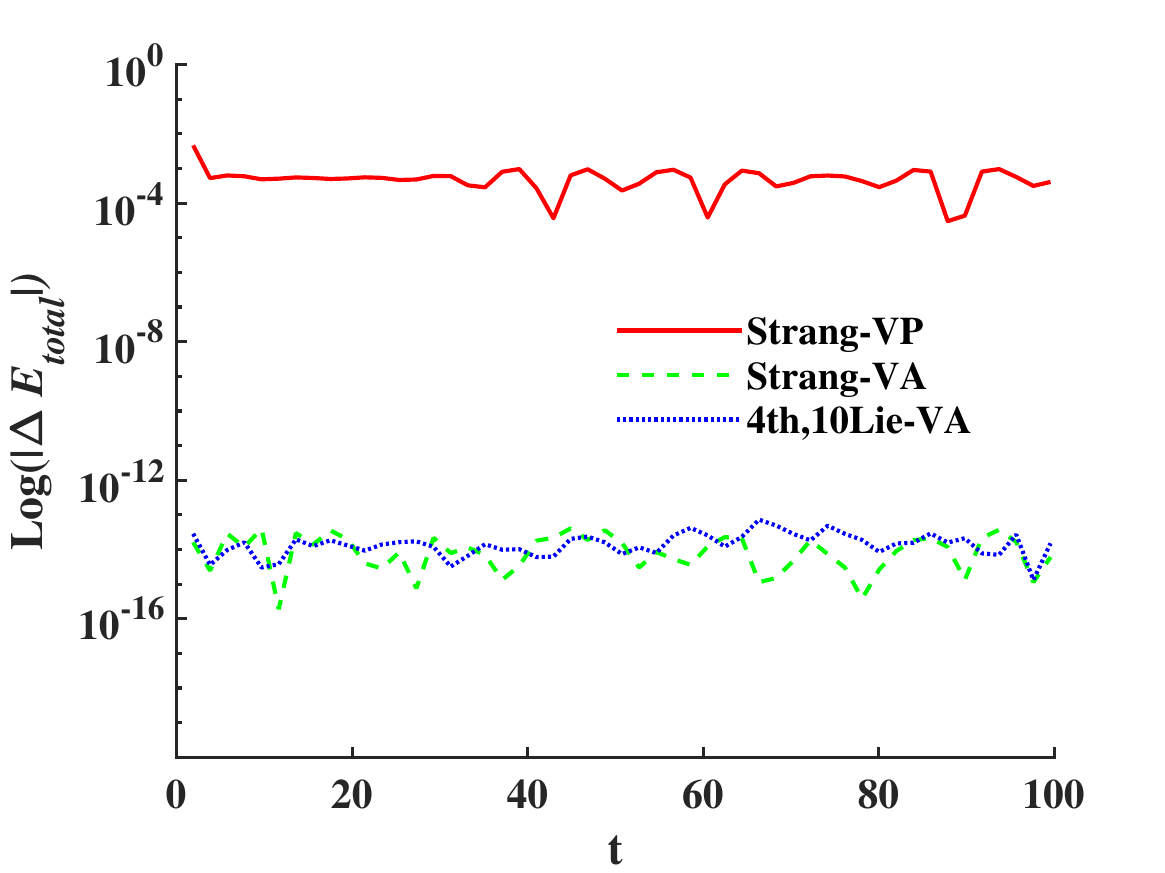}}
			\centerline{\textbf{(b)}}
		\end{minipage}
	\end{center}
	\caption{Strong Landau damping: $N_x=N_v=128$ and $\text{CFL}=10$. Time evolution of the electric energy (a) and the absolute value of the relative deviation of the total energy (b).}\label{electric-s}
	\label{electric-AP-s}
\end{figure}


In Fig.~\ref{electric-AP-s}, we plot the time evolution of the electric energy $E_{\bm{E}}$ and the total energy deviation $\Delta E_{total}=[E_{total}(t)-E_{total}(0)]/E_{total}(0)$ at $\text{CFL}=10$. In Fig.~\ref{electric-AP-s}(a), the electric‐energy curve computed by the CSLDG method for the Vlasov–Poisson system serves as the reference benchmark.  The results show that the decay or growth rate of electric energy at different splitting schemes are consistent with reference VP results. As shown in Fig.~\ref{electric-AP-s}(b), the ECSLDG method conerves the total energy when combined with either the second-order Strang splitting scheme or the fourth-order $10Lie$ scheme.



\begin{figure}[htbp]
	\begin{center}
		\begin{minipage}{0.49\linewidth}
			\centerline{\includegraphics[width=1\linewidth]{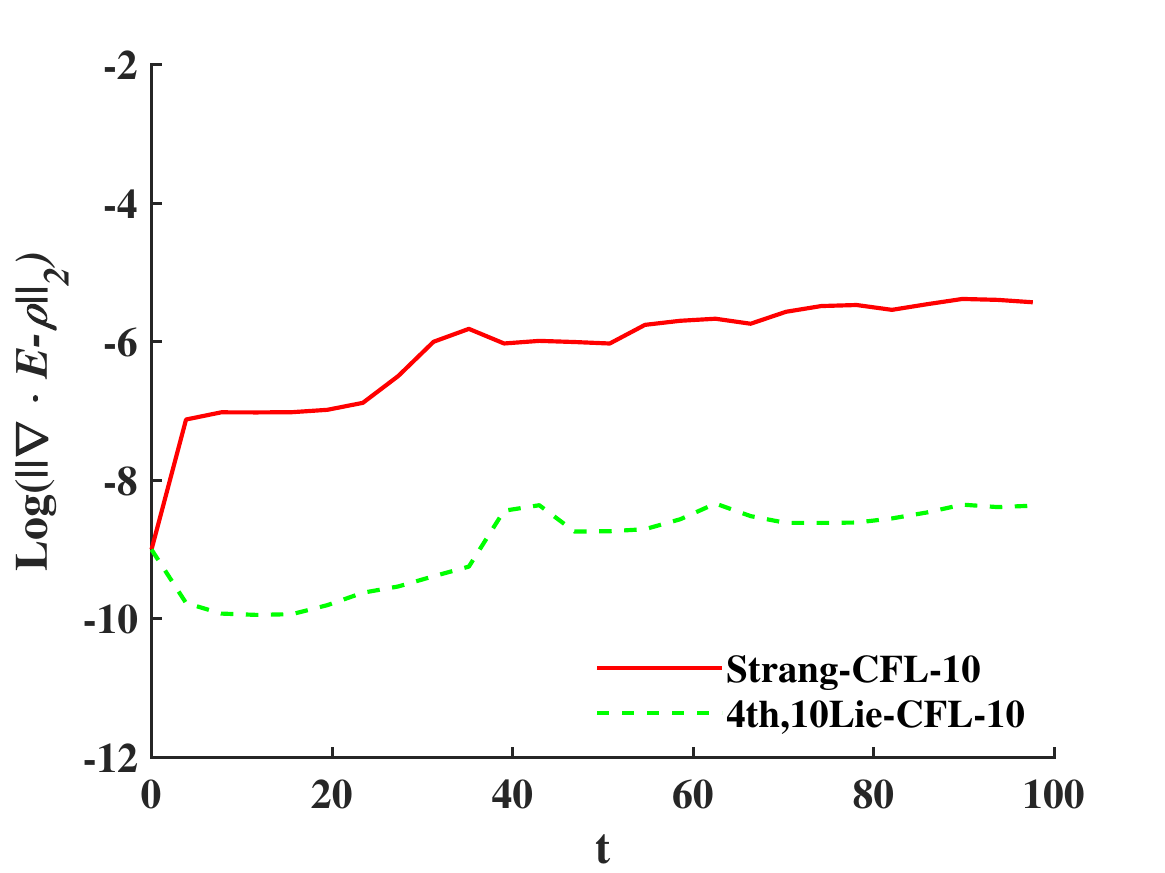}}
			\centerline{\textbf{(a)}}
		\end{minipage}
		\hfill
		\begin{minipage}{0.49\linewidth}
			\centerline{\includegraphics[width=1\linewidth]{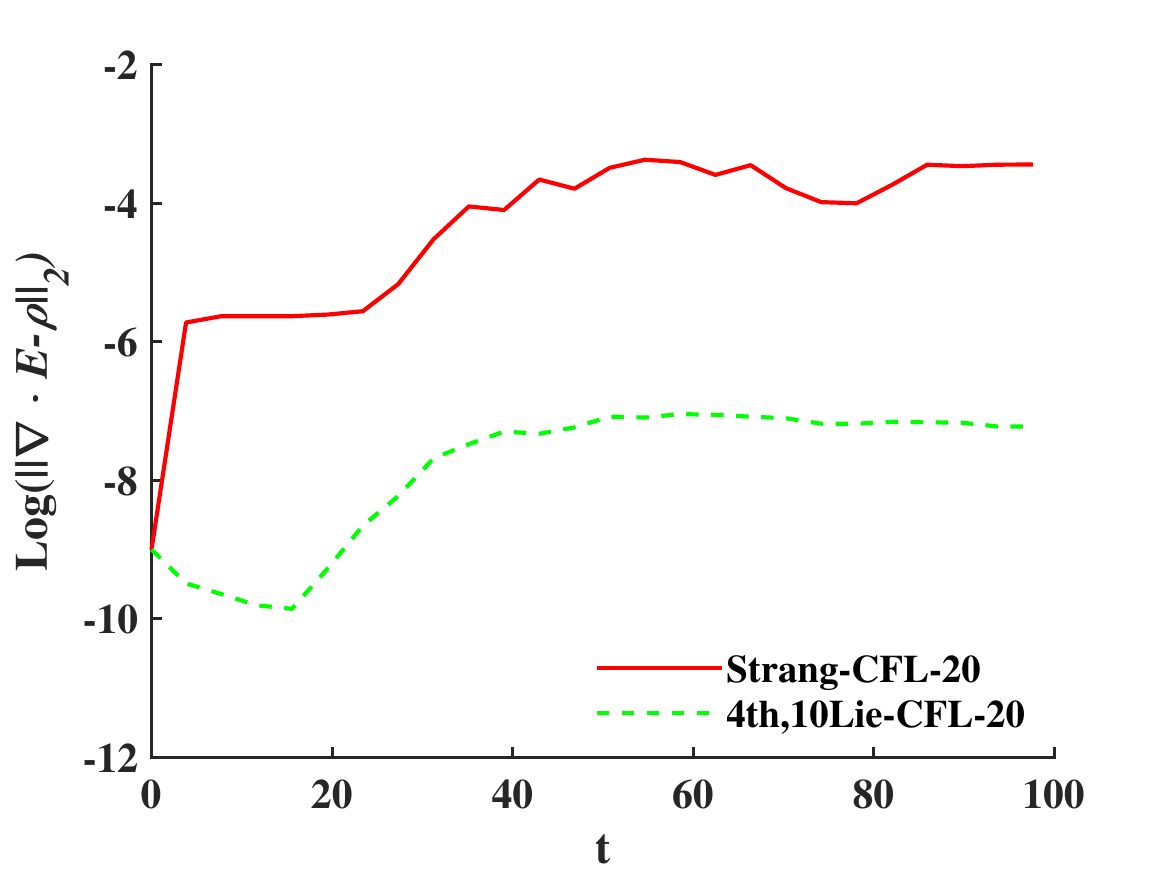}}
			\centerline{\textbf{(b)}}
		\end{minipage}
	\end{center}
	\caption{Strong Landau damping: $N_x=N_v=128$. Temporal evolution of Gauss law's residuals for the ECSLDG method with different splitting schemes at $\text{CFL}=10$ (a) and $\text{CFL}=20$ (b). }
	\label{gauss-order-s}
\end{figure}

Although higher temporal accuracy only slightly influences electric‐energy evolution and total‐energy conservation, it markedly improves the Gauss’s law residuals.  In Fig.~\ref{gauss-order-s}, we plot the time evolution of Gauss's law residuals for the second-order Strang splitting scheme and the fourth-order $10Lie$ splitting scheme at $\text{CFL}=10$ and $\text{CFL}=20$. At $\text{CFL}=10$, the residuals associated with the fourth-order scheme decrease markedly compared to the second-order scheme. 
Even with the larger time step $\text{CFL}=20$, the residuals for the fourth-order scheme remain similar to those at $\text{CFL}=10$, whereas those for the second-order scheme exhibit a clear growth trend.  This demonstrates that high-order splitting schemes substantially enhance the long-term compliance of the ECSLDG method with Gauss's law when using large time steps.


\begin{figure}[htbp]
	\begin{center}
		\begin{minipage}{0.49\linewidth}
			\centerline{\includegraphics[width=1\linewidth]{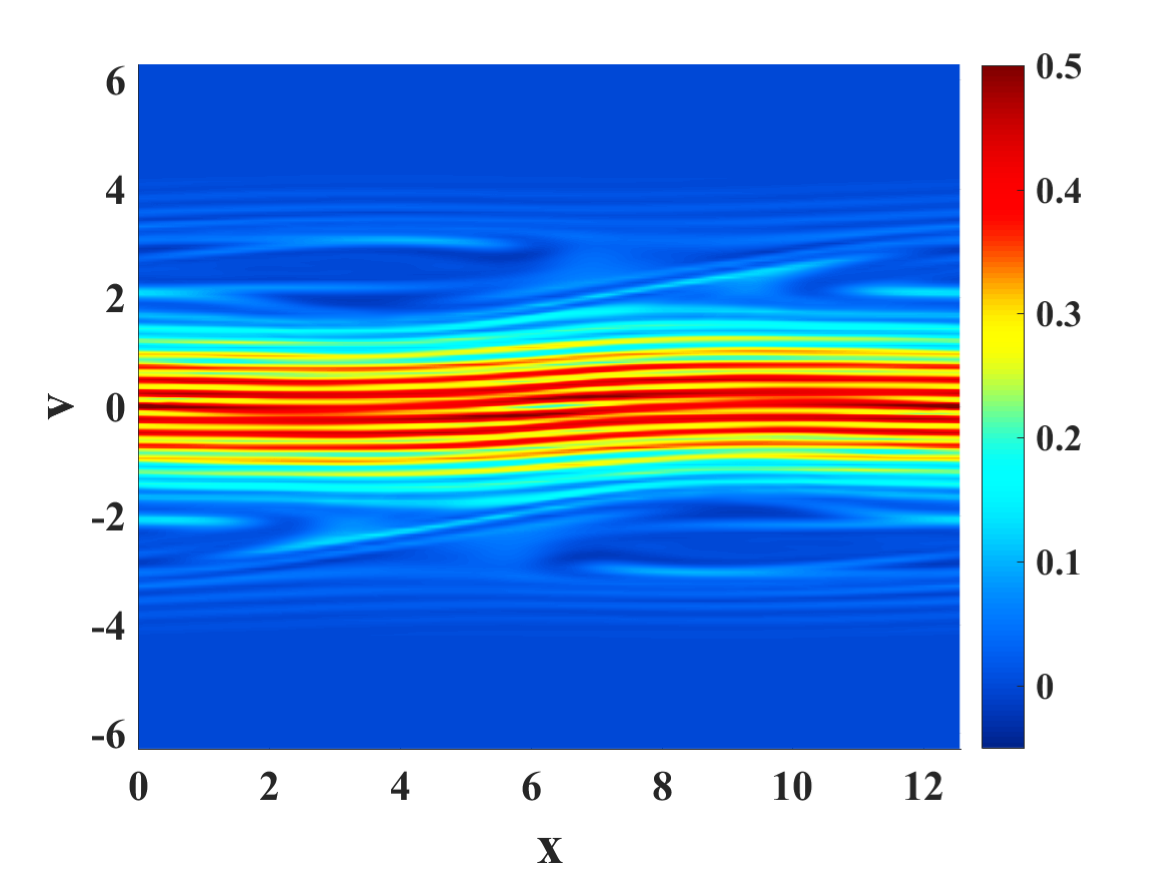}}
			\centerline{\textbf{(a)}}
		\end{minipage}
		\hfill
		\begin{minipage}{0.49\linewidth}
			\centerline{\includegraphics[width=1\linewidth]{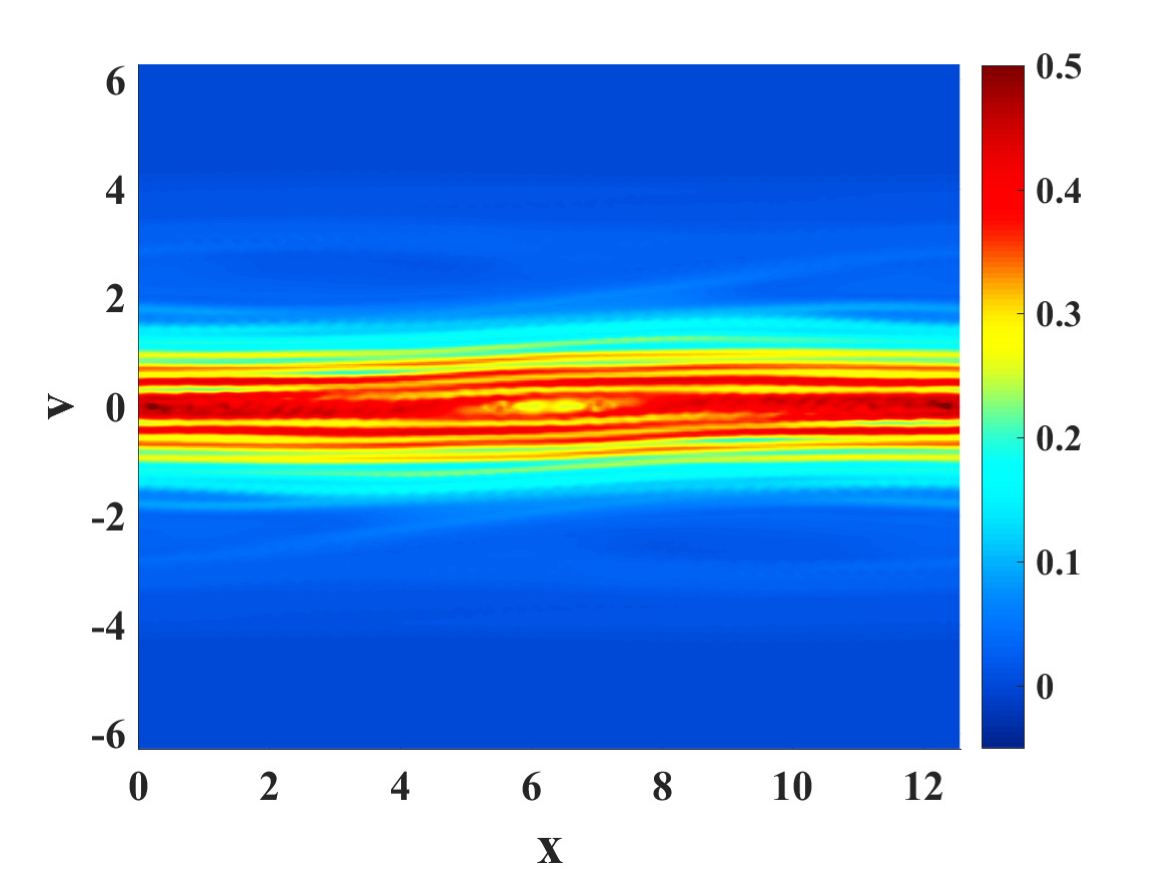}}
			\centerline{\textbf{(b)}}
		\end{minipage}
        \vfill
        \begin{minipage}{0.49\linewidth}
			\centerline{\includegraphics[width=1\linewidth]{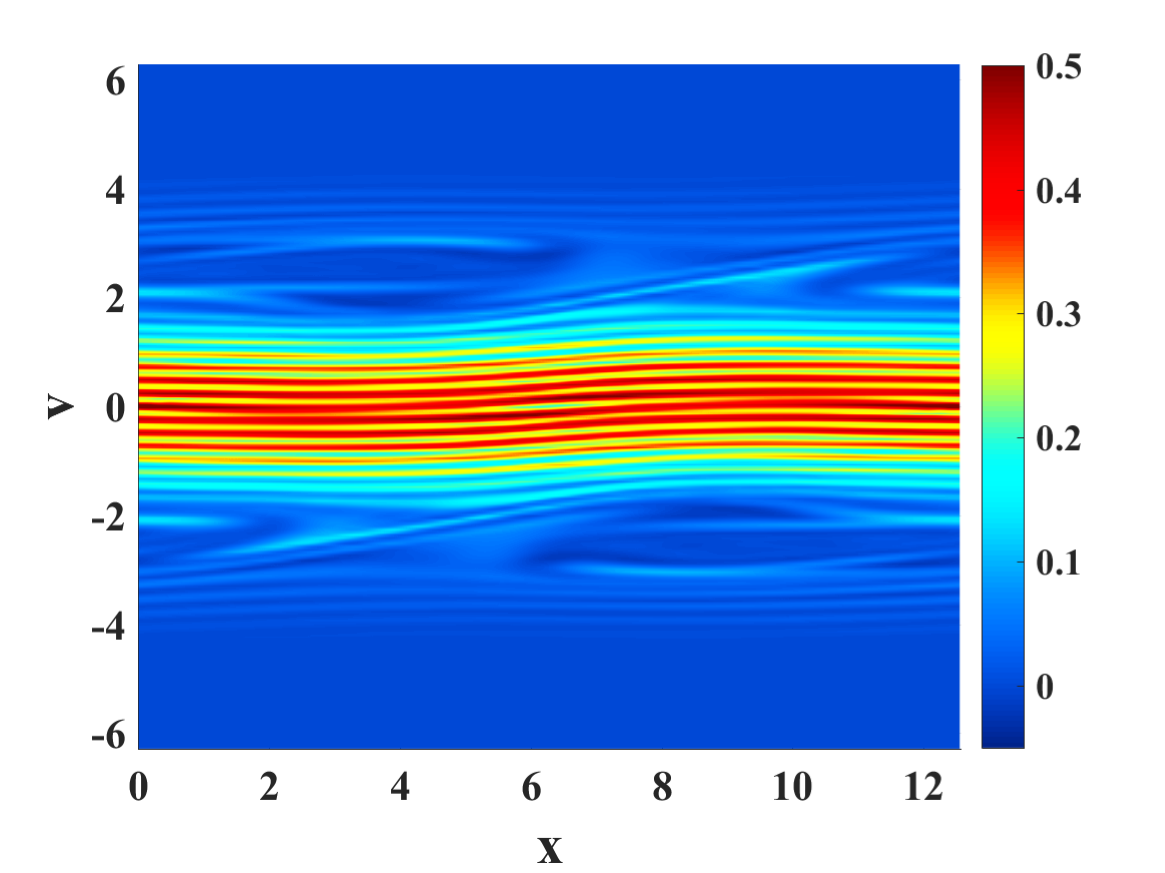}}
			\centerline{\textbf{(c)}}
		\end{minipage}
		\hfill
		\begin{minipage}{0.49\linewidth}
			\centerline{\includegraphics[width=1\linewidth]{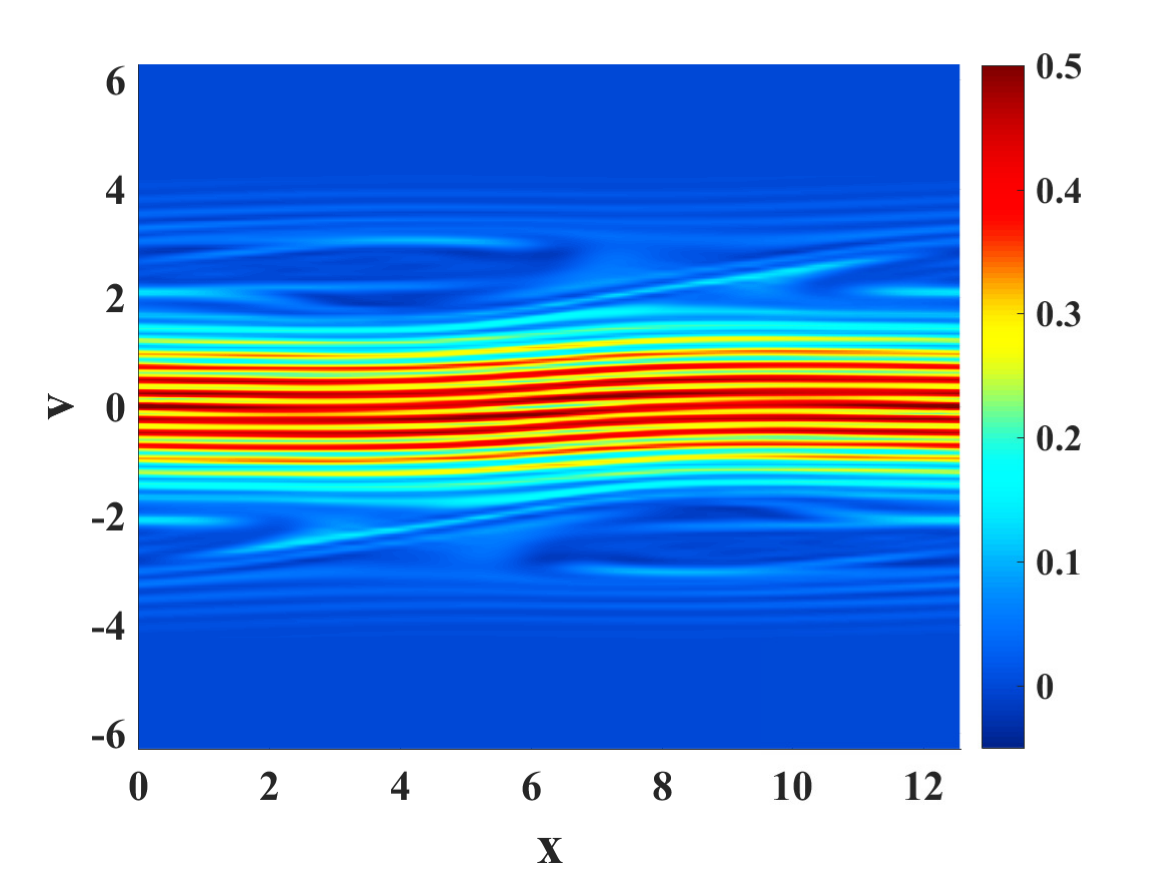}}
			\centerline{\textbf{(d)}}
		\end{minipage}
	\end{center}
	\caption{Strong Landau damping: $N_x=N_v=128$ and $T=50$. Phase space plots of the velocity distribution function simulated from CSLDG method with second-order splitting method $SS$ for the VP system at $\text{CFL}=1$ (a) and ECSLDG method with fourth-order splitting scheme $SS_3$ (b), $SS_{13}$ (c) and $10Lie$ (d) for the VA system at $\text{CFL}=20$.}
	\label{plot-method-s}
\end{figure}


In our numerical tests, we found that although both fourth-order splitting schemes achieve the expected convergence rate, their performances differ significantly. To further illustrate this, Fig.~\ref{plot-method-s} presents phase space plots of velocity distribution function generated by ECSLDG method combined with different fourth-order splitting schemes described in Section \ref{hgspva}. Using VP results as a reference, we observe that at $\text{CFL}=20$, the most efficient fourth-order scheme $SS_3$ (corresponding to $m=1$ in Eq.(\ref{eqss})) produces some unphysical oscillations, as shown in Fig.~\ref{plot-method-s}(b).
As the number of splitting stages increases to $SS_{13}$ (corresponding to $m=6$ in Eq.~(\ref{eqss})), the phase space structure gradually approaches the benchmark result, as shown in Fig.~\ref{plot-method-s}(c). However, the increased number of stages inevitably leads to a significant rise in computational cost.
Encouragingly, the fourth-order $10Lie$ scheme (Eq.~(\ref{10lie})) achieves comparable accuracy to $SS_{13}$ while reducing computational cost by approximately 50\%.

The above tests demonstrate that 
the use of high-order splitting schemes significantly enhances the ability of the ECSLDG method to preserve Gauss's law and produce more high-fidelity results, especially when using large CFL numbers. Based on these observations, we choose the fourth-order splitting scheme $10Lie$ for the following numerical examples.

\begin{figure}[htbp]
	\begin{center}
		\begin{minipage}{0.49\linewidth}
			\centerline{\includegraphics[width=1\linewidth]{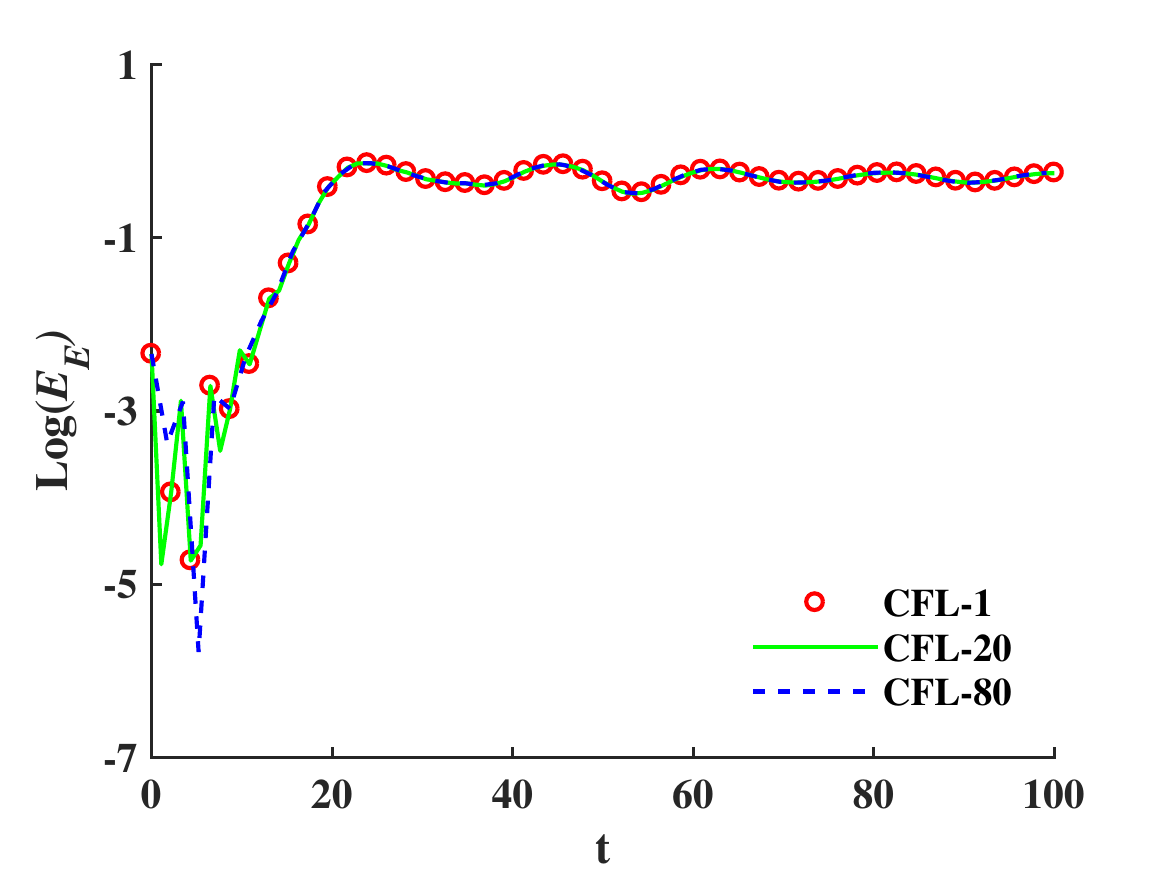}}
			\centerline{\textbf{(a)}}
		\end{minipage}
		\hfill
		\begin{minipage}{0.49\linewidth}
			\centerline{\includegraphics[width=1\linewidth]{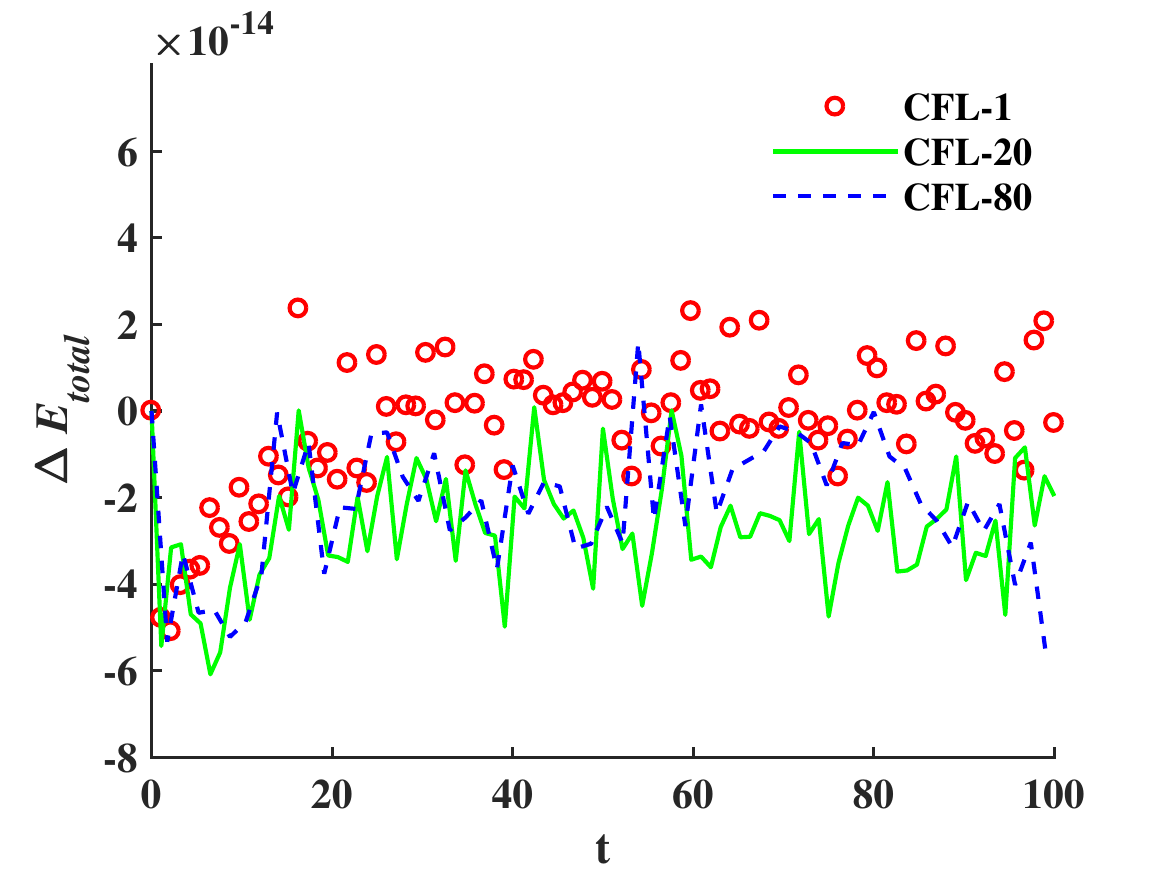}}
			\centerline{\textbf{(b)}}
		\end{minipage}
		\vfill
		\begin{minipage}{0.49\linewidth}
			\centerline{\includegraphics[width=1\linewidth]{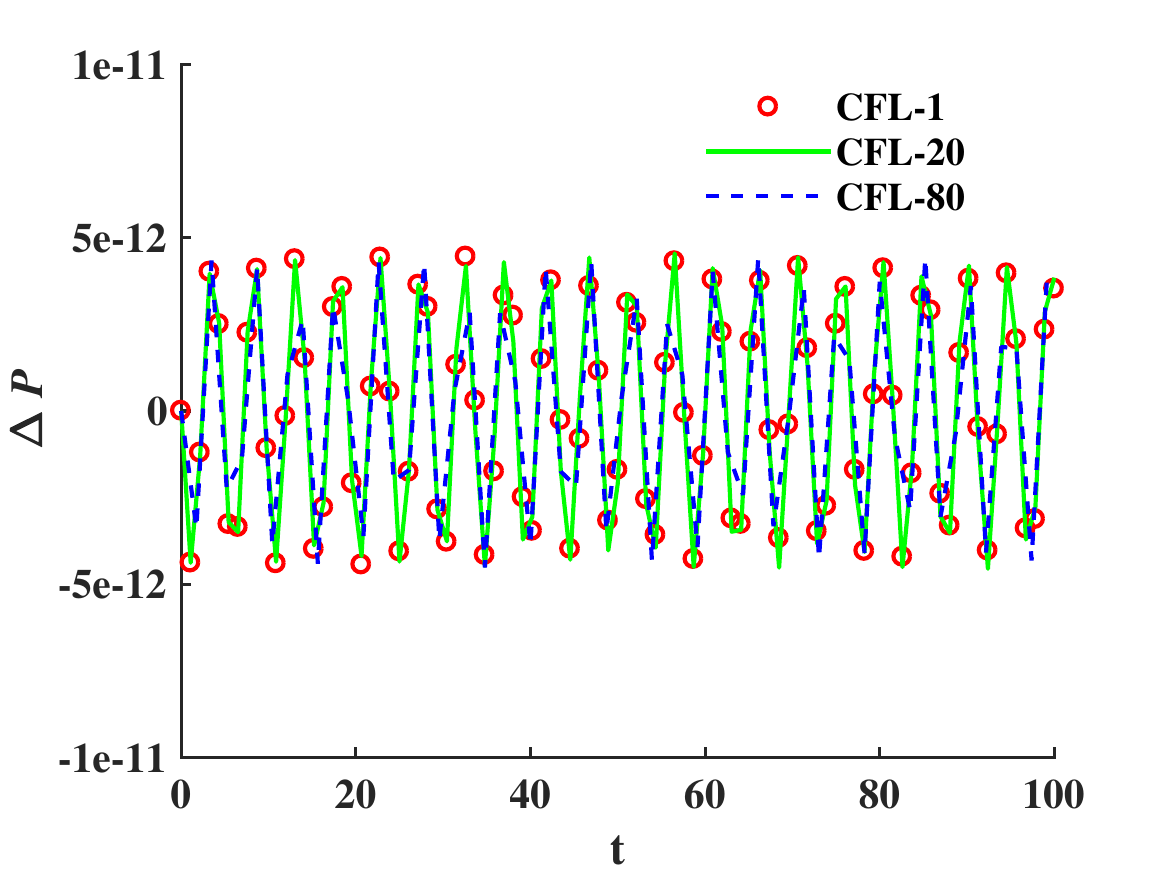}}
			\centerline{\textbf{(c)}}
		\end{minipage}
		\hfill
		\begin{minipage}{0.49\linewidth}
			\centerline{\includegraphics[width=1\linewidth]{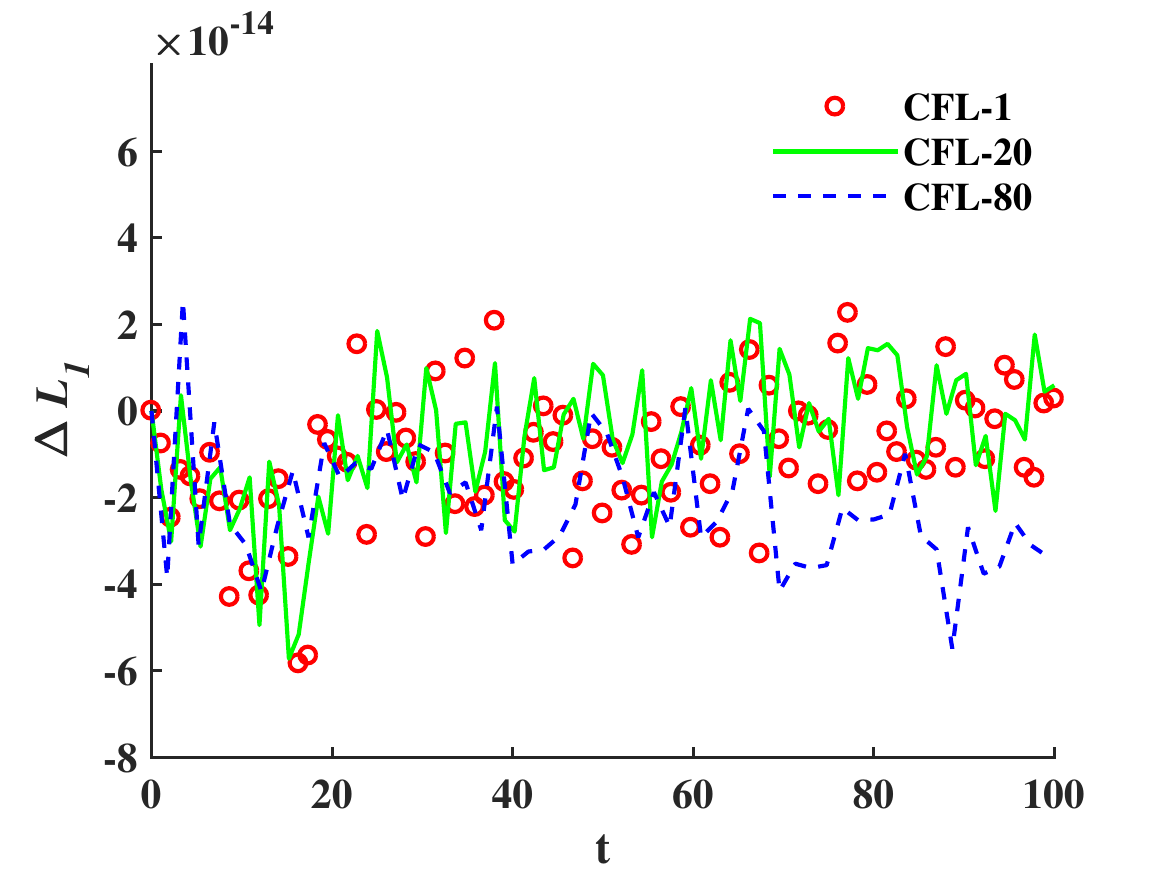}}
			\centerline{\textbf{(d)}}
		\end{minipage}
	\end{center}
	\caption{Two-stream instability I: $N_x=N_v=128$. Time evolution of the electric energy (a), relative deviation of the total energy (b), deviation of momentum (c) and relative deviation of mass (d) with different CFL.}
	\label{Physics-cfl-st2}
\end{figure}

\subsection{Two-stream instability I \label{test3}}
In this section, we apply the proposed ECSLDG method to simulate the two-stream instability for two main purposes. First, we test the conservation properties of the fourth-order ECSLDG method using different CFL numbers. Second, we investigate the advantages of the ECSLDG method with high-order splitting schemes compared to the second-order scheme. 

The initial velocity distribution function is given by \cite{CAI2018529}
\begin{align*}
f_0(x,v)=\frac{2}{7\sqrt{2\pi}}(1+5v^2)(1+\alpha((\cos(2\kappa x)+\cos(3\kappa x))/1.2+\cos(\kappa x)))\exp(-\frac{v^2}{2}),
\end{align*}
where $\alpha=0.01$, $\kappa=0.5$. The computational domain is set to $[0,{2\pi}/{\kappa}]\times [-v_m,v_m]$.
Unless otherwise stated, simulations employ $v_m=10$ and a grid of $128 \times 128$.

First, we investigate conservative properties of the ECSLDG method with different CFL numbers. In Fig.~\ref{Physics-cfl-st2}, we present the time evolution of the electric energy, relative deviation of the total energy, deviation of momentum $\Delta P(t)=P(t)-P(0)$ and relative deviation of mass $\Delta L_1(t)=[L_1(t)-L_1(0)]/L_1(0)$ for different CFL numbers. As shown in Fig.~\ref{Physics-cfl-st2}(a), the ECSLDG method with $\text{CFL}=80$ 
produces nearly the same electric energy as the one with $\text{CFL}=1$. Moreover, even with $\text{CFL}=80$, the ECSLDG method still maintains the relative energy deviation at $O(10^{-14})$ and the relative mass deviation at $O(10^{-14})$. Additionally, the momentum is also well preserved. Clearly, the conservation properties of the proposed ECSLDG method are independent of the temporal resolution.
\begin{figure}[htbp]
	\begin{center}
		\begin{minipage}{0.475\linewidth}
			\centerline{\includegraphics[width=1\linewidth]{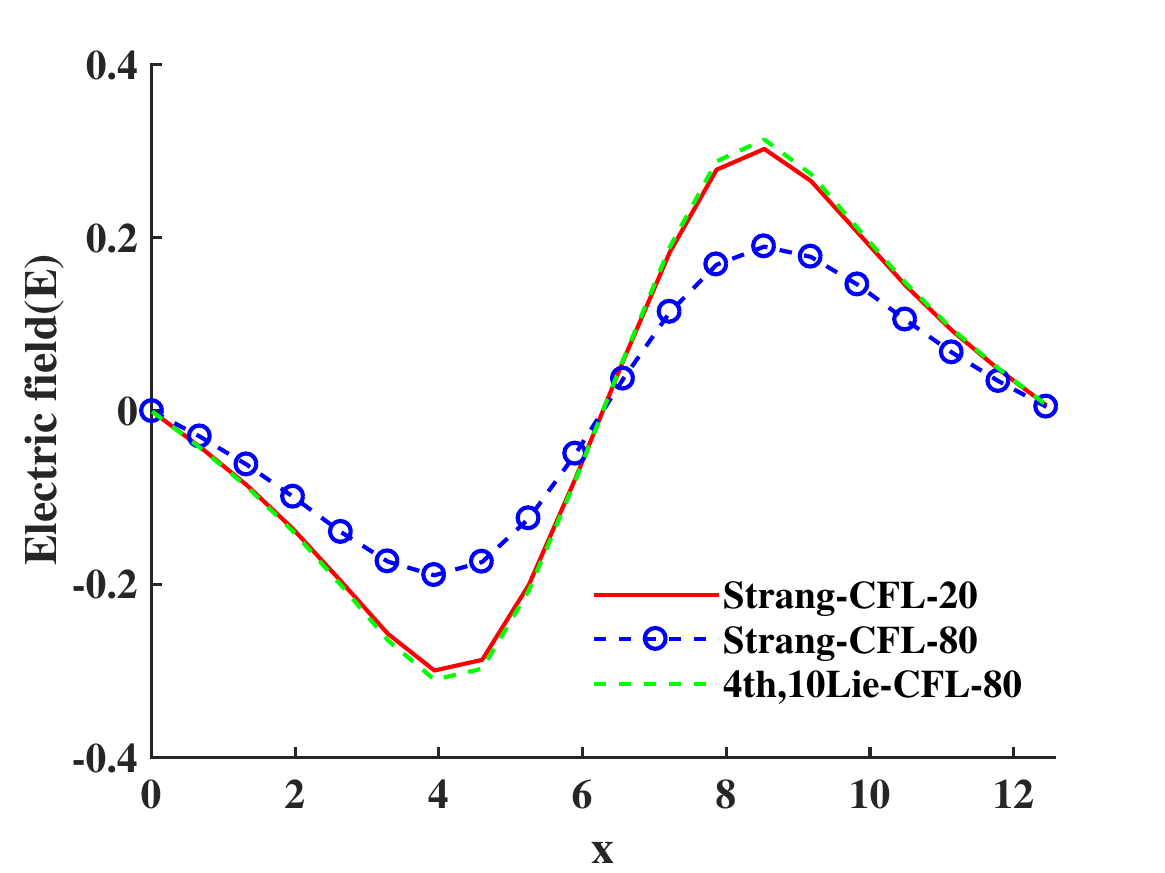}}
			\centerline{\textbf{(a)}}
		\end{minipage}
		\hfill
		\begin{minipage}{0.475\linewidth}
			\centerline{\includegraphics[width=1\linewidth]{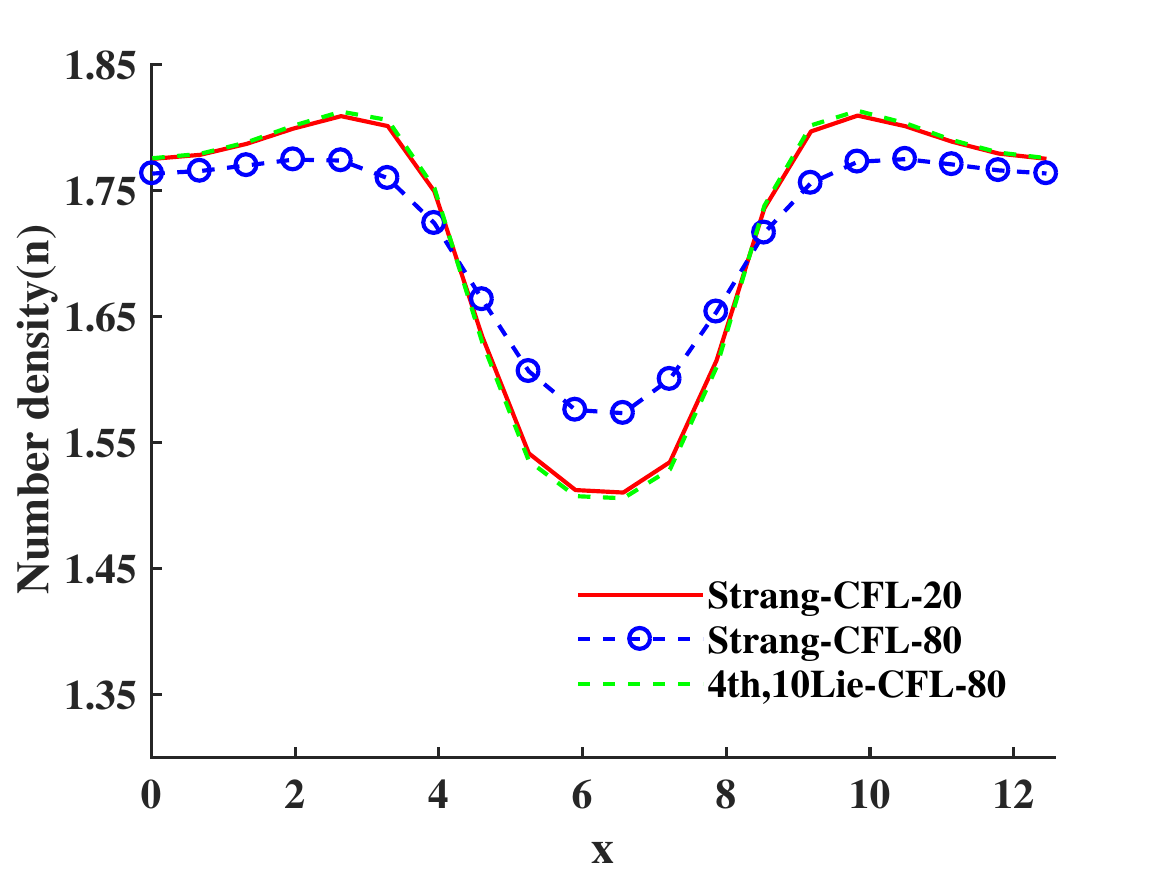}}
			\centerline{\textbf{(b)}}
		\end{minipage}
	\end{center}
	\caption{Two-stream instability I: $N_x=N_v=128$ and $T=20$. Plot of the electric field (a) and the number density (b) using the ECSLDG method with different splitting methods.}
	\label{electric-temporal-st2}
\end{figure}

Then, we examine the performance of the ECSLDG method with high-order splitting schemes. In Fig.~\ref{electric-temporal-st2}, we show the electric field and number density produced by our ECSLDG method combined with different splitting schemes at $T = 20$. With a large time step (CFL = 80), the ECSLDG method with the fourth-order splitting scheme yields results nearly identical to the second-order scheme with CFL = 20, while the second-order scheme with CFL = 80 shows noticeable deviations.

\begin{figure}[htbp]
	\begin{center}
		\begin{minipage}{0.3\linewidth}
			\centerline{\includegraphics[width=1.2\linewidth]{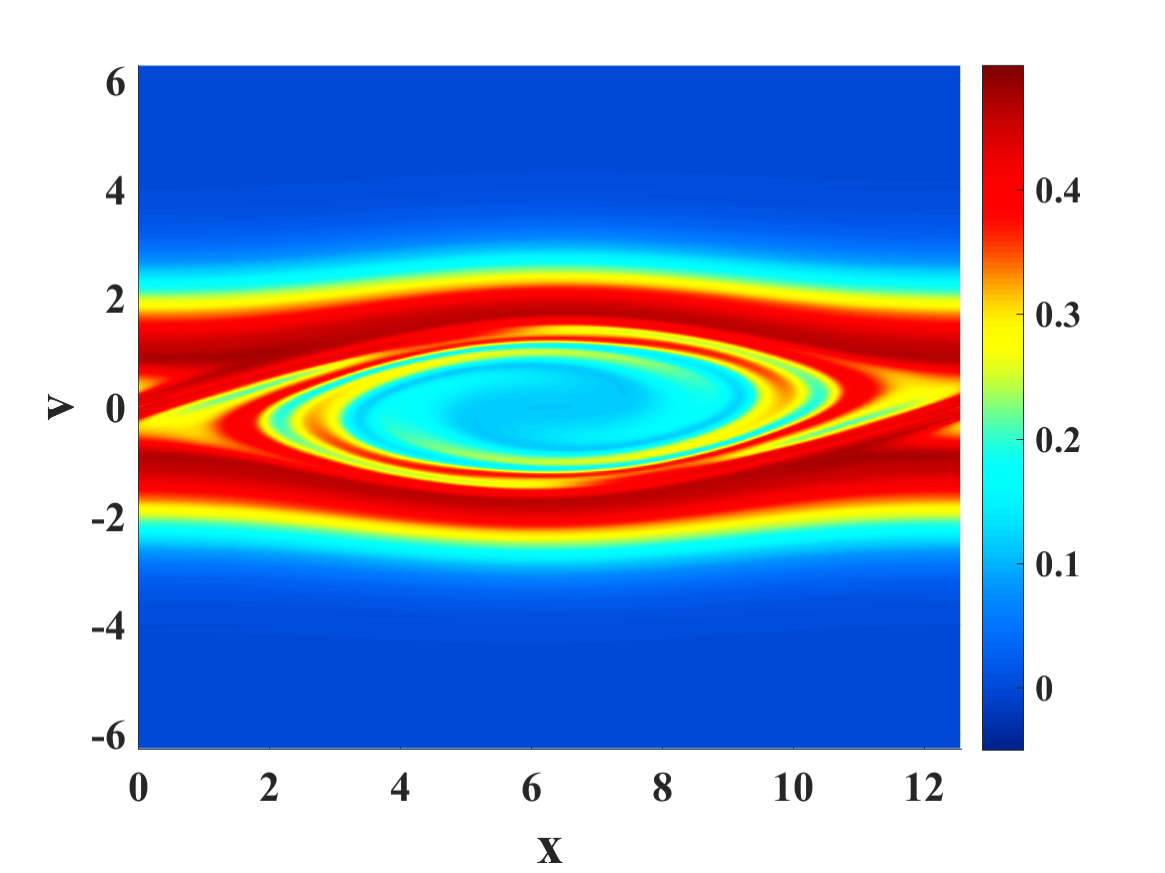}}
			\centerline{\textbf{(a)}}
		\end{minipage}
		\hfill
		\begin{minipage}{0.3\linewidth}
			\centerline{\includegraphics[width=1.2\linewidth]{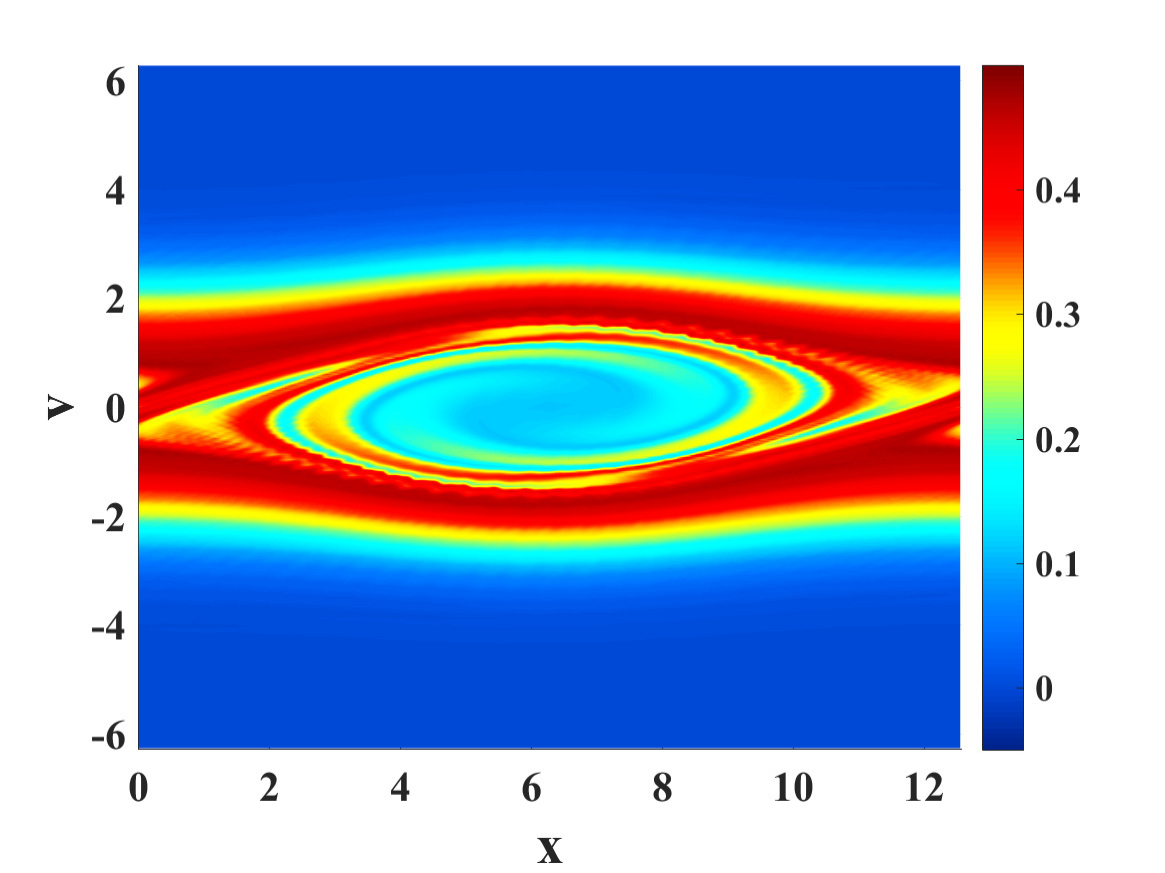}}
			\centerline{\textbf{(b)}}
		\end{minipage}
		\hfill
		\begin{minipage}{0.3\linewidth}
			\centerline{\includegraphics[width=1.2\linewidth]{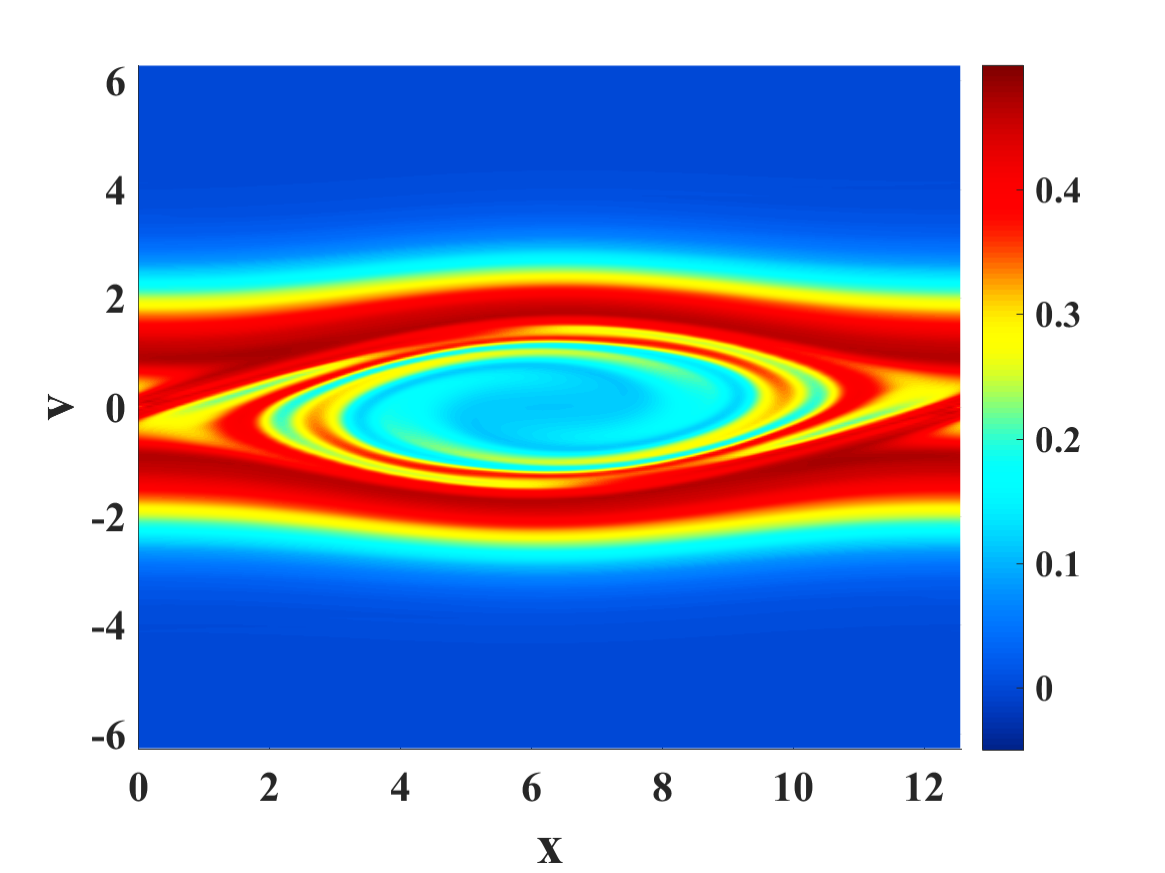}}
			\centerline{\textbf{(c)}}
		\end{minipage}
	\end{center}
        \caption{Two-stream instability I: $N_x=N_v=128$, CFL = 20 and $T=50$. Phase space plots of the velocity distribution function simulated from CSLDG method with second-order Strang splitting method for the VP system (a) and ECSLDG method with the second-order Strang splitting method (b) and the fourth-order splitting method $10Lie$ (c) for the VA system.}
	\label{plot-st2}
\end{figure}

The distinction between the second-order and fourth-order splitting schemes can also be observed in the behavior of the microscopic distribution function.
Fig.~\ref{plot-st2} presents the phase space plots of the velocity distribution function produced from the ECSLDG method combined with different splitting schemes at $T=50$. The phase plot simulated from the CSLDG method with the second-order Strang splitting scheme  for the VP system is used as a reference benchmark. When CFL = 20, some unphysical oscillations appear in the phase space plot of the velocity distribution function generated by the ECSLDG method with the second-order Strang splitting scheme, while the plot generated by the ECSLDG method with the fourth-order splitting scheme $10Lie$ is nearly identical to the reference plot.

These results confirm that the proposed ECSLDG method conserves both mass and total energy for the VA system, regardless of temporal resolution. Furthermore, by employing high-order splitting schemes, the ECSLDG method maintains accuracy even at large CFL numbers, making it well suited for long-time simulations.

\subsection{Two-stream instability II \label{test4}}
In this section, the ECSLDG method is employed to simulate a symmetric two-stream instability for three primary objectives.
First, we examine the conservation properties of the proposed ECSLDG method under different spatial resolutions.
Second, we examine the physical performance under different orders of spatial and temporal discretizations.
Finally, we assess the robustness of the method when the time step size is insufficient to resolve the plasma oscillation period.

The initial condition for the symmetric two-stream instability is given by \cite{CROUSEILLES20101927,UT2008}
\begin{align*}
f_0(x,v)=\displaystyle\frac{1}{2\sqrt{2\pi}v_t}\biggl(1+\alpha\cos(\kappa x))(\exp(-(v-v_0)^2/2v_t^2)+\exp(-(v+v_0)^2/2v_t^2)\biggr),
\end{align*}
where $\alpha=0.05$, $\kappa=2/13$, $v_0=0.99$ and $v_t=0.3$. The computational domain is set to $\displaystyle [0,13\pi]\times [-v_m,v_m]$. Here, we set $v_m=10$ and employ a  mesh grid of $128 \times 128$. 
\begin{figure}[htbp]
	\begin{center}
		\begin{minipage}{0.49\linewidth}
			\centerline{\includegraphics[width=1\linewidth]{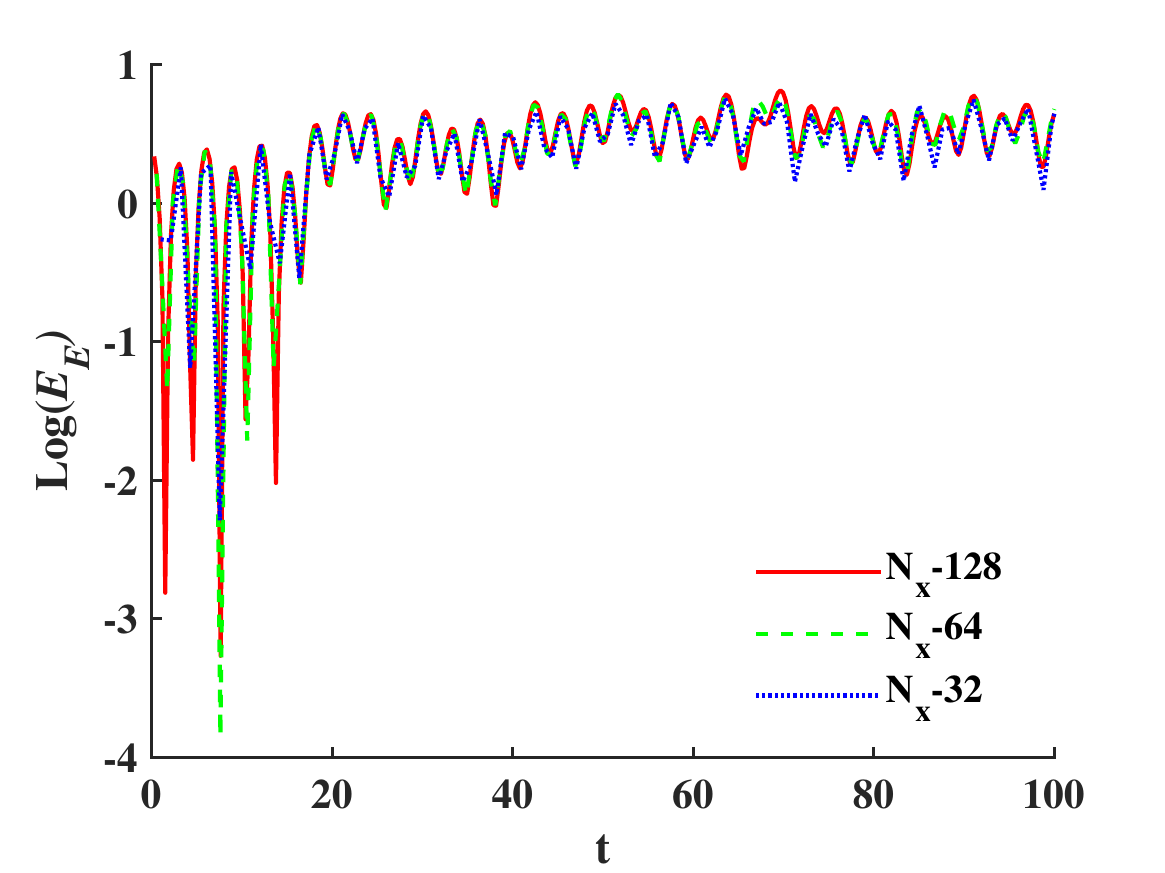}}
			\centerline{\textbf{(a)}}
		\end{minipage}
		\hfill
		\begin{minipage}{0.49\linewidth}
			\centerline{\includegraphics[width=1\linewidth]{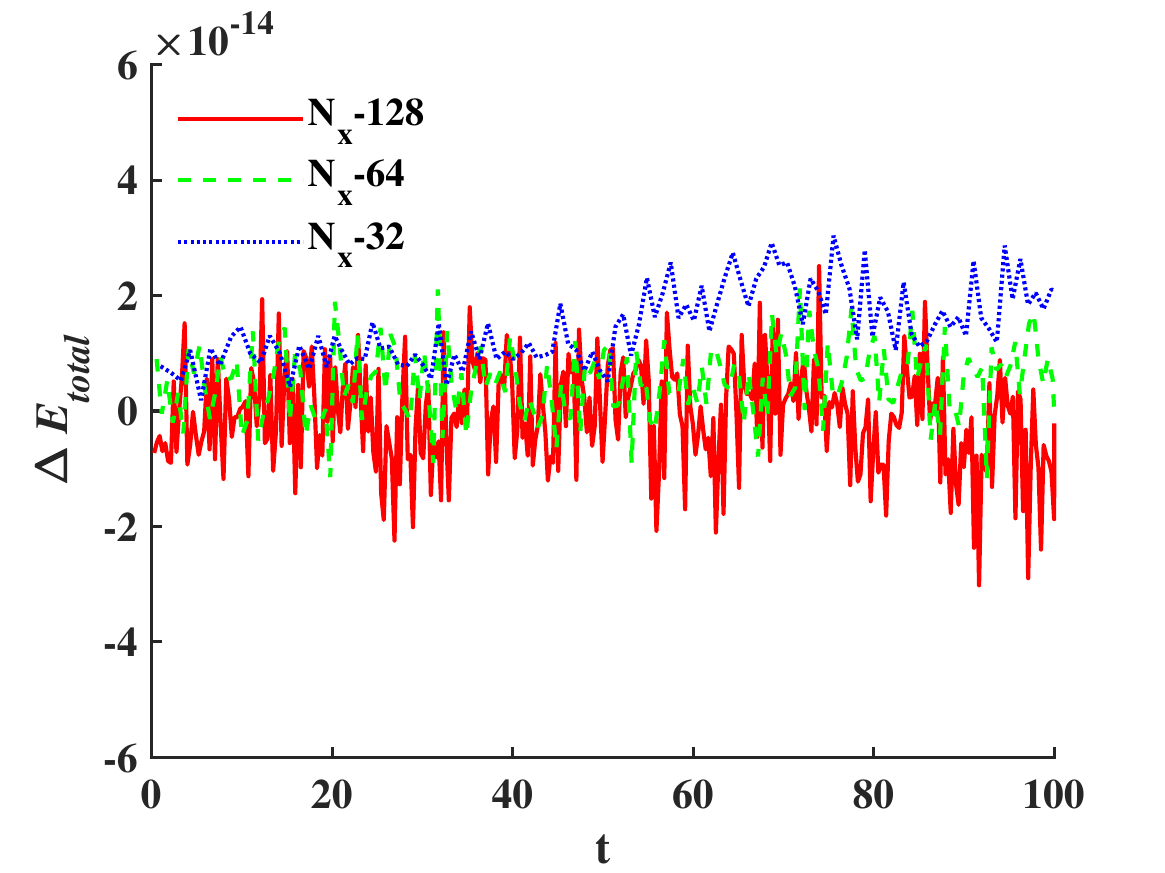}}
			\centerline{\textbf{(b)}}
		\end{minipage}
		\vfill
		\begin{minipage}{0.49\linewidth}
			\centerline{\includegraphics[width=1\linewidth]{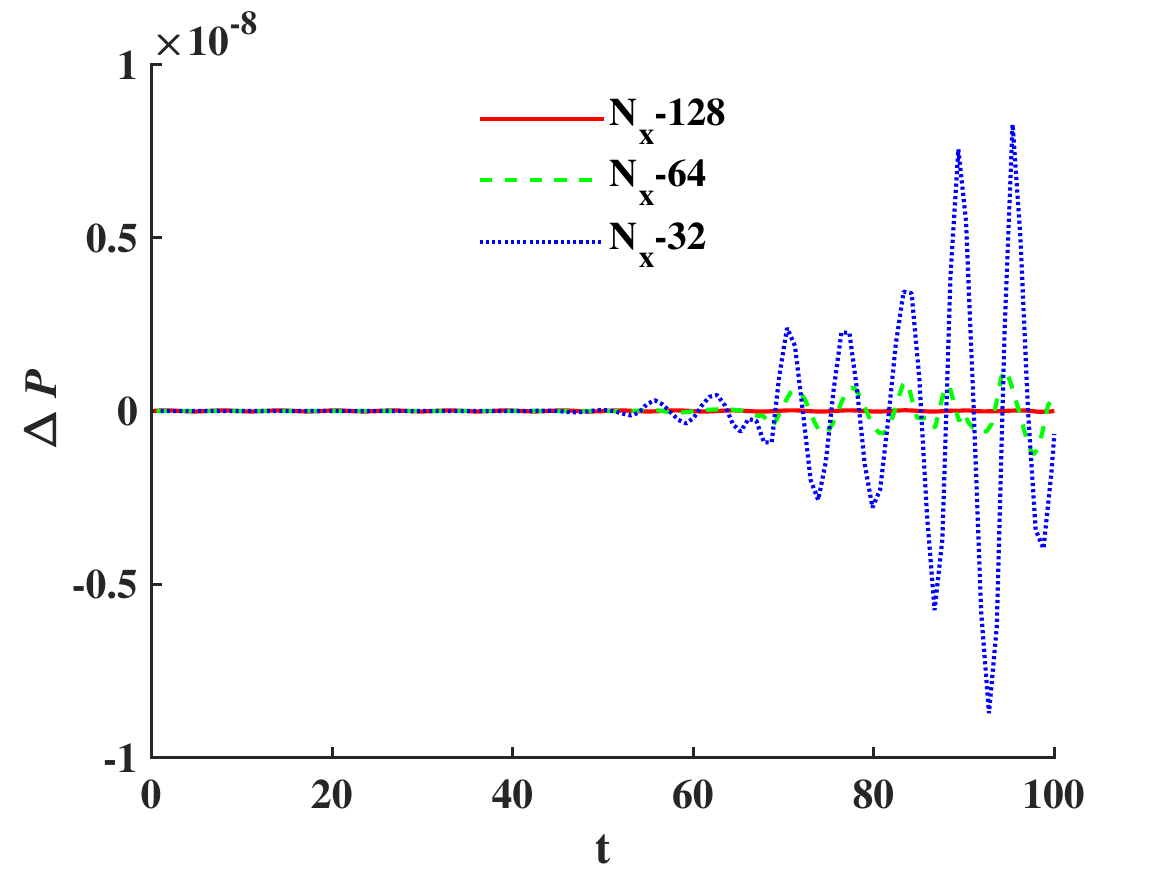}}
			\centerline{\textbf{(c)}}
		\end{minipage}
		\hfill
		\begin{minipage}{0.49\linewidth}
			\centerline{\includegraphics[width=1\linewidth]{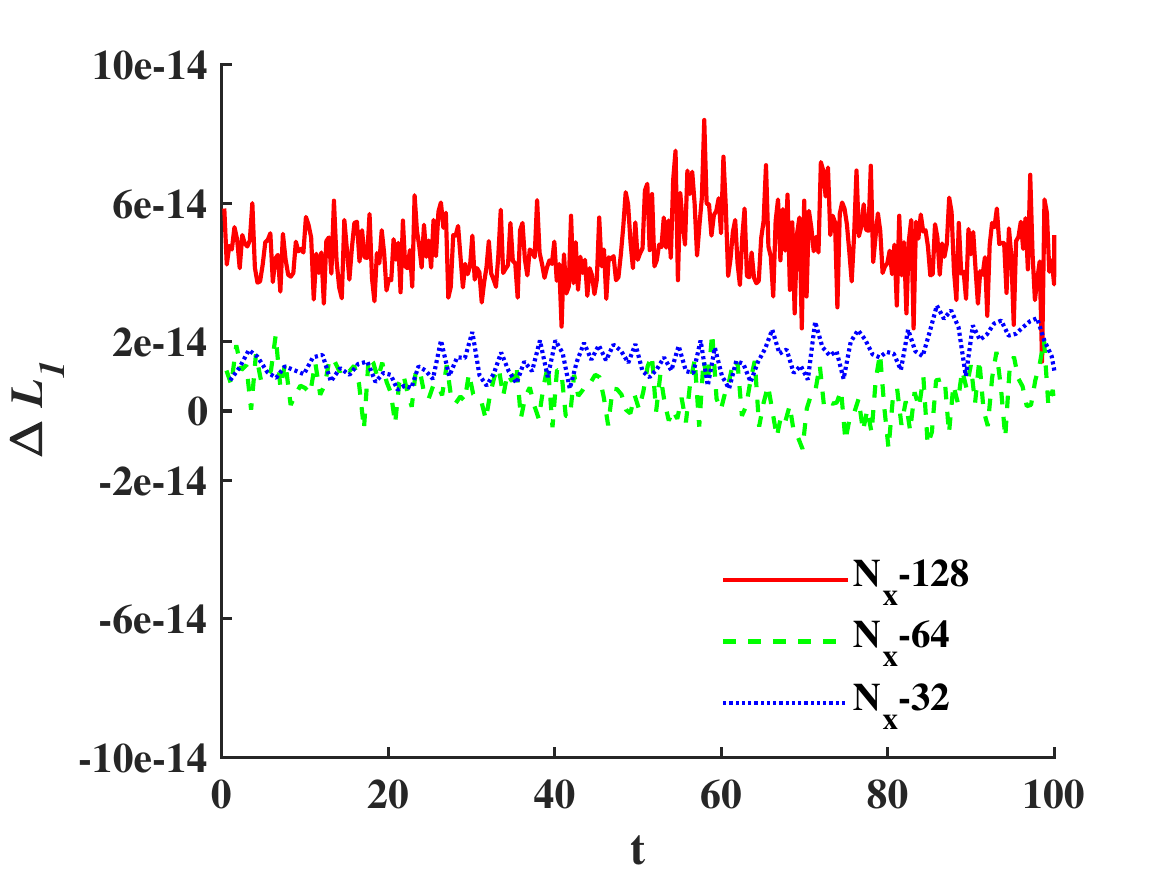}}
			\centerline{\textbf{(d)}}
		\end{minipage}
	\end{center}
	\caption{Two-stream instability II: $\lambda =1$ and CFL = 10. Time evolution of the electric energy (a), relative deviation of total energy (b), deviation of momentum
		(c) and relative deviation of mass (d) with different $N_x$.}
	\label{electric-Nx-st3}
\end{figure}

First, we investigate the conservation properties of the proposed ECSLDG method under different physical grid resolutions $N_x$.
Fig.~\ref{electric-Nx-st3} presents the evolution of the electric energy $E_{\bm{E}}$, the relative deviation of the total energy $\Delta E_{total}$, the relative deviation of mass $\Delta L_1$, and the deviation of momentum $\Delta P$ at $\mathrm{CFL}=10$.
As shown in Fig.~\ref{electric-Nx-st3}(a), the ECSLDG method with a coarse grid of $N_x=32$ can accurately capture the evolution of the electric field, closely matching results obtained with a finer grid of $N_x=128$.
Moreover, the results demonstrate that the ECSLDG method conserves the mass and total energy of the VA system, independent of the spatial resolution. Even with the coarse grid of $N_x=32$, the momentum deviation in this simulation is quite small.

\begin{figure}[htbp]
	\begin{center}
		\begin{minipage}{0.3\linewidth}
			\centerline{\includegraphics[width=1.2\linewidth]{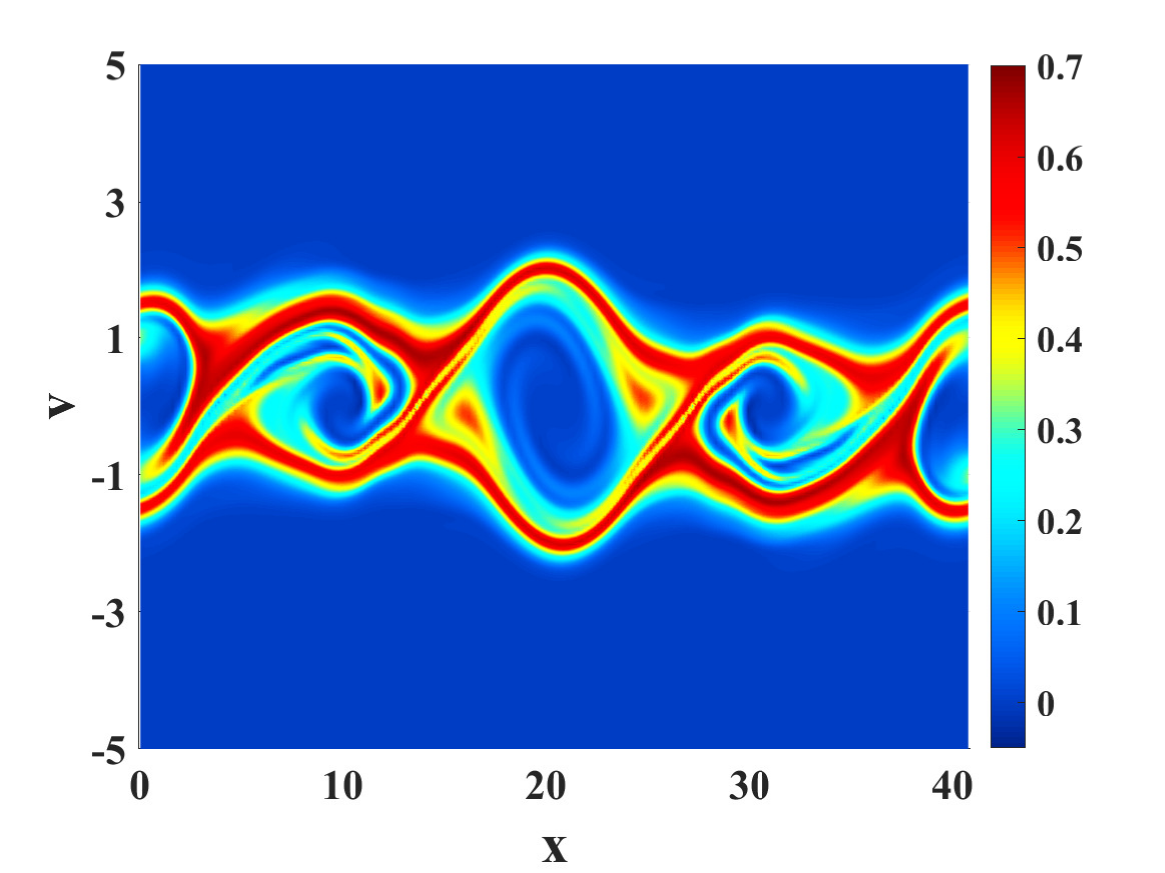}}
			\centerline{\textbf{$P^1,128\times 256,\mathrm{VP}$}}
		\end{minipage}
		\hfill
		\begin{minipage}{0.3\linewidth}
			\centerline{\includegraphics[width=1.2\linewidth]{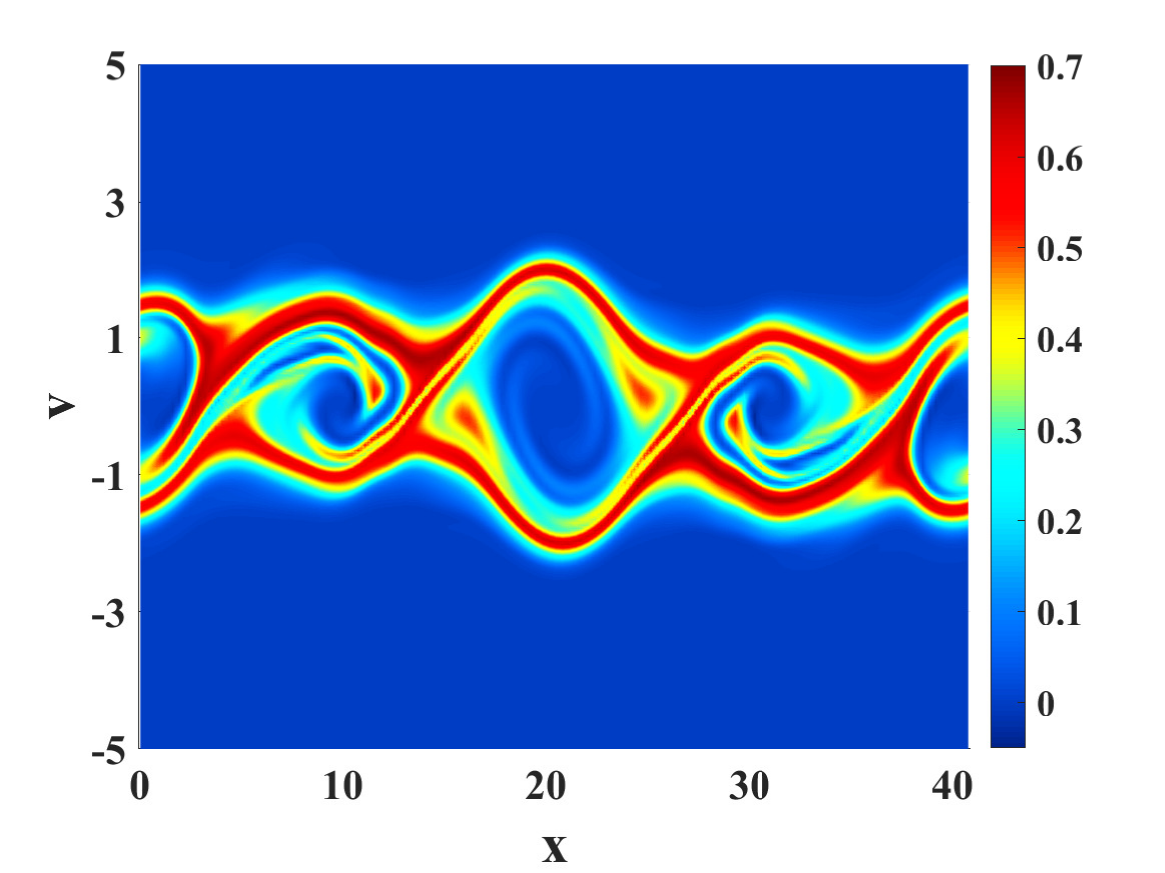}}
			\centerline{\textbf{$P^1,128\times 256,\mathrm{VA}$}}
		\end{minipage}
            \hfill
		\begin{minipage}{0.3\linewidth}
			\centerline{\includegraphics[width=1.2\linewidth]{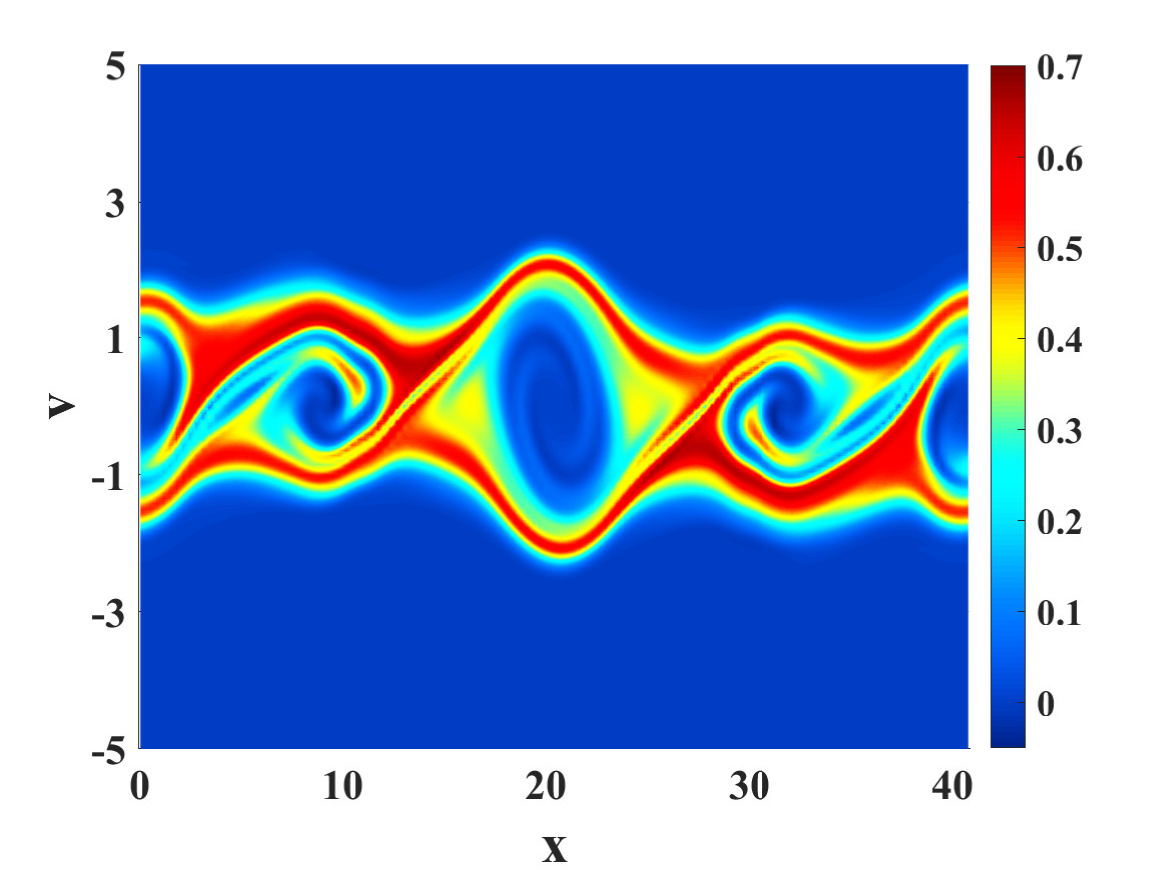}}
			\centerline{\textbf{$P^1,128\times 256,\mathrm{VA}$}}
		\end{minipage}
            \vfill
            \begin{minipage}{0.3\linewidth}			\centerline{\includegraphics[width=1.2\linewidth]{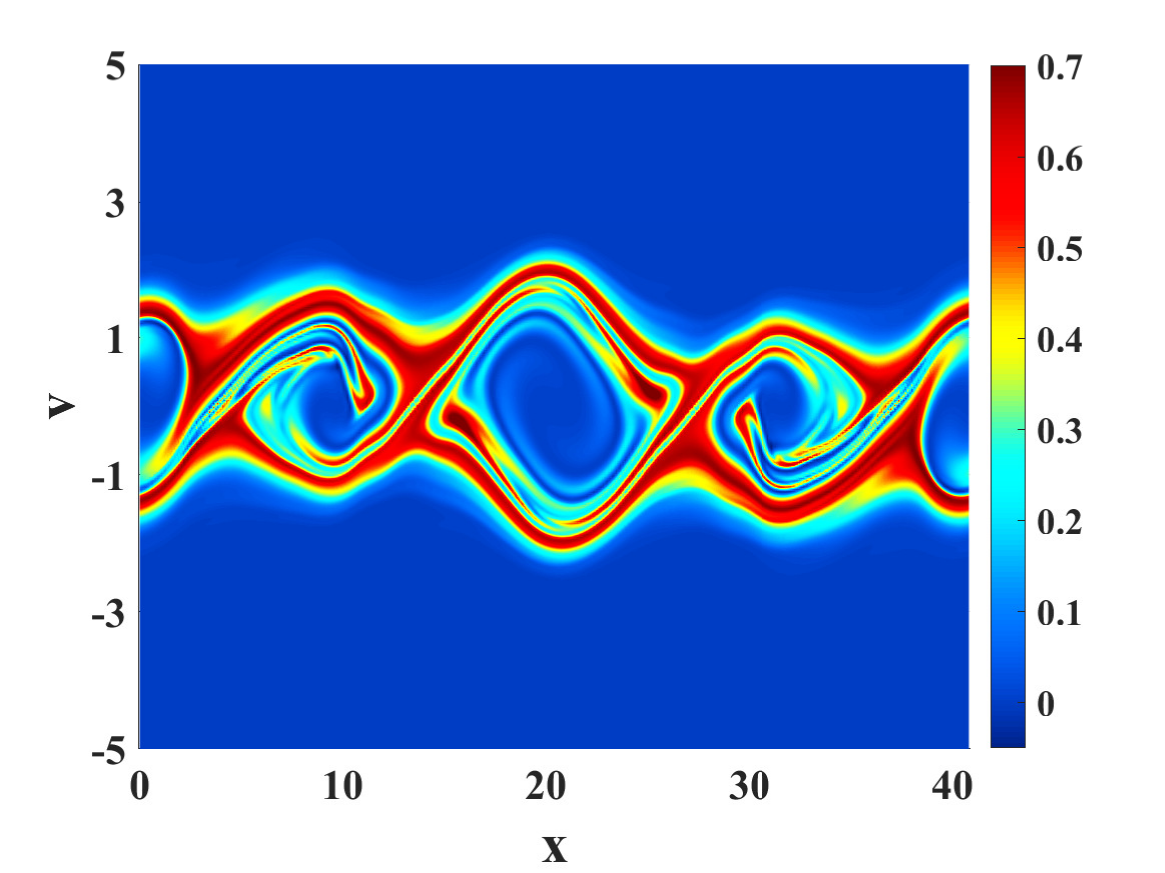}}
			\centerline{\textbf{$P^2,128\times 256,\mathrm{VP}$}}
		\end{minipage}
		\hfill
		\begin{minipage}{0.3\linewidth}
			\centerline{\includegraphics[width=1.2\linewidth]{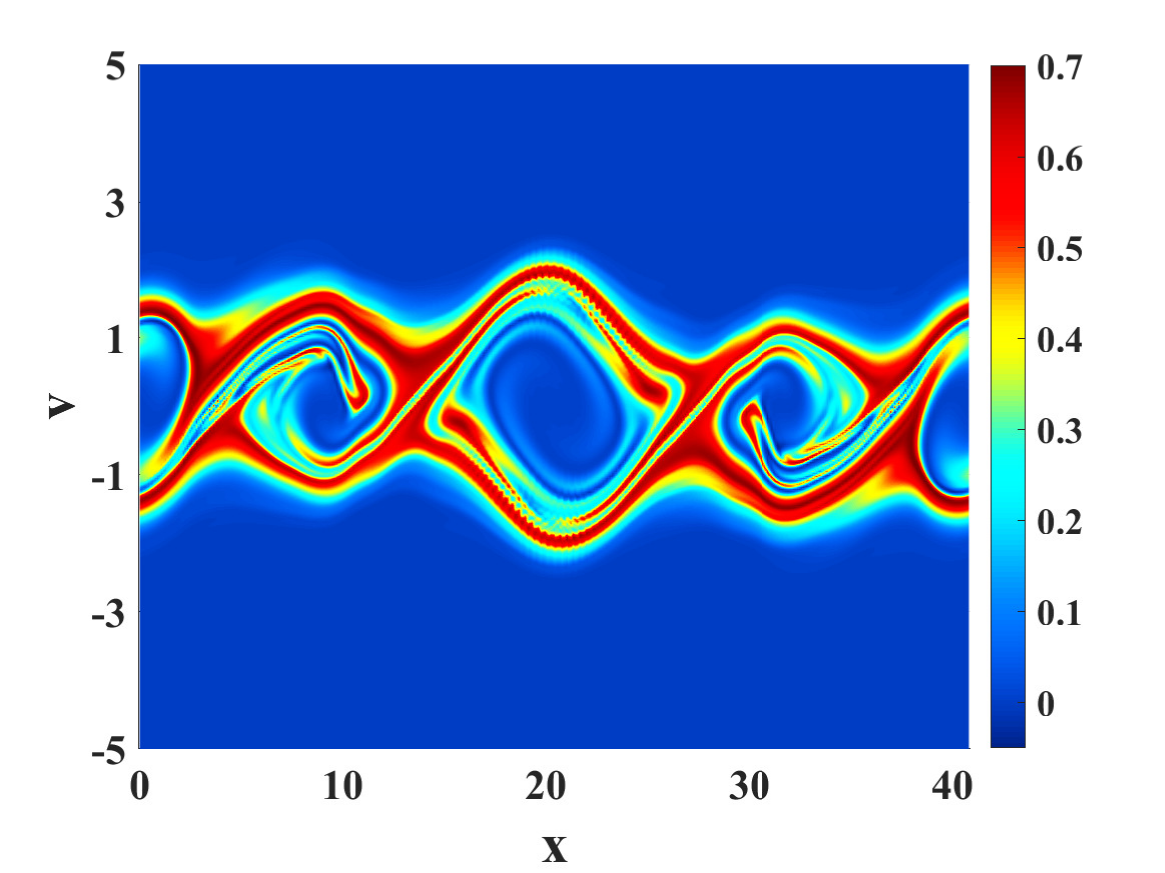}}
			\centerline{\textbf{$P^2,128\times 256,\mathrm{VA}$}}
		\end{minipage}
            \hfill
		\begin{minipage}{0.3\linewidth}
			\centerline{\includegraphics[width=1.2\linewidth]{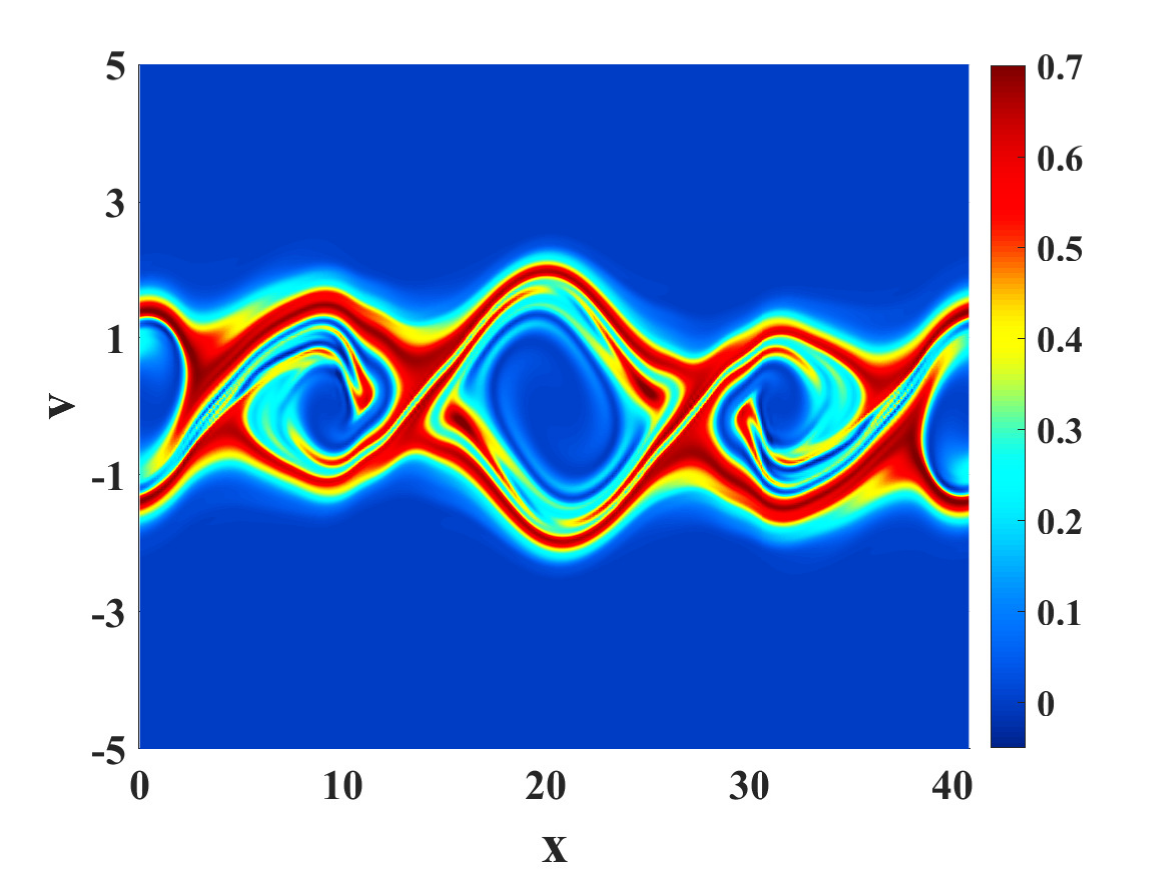}}
			\centerline{\textbf{$P^2,128\times 256,\mathrm{VA}$}}
		\end{minipage}
            \vfill
            \begin{minipage}{0.3\linewidth}			\centerline{\includegraphics[width=1.2\linewidth]{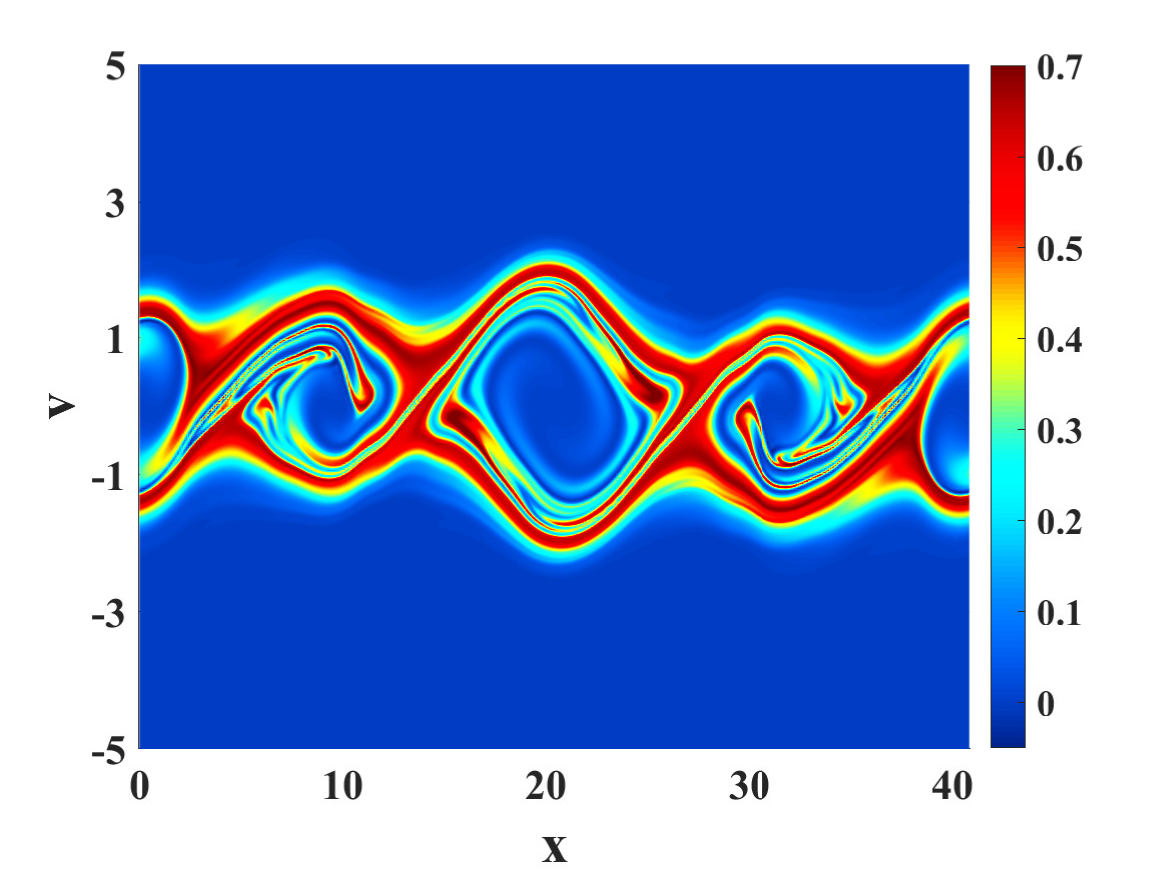}}
			\centerline{\textbf{$P^3,128\times 256,\mathrm{VP}$}}
		\end{minipage}
		\hfill
		\begin{minipage}{0.3\linewidth}
			\centerline{\includegraphics[width=1.2\linewidth]{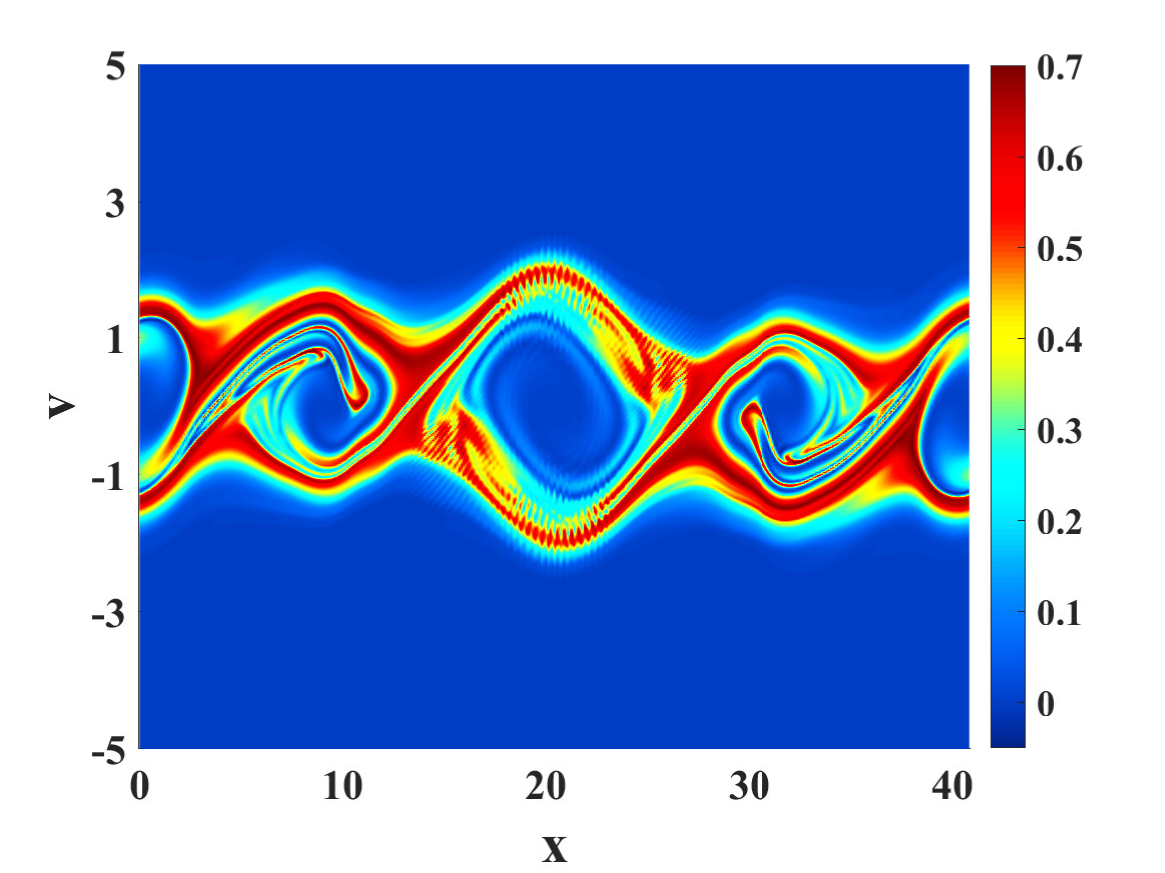}}
			\centerline{\textbf{$P^3,128\times 256,\mathrm{VA}$}}
		\end{minipage}
            \hfill
		\begin{minipage}{0.3\linewidth}
			\centerline{\includegraphics[width=1.2\linewidth]{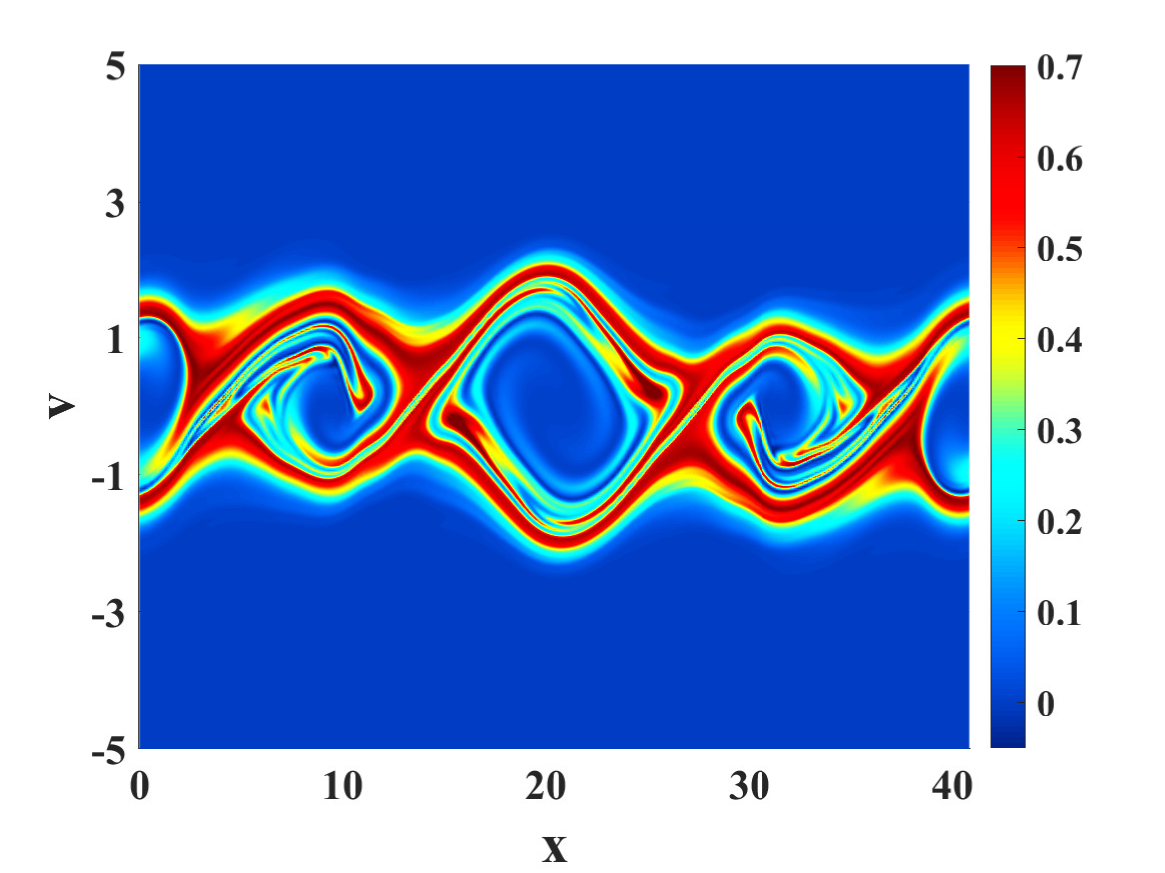}}
			\centerline{\textbf{$P^3,128\times 256,\mathrm{VA}$}}
		\end{minipage}
        
	\end{center}
	\caption{Two-stream instability II: $\lambda = 1$, $\Delta t=0.2$, $N_x=128$, $N_v=256$, and $T=50$. Phase space plots of the velocity distribution function for the VP system using the second-order Strang splitting method  (left), the ECSLDG method using the second-order Strang splitting method (middle) and the ECSLDG method using the fourth-order method $10Lie$ (right).  }
	\label{plot-st3}
\end{figure}
Then, we examine the impact of consistent high-order accuracy in both space and time on capturing physical phenomena. Figure~\ref{plot-st3} shows phase space plots of the velocity distribution function obtained using the ECSLDG method with varying spatial polynomial orders and temporal splitting schemes. The CSLDG result for the VP system serves as a reference.  The results demonstrate that, for a fixed splitting scheme,  increasing spatial accuracy of ECSLDG significantly enhances the resolution of fine structures. However, when the second-order Strang splitting scheme is used, some unphysical oscillations persist in the phase space plots, even with higher-order spatial discretizations. In contrast, the fourth-order scheme ($10Lie$) effectively eliminates these oscillations, yielding results that closely match the reference solution.

\begin{figure}[htbp]
	\begin{center}
		\begin{minipage}{0.49\linewidth}
			\centerline{\includegraphics[width=1\linewidth]{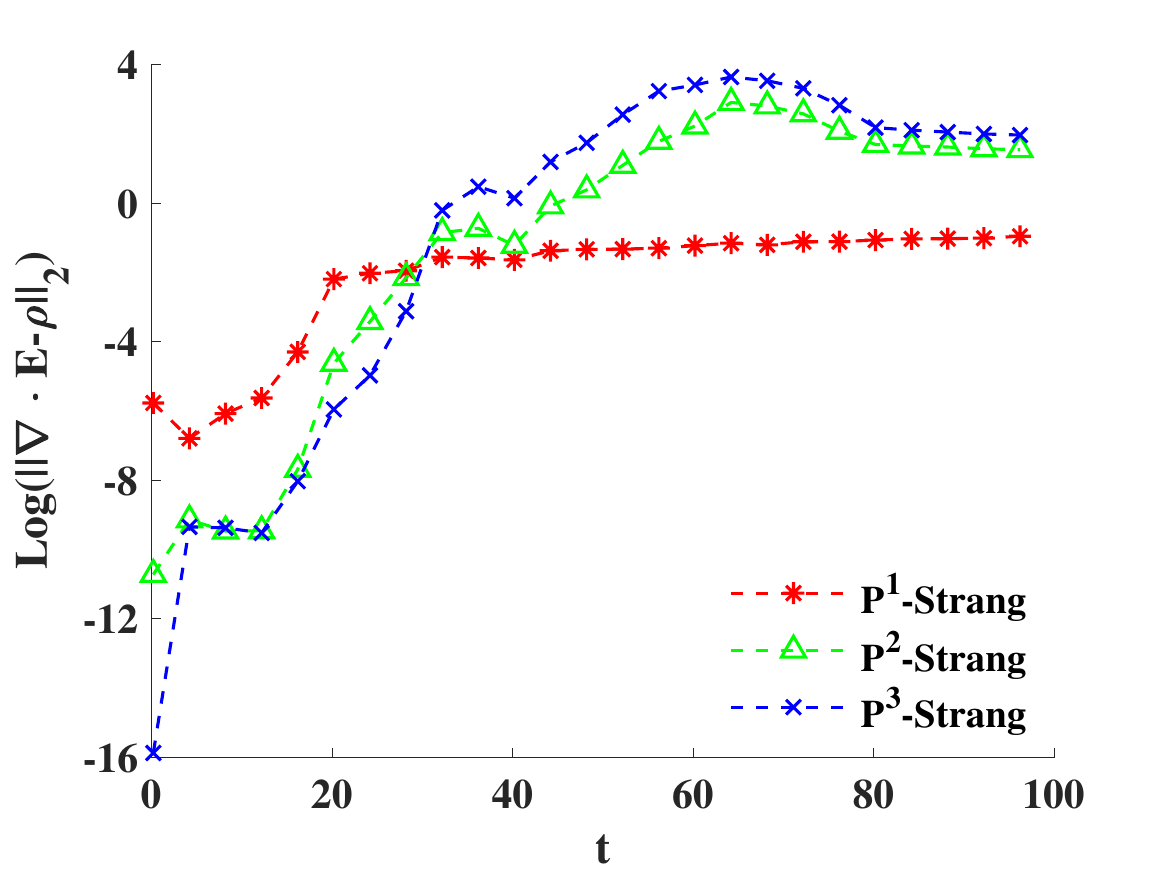}}
			\centerline{\textbf{(a)}}
		\end{minipage}
		\hfill
		\begin{minipage}{0.49\linewidth}
			\centerline{\includegraphics[width=1\linewidth]{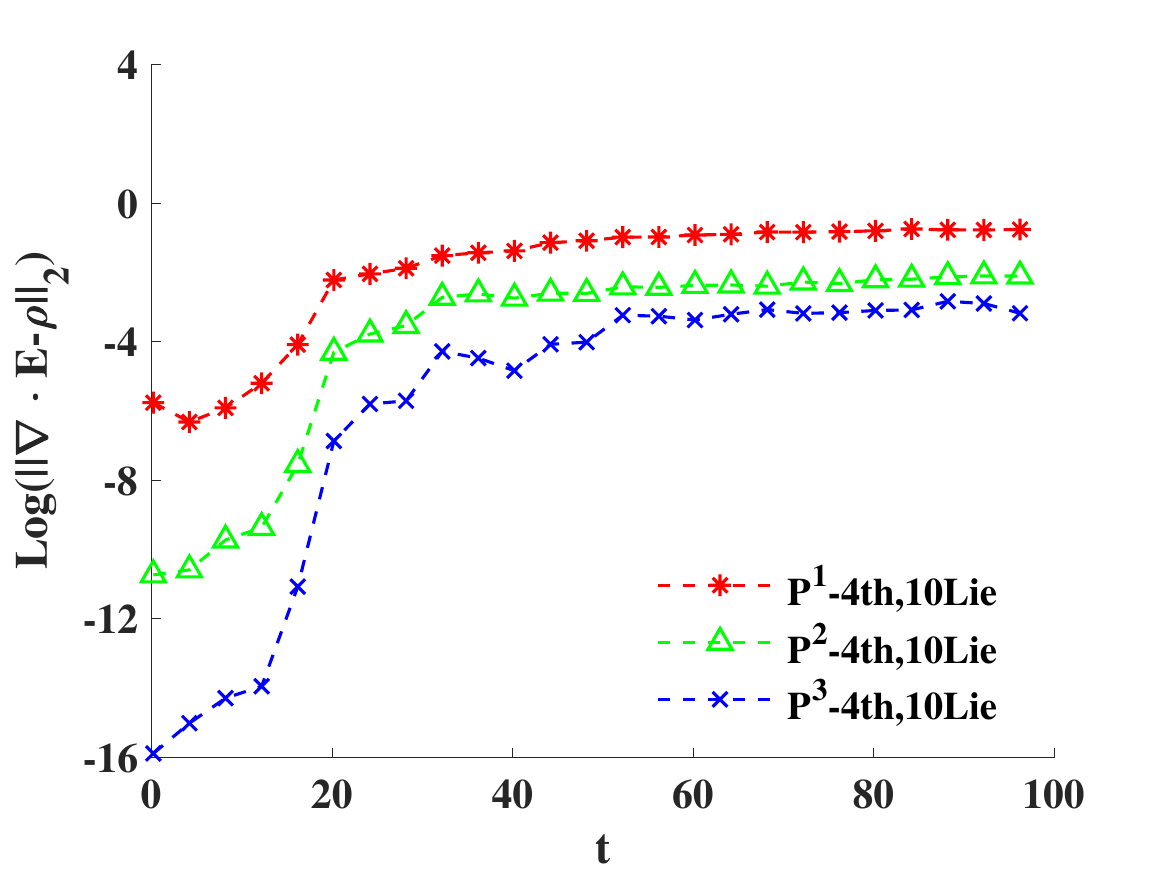}}
			\centerline{\textbf{(b)}}
		\end{minipage}
	\end{center}
	\caption{Two-stream instability II: $N_x=128,N_v=256$ and $\Delta t=0.2$. Time evolution of Gauss law's residuals under different spatial accuracies when using the ECSLDG method combined with the second-order Strang splitting scheme (a) and the fourth-order splitting scheme $10Lie$ (b). }
	\label{gauss-st3}
\end{figure}

\begin{figure}[htbp]
	\begin{center}
		\begin{minipage}{0.49\linewidth}
			\centerline{\includegraphics[width=1\linewidth]{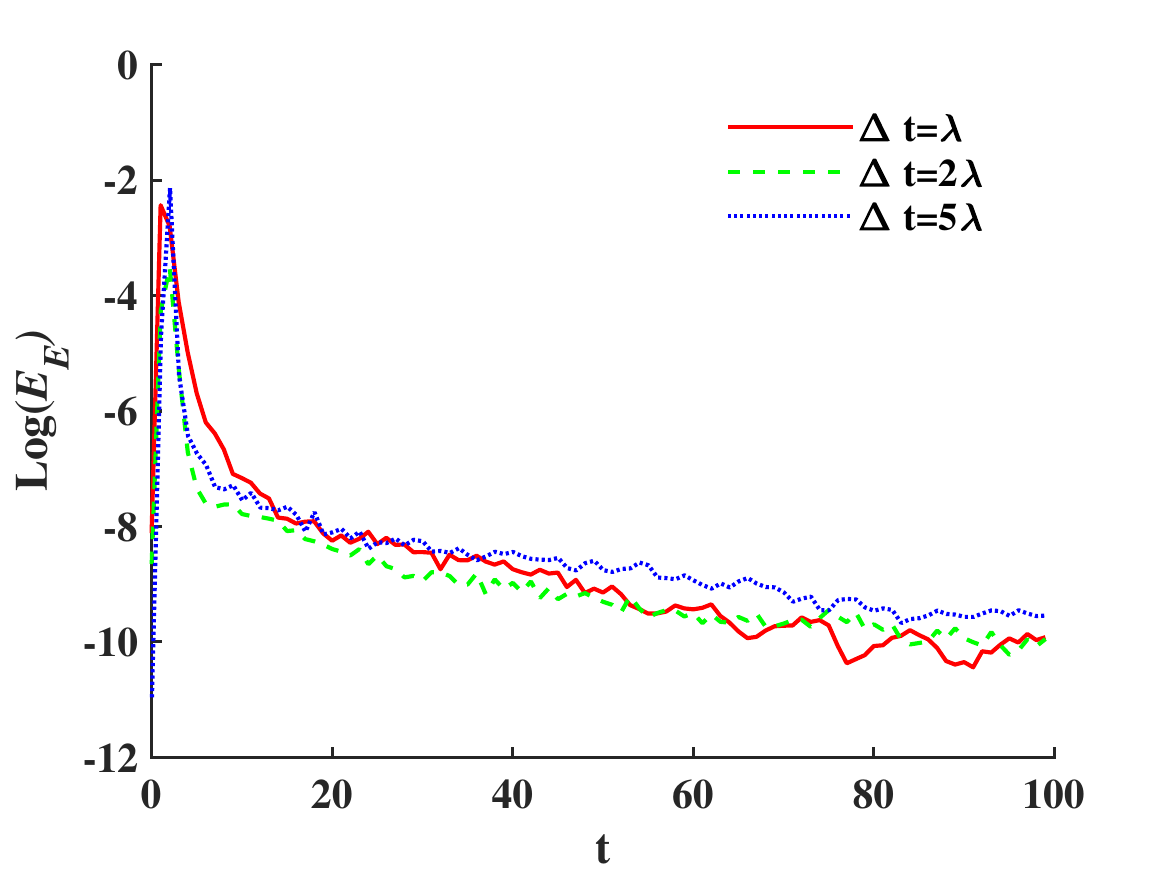}}
			\centerline{\textbf{(a)}}
		\end{minipage}
		\hfill
		\begin{minipage}{0.49\linewidth}
			\centerline{\includegraphics[width=1\linewidth]{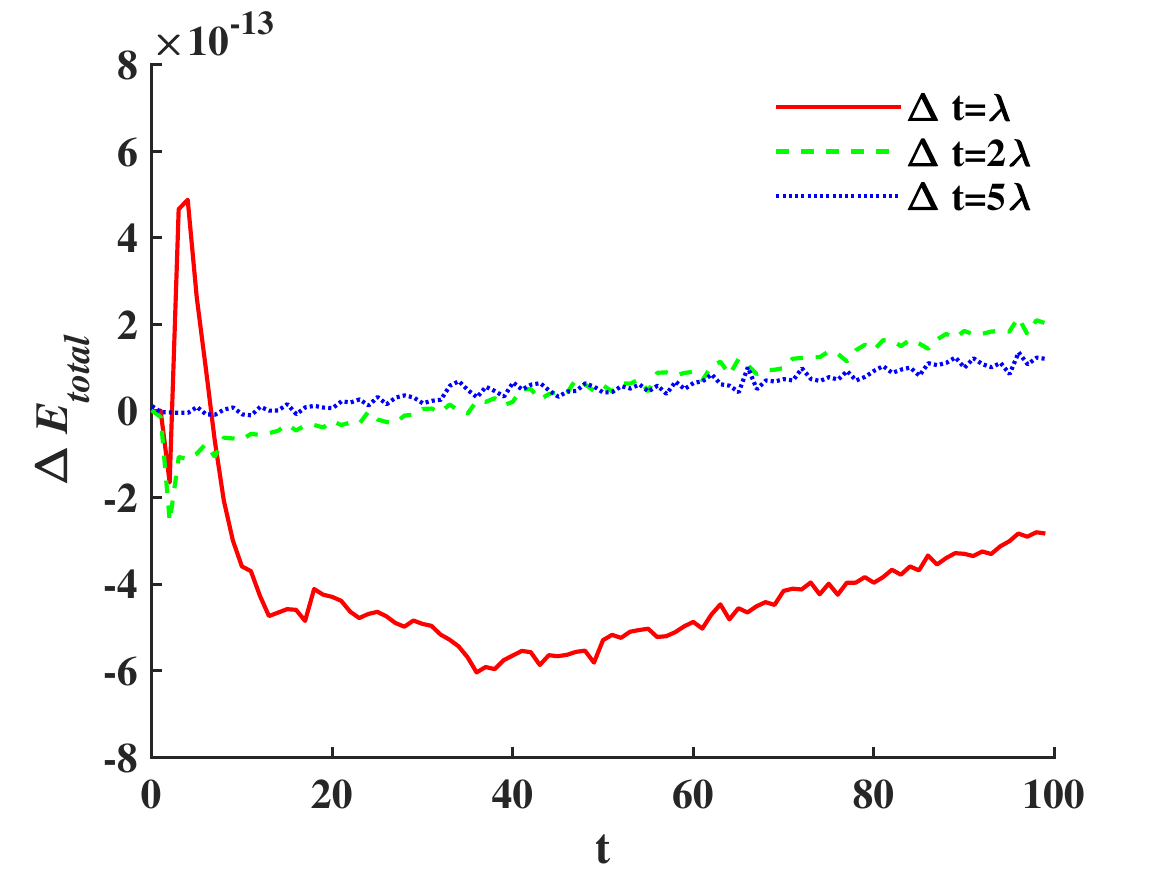}}
			\centerline{\textbf{(b)}}
		\end{minipage}
	\end{center}
	\caption{Two-stream instability II: $N_x=N_v=128$ and $\lambda = 0.01$.  Time evolution of the electric energy (a) and relative deviation of total energy (b) with different 
$\Delta t$. }
	\label{lambda-st3}
\end{figure}
To further investigate the cause of the unphysical oscillations observed when using the second-order Strang splitting scheme in the ECSLDG method, we examine the behavior of Gauss's law residuals under different splitting schemes. 
As shown in Fig. \ref{gauss-st3}(a), when using the second-order Strang splitting scheme, the Gauss residuals increase as the spatial resolution improves, leading to unphysical charge separation and potentially triggering unphysical oscillations observed in the phase space  plots (see Fig. \ref{plot-st3}(middle)). Encouragingly, Fig.~\ref{gauss-st3}(b) shows that when the fourth-order Lie splitting scheme ($10Lie$) is employed, the Gauss residuals decrease with improved spatial accuracy and remain low and stable over time.

Finally, we investigate the robustness of the proposed ECSLDG method when the temporal resolution is insufficient to resolve the plasma oscillation period. Traditional kinetic methods typically require fine temporal resolutions to resolve the plasma period due to stability constraints. As shown in Fig.~\ref{lambda-st3}(a), even when $\Delta t \ge \lambda$, the electric field behavior remains in good agreement with the fully resolved case. Moreover, the ECSLDG method conserves the total energy across all simulations, making it particularly advantageous for multiscale simulations.

The above tests demonstrate that the proposed ECSLDG method conserves total energy independently of spatial resolution. With high-order accuracy in both space and time, it provides more accurate and physically consistent descriptions of electrostatic phenomena. Importantly, the ECSLDG method does not require resolving the plasma oscillation period, highlighting its robustness and potential for efficient multiscale plasma simulations.

\section{Conclusion} \label{section:conclusion}


We have developed a high-order energy-conserving semi-Lagrangian discontinuous Galerkin (ECSLDG) method for the Vlasov-Ampère system. The proposed method combines the efficiency of explicit schemes with the energy conservation and unconditional stability properties of implicit schemes, without relying on nonlinear iteration. 
By integrating the high-order spatial accuracy of the discontinuous Galerkin approach with high-order temporal accuracy using suitable operator splitting techniques, the ECSLDG method achieves high-order accuracy in both space and time. Numerical simulations demonstrate that the high-order splitting method significantly improves Gauss's law enforcement, providing more accurate results compared to second-order schemes, especially when using a large CFL number.

While 1D1V simulations are presented in this paper, the method can be easily extended to multidimensional simulations using the operator splitting approach. Future work will focus on further exploring the scheme’s ability to conserve Gauss’s law and extending the method to the full Vlasov-Maxwell system.

\section*{Acknowledgments} 
This work was supported by the National Natural Science Foundation of China (No.~12201052), the Guangdong Provincial Key Laboratory of Inter-disciplinary Research and Application for Data Science of project code 2022B1212010006, BNBU Research Grant with No.~of UICR0700035-22 at Beijing Normal-Hong Kong Baptist University,Zhuhai, PR China, Guangdong basic and applied basic research foundation[2025A1515012182], National Key Laboratory for Computaional Physics [6142A05230201]. 
Liu is supported by a PostDoctoral Fellowship 1252224N from Research Foundation-Flanders (FWO).
Zheng was partially supported by the National Natural Science Foundation of China (No.~12471406), Science and Technology Commission of Shanghai Municipality (No.~22DZ2229014).
\clearpage
\bibliographystyle{unsrt} 
\bibliography{references}
\end{document}